\documentclass[10pt,a4paper]{article}

\usepackage{amsmath,amssymb,amsfonts,amsthm}
\usepackage{mathtools}
\usepackage{booktabs}
\usepackage{tikz}
\usetikzlibrary{arrows.meta,calc,decorations.pathreplacing}
\usepackage[hidelinks]{hyperref}
\usepackage{geometry}
\usepackage[expansion=false]{microtype}
\usepackage{cite}

\theoremstyle{plain}
\newtheorem{theorem}{Theorem}[section]
\newtheorem{proposition}[theorem]{Proposition}
\newtheorem{lemma}[theorem]{Lemma}
\newtheorem{corollary}[theorem]{Corollary}

\theoremstyle{definition}
\newtheorem{definition}[theorem]{Definition}

\theoremstyle{remark}
\newtheorem{remark}[theorem]{Remark}

\newtheorem{hypinner}{Hypothesis}

\newcommand{\ppart}[1]{\smallskip\noindent\textit{#1}\ }

\allowdisplaybreaks
\newcommand{\RR}{\mathbb{R}}
\newcommand{\ZZ}{\mathbb{Z}}
\newcommand{\EE}{\mathbb{E}}
\newcommand{\PP}{\mathbb{P}}
\newcommand{\Sph}[1]{S^{#1}}
\DeclareMathOperator{\rank}{rank}
\DeclareMathOperator{\im}{im}
\DeclareMathOperator{\spn}{span}

\DeclareMathOperator{\tr}{tr}
\DeclareMathOperator{\diag}{diag}
\DeclareMathOperator{\sgn}{sgn}

\DeclareMathOperator{\Var}{Var}
\DeclareMathOperator{\Cov}{Cov}
\newcommand{\transp}{^{\!\top}}

\newcommand{\Wish}{\mathcal{W}}
\newcommand{\Norm}{\mathcal{N}}
\newcommand{\Haus}[1]{\mathcal{H}^{#1}}       
\newcommand{\ip}[2]{\langle #1,\, #2\rangle}
\newcommand{\norm}[2][]{{#1\lVert #2 #1\rVert}}
\newcommand{\abs}[2][]{{#1\lvert #2 #1\rvert}}
\newcommand{\dd}{\mathop{}\!d}

\newcommand{\dimin}{D}                  
\newcommand{\ngen}{H}                   

\newcommand{\arr}{\mathcal{A}_{\ngen}}           

\newcommand{\Regge}{\mu_{\mathrm{Regge}}}

\newcommand{\sigw}{\sigma_w}                     
\newcommand{\wallis}[1]{\mathsf{I}_{#1}}

\newcommand{\Zlim}{Z_\infty}                     
\usepackage{enumitem}
\setlist[enumerate]{topsep=2pt,itemsep=1pt,parsep=0pt}
\setlist[itemize]{topsep=2pt,itemsep=1pt,parsep=0pt}
\AtBeginDocument{%
  \setlength{\abovedisplayskip}{4.5pt plus 1.5pt minus 1.5pt}%
  \setlength{\belowdisplayskip}{4.5pt plus 1.5pt minus 1.5pt}%
  \setlength{\abovedisplayshortskip}{2pt plus 1pt}%
  \setlength{\belowdisplayshortskip}{3pt plus 1pt minus 1pt}%
}

\newcommand{\inner}[2]{\ip{#1}{#2}}
\newcommand{\Unif}{\mathrm{Unif}}
\newcommand{\CC}{\mathbb{C}}
\newcommand{\Bin}{\mathrm{Bin}}
\newcommand{\indep}{\perp\!\!\!\perp}

\title{\textbf{Antipodal Parity and Second-Chaos Regge Curvature on Random Hyperplane Arrangements}}
\author{William Nelson\\[2pt] Google DeepMind\\[2pt] \texttt{wtnelson@google.com}}
\date{}

\begin{document}
\maketitle

\begin{abstract}
Limit theorems are proved for the signed Regge curvature of independent
oriented hyperplanes in Euclidean space under the piecewise-flat metric induced
by unit linear ramps on their positive half-spaces. Provided each hyperplane
law is invariant under orientation reversal, the curvature has zero mean and
zero first Hoeffding component at every finite sample size. In particular, it
is uncorrelated with all square-integrable unoriented functionals. For isotropic
Gaussian normals on bounded windows, the rescaled curvature is approximated in
mean square by a degenerate U-statistic and converges to a second-chaos series
diagonalised on odd spherical harmonics. In the plane, the global cone-angle
deficit admits a closed-form limit law universal up to scale under finite
fourth moments.
\end{abstract}

\section{Introduction}
\label{sec:intro}

\subsection{A signed curvature functional of a random hyperplane tessellation}
\label{subsec:intro-hook}

Whilst cell counts, skeleton volumes and intrinsic volumes of a random
hyperplane tessellation~\cite{SchneiderWeil2008,HugSchneider2024} depend on the
hyperplanes only as unoriented sets, one may ask how the asymptotic behaviour
changes for geometric functionals that also depend on the orientation of the
hyperplanes. What happens to the mean, the fluctuation scale and the
limiting law when reversing a single hyperplane reverses the local contribution
along its intersections?

To make this concrete, let $\Pi_1, \dots, \Pi_{\ngen}$ be $\ngen$ independent
random affine hyperplanes in $\RR^{\dimin}$, $\Pi_h = \{x \in \RR^{\dimin} :
\ip{w_h}{x} + b_h = 0\}$, with normal vectors $w_h \in \RR^{\dimin}$ and offsets
$b_h \in \RR$. Each generator $(w_h, b_h)$ defines open half-spaces $\Pi_h^\pm
:= \{x \in \RR^{\dimin} : \pm(\ip{w_h}{x} + b_h) > 0\}$. Whilst the unoriented
hyperplanes $\{\Pi_h\}$ determine the open cells $\tau$ of the tessellation, the
positive half-spaces $\Pi_h^+$ equip each cell with the constant Riemannian
metric
\begin{equation}\label{eq:intro-metric}
  G(\tau) := I_{\dimin} + \sum_{h \in A(\tau)} w_h w_h\transp,
  \qquad
  A(\tau) := \{1 \le h \le \ngen : \tau \subseteq \Pi_h^+\}.
\end{equation}
This metric is pulled back from the Euclidean metric on $\RR^{\dimin + \ngen}$
by the piecewise-affine ramp lift
\begin{equation}\label{eq:intro-map}
  \Phi_{W,b}\colon \RR^{\dimin} \to \RR^{\dimin + \ngen},
  \qquad
  \Phi_{W,b}(x) = \bigl(x,\, (w_1\transp x + b_1)_+,\, \dots,\, (w_{\ngen}\transp x + b_{\ngen})_+\bigr).
\end{equation}
In machine learning, $\Phi_{W,b}$ is the graph of a single rectified linear unit
(ReLU) layer, whose one-sided ramps $(w_h\transp x + b_h)_+$ break orientation
symmetry, unlike combinatorial region counts and Euclidean skeleton volumes of
ReLU arrangements (Section~\ref{subsec:intro-literature}). Because adjacent cell
metrics agree tangentially along their common wall, $(\RR^{\dimin}, G)$ is an
intrinsic polyhedral manifold in the sense of Regge~\cite{Regge1961} and
Cheeger, M\"uller and Schrader~\cite{CheegerMullerSchrader1984}: it is flat on
cells, has zero crease curvature across walls, and concentrates all intrinsic
curvature on the codimension-$2$ \emph{bones} $B_{cc'} = \Pi_c \cap \Pi_{c'}$.
Around a point $x \in B_{cc'}$ (for a pair $B = \{c, c'\}$) away from the
remaining hyperplanes, $\Pi_c^\pm$ and $\Pi_{c'}^\pm$ cut out four transverse
quadrants $\Pi_c^{\sigma_c} \cap \Pi_{c'}^{\sigma_{c'}}$ indexed by $\sigma =
(\sigma_c, \sigma_{c'}) \in \{\pm 1\}^2$, with cell metrics
\[
  G_\sigma(x)
  = G_B(x) + \mathbf{1}\{\sigma_c = +1\}\,w_c w_c\transp
           + \mathbf{1}\{\sigma_{c'} = +1\}\,w_{c'} w_{c'}\transp,
  \qquad
  G_B(x)
  := I_{\dimin} + \sum_{\substack{h \notin \{c, c'\} \\ x \in \Pi_h^+}} w_h w_h\transp.
\]
Writing $\theta_\sigma(x)$ for the $G_\sigma(x)$-dihedral angle of quadrant
$\sigma$ (Definition~\ref{def:arrangement}), the pointwise angular deficit and
the \emph{ambient Regge curvature measure} on a Borel window $\Omega \subseteq
\RR^{\dimin}$ are
\[
  \varepsilon_{cc'}(x)
  := 2\pi - \sum_{\sigma \in \{\pm 1\}^2} \theta_\sigma(x),
  \qquad
  \Regge(\Omega)
  := \sum_{1 \le c < c' \le \ngen}
     \int_{B_{cc'} \cap \Omega} \varepsilon_{cc'}(x)\, \dd\Haus{\dimin-2}(x).
\]
Because $G_B(x)$ sums over all active background half-spaces $\Pi_h^+$ ($h
\notin \{c, c'\}$), $\varepsilon_{cc'}(x)$ depends on all $\ngen$ generators
rather than on $(\Pi_c, \Pi_{c'})$ alone. In particular, $\Regge(\Omega)$ is not
itself a pairwise U-statistic.

For $\dimin \ge 3$ the bones are unbounded $(\dimin - 2)$-flats, so we restrict
to windows $\Omega$ of bounded $(\dimin-2)$-slice volume. In Regge
calculus~\cite{Regge1961} and Cheeger, M\"uller and
Schrader~\cite{CheegerMullerSchrader1984}, intrinsic polyhedral curvature weights
$\varepsilon_{cc'}(x)$ by the $G$-induced bone volume element
$\sqrt{\det(G_B(x)|_{TB_{cc'}})}\,\dd\Haus{\dimin-2}(x)$ (which coincides with
$\dd\Haus{0}$ in the plane). Integrating instead against ambient
$\dd\Haus{\dimin-2}$ in $\Regge(\Omega)$ serves two purposes when $\dimin \ge 3$.
Firstly, it isolates the transverse deficit $\varepsilon_{cc'}(x)$ from the
$(\ngen\sigw^2/2)^{(\dimin-2)/2}$ tangential metric stretching along $TB_{cc'}$.
Secondly, the reference measure on $B_{cc'}$ depends on $(\Pi_c, \Pi_{c'})$
alone, so bone-wise concentration integrates directly to a pairwise U-statistic.
Remark~\ref{rem:cms-intrinsic-equivalence} shows that our finite-$\ngen$ parity
results carry over verbatim to the intrinsically weighted
Cheeger--M\"uller--Schrader measure, whilst leaving its quantitative $L^2$
transfer to future work.

In the plane ($\dimin = 2$), the bones are single vertices $x_{cc'}$.
Thus the polyhedral Gauss--Bonnet
identity~\cite{CohnVossen1935,Banchoff1967} reduces the global sum over all
$\binom{\ngen}{2}$ vertices to a boundary deficit at infinity, $\Regge(\RR^2) =
\sum_{1 \le c < c' \le \ngen} \varepsilon_{cc'}(x_{cc'}) = 2\pi -
\Theta_\infty(W)$. Here the \emph{far-field cone angle} $\Theta_\infty(W) =
\sum_{k=1}^{2\ngen}\alpha_k^\infty(W)$ sums the intrinsic opening angles of the
$2\ngen$ unbounded wedges in the metric~\eqref{eq:intro-metric} (an
intrinsic-metric counterpart of Euclidean conical angle
sums~\cite{HugSchneider2016}) and does not depend on the offsets $b$.

\subsection{Orientation parity and the pairwise limit kernel}
\label{subsec:intro-idea}

Fix a bone $B_{cc'}$ ($B = \{c, c'\}$) and a point $x \in B_{cc'}$ off the
remaining hyperplanes. Reversing $(w_c, b_c) \mapsto (-w_c, -b_c)$ leaves the
unoriented tessellation, the bone $B_{cc'}$ and that bone's background Gram
matrix $G_B(x)$ unchanged (since $G_B(x)$ sums only over $h \notin \{c, c'\}$),
whilst swapping $\Pi_c^+ \leftrightarrow \Pi_c^-$ shifts $+w_c w_c\transp$ to
the supplementary quadrants across $\Pi_c$ and negates $\varepsilon_{cc'}(x)$
pointwise (Lemma~\ref{lem:sign-flip-oddness}). On bones not containing $\Pi_c$,
however, flipping $(w_c, b_c)$ changes $\mathbf{1}\{x \in \Pi_c^+\}$ inside
their background metric. For independent generators $\xi_h := (w_h, b_h)$,
recall the orthogonal Hoeffding (ANOVA) decomposition~\cite{Hoeffding1948} $X =
\EE[X] + \sum_{k=1}^{\ngen} P_k(X)$ into $k$-generator projections
(Definition~\ref{def:bone-conditioning}), with $P_1(X) =
\sum_{j=1}^{\ngen}(\EE[X \mid \xi_j] - \EE[X])$. Provided $\xi_h \stackrel{d}{=}
-\xi_h$, for every bone $B_{cc'}$ and any single conditioned generator $\xi_j$
one can flip a bone generator in $\{c, c'\} \setminus \{j\}$ without altering
$\xi_j$, negating $Y_{cc'}(\Omega) := \int_{B_{cc'} \cap \Omega}
\varepsilon_{cc'}(x)\,\dd\Haus{\dimin-2}(x)$. Thus $\EE[Y_{cc'}(\Omega) \mid
\xi_j] = 0$, so both $\EE[\Regge(\Omega)]$ and $P_1(\Regge(\Omega))$ vanish
identically at every finite $\ngen$. The leading fluctuations are pushed into
$P_2(\Regge(\Omega))$ at scale $\ngen^{-1}$, with non-Gaussian second-chaos
limits~\cite{RubinVitale1980,Janson1997}.

To identify the limit kernel despite the background coupling in $G_B(x)$,
Section~\ref{sec:local-geometry} applies a $2 \times 2$ Woodbury reduction to
express the four angles $\theta_\sigma(x)$, and hence an exact four-arccos
formula for $\varepsilon_{cc'}(x)$ (Theorem~\ref{thm:four-arccos}), solely
through the overlap matrix $M_B(x) := L_B\transp G_B(x)^{-1} L_B$ ($L_B := [w_c
\;\; w_{c'}]$). Expanding around $M_B = 0$ cancels zeroth- and first-order terms
across the four quadrants, leaving the quadratic proxy
(Proposition~\ref{prop:deficit-integral})
\[
  \varepsilon_{cc'}(x)
  = -\frac{1}{4}\,M_{B,12}(x)\sqrt{\det M_B(x)}\,
    \bigl(1 + O(\tr M_B(x))\bigr).
\]
Conditionally on $(\xi_c, \xi_{c'})$ at any $x \in B_{cc'} \setminus \{0\}$,
central symmetry $\xi_h \stackrel{d}{=} -\xi_h$ ($h \notin \{c, c'\}$) flips
$\{x \in \Pi_h^+\}$ whilst preserving $w_h w_h\transp$, so the indicators
$\mathbf{1}\{x \in \Pi_h^+\}$ are i.i.d.\ $\mathrm{Bernoulli}(1/2)$ variables
independent of $w_h w_h\transp$ (Lemma~\ref{lem:active-count}). Under isotropic
Gaussian normals $\Norm(0, \sigw^2 I_{\dimin})$, $G_B(x) - I_{\dimin}$ is a
$\Bin(\ngen - 2, \frac{1}{2})$-thinned Wishart matrix concentrating around
$\frac{1}{2}\ngen\sigw^2 I_{\dimin}$. Thus $M_B(x) \approx
\frac{2}{\ngen\sigw^2}L_B\transp L_B$ has scale $\Theta(\ngen^{-1})$, each bone
deficit $\varepsilon_{cc'}(x)$ has scale $\Theta(\ngen^{-2})$, and the rescaled
deficit $\frac{\ngen^2}{4}\varepsilon_{cc'}(x)$ concentrates around the pairwise
kernel
\begin{equation}\label{eq:intro-Z-infinity}
  \Zlim(w_c, w_{c'})
  := -\frac{\ip{w_c}{w_{c'}}}{4\sigw^4}
     \sqrt{\norm{w_c}^2 \norm{w_{c'}}^2 - \ip{w_c}{w_{c'}}^2}
   = -\frac{\ip{w_c}{w_{c'}}\,\norm{w_c \wedge w_{c'}}}{4\sigw^4}.
\end{equation}
Note that single-bone concentration alone would not decouple the
$\binom{\ngen}{2}$ bones. In an unsigned sum, $O(\ngen^{-3})$ background
perturbations per bone add coherently across bones to produce
$\Theta(\ngen^{-1/2})$ first-order fluctuations. Here orientation oddness
annihilates all projections $\pi_A Y_{cc'}(\Omega)$ with $\{c, c'\}
\not\subseteq A$ (Lemma~\ref{lem:hoeffding-support}). Thus $P_2(\Regge(\Omega))$
is a degenerate order-$2$ U-statistic with an $\ngen$-dependent kernel
approximating $\frac{4}{\ngen^2}\Zlim$, whilst disjoint bones in $\sum_{k \ge 3}
P_k(\Regge(\Omega))$ share only projections of order $\abs{A} \ge 4$. Bounding
single- and double-generator replacement differences via a two-step Efron--Stein
inequality~\cite{EfronStein1981,BLM2013} (accommodating indicator jumps across
$\Pi_h$ and general symmetric offsets $b_h$) shows that
$\frac{\ngen}{4}\sum_{k \ge 3} P_k(\Regge(\Omega))$ vanishes in $L^2(\PP)$ at
rate $O(\ngen^{-1/2})$, reducing the limit law of $\frac{\ngen}{4}\Regge(\Omega)$
to a degenerate order-$2$ U-statistic with kernel $\Zlim$.

\subsection{Main results}
\label{subsec:intro-main-results}

In this paper we establish four main results, working at fixed dimension $\dimin
\ge 2$ as $\ngen \to \infty$ (the precise standing models are collected in
Definition~\ref{def:standing-assumptions}). They concern, in turn, parity
cancellation at finite $\ngen$ (Section~\ref{sec:parity-kernel}), the
second-chaos transfer theorem (Section~\ref{sec:transfer}), the planar limit
law (Section~\ref{sec:planar}) and the all-dimensional spectrum and Kolmogorov
rate (Section~\ref{sec:high-dim}).

For a Borel window $\Omega \subseteq
\RR^{\dimin}$, write $\ell_{cc'}(\Omega) := \Haus{\dimin-2}(B_{cc'} \cap
\Omega)$, let $\beta_\Omega := \sup_F \Haus{\dimin-2}(F \cap \Omega)$ over
affine $(\dimin-2)$-flats, and set $U_{\ngen,\infty}(\Omega) :=
\frac{1}{\ngen}\sum_{c<c'}\Zlim(w_c, w_{c'})\ell_{cc'}(\Omega)$ (abbreviated
$U_{\ngen,\infty}(W)$ when $\ell_{cc'} \equiv 1$).

Firstly, provided each generator satisfies joint central symmetry $(w_h, b_h)
\stackrel{d}{=} (-w_h, -b_h)$, we have $\EE[\Regge(\Omega)] = P_1(\Regge(\Omega))
= 0$ at every finite $\ngen$. In fact, $\Regge(\Omega)$ is uncorrelated with every
square-integrable functional of the unoriented tessellation
(Propositions~\ref{thm:first-chaos-annihilation}
and~\ref{thm:parity-orthogonality}). In particular, in the plane
$\EE[\Theta_\infty(W)] = 2\pi$ without moment assumptions
(Corollary~\ref{cor:planar-moments}). The limiting objects inherit this
orthogonality: $U_{\ngen,\infty}(W)$ under isotropic elliptically symmetric
normals, and $\frac{\ngen}{4}\Regge(\Omega)$ under Gaussian normals on any
window with $\beta_\Omega < \infty$, are asymptotically independent of every
converging sequence of unoriented functionals
(Corollary~\ref{cor:parity-asymp-indep}).

Secondly, for isotropic Gaussian normals $\Norm(0, \sigw^2 I_{\dimin})$ and
symmetric i.i.d.\ offsets on any window with $\beta_\Omega < \infty$, one can
compare $\frac{\ngen}{4}\Regge(\Omega) =
\frac{1}{\ngen}\sum_{c<c'}\bigl(\frac{\ngen^2}{4}Y_{cc'}(\Omega)\bigr)$ directly
with $U_{\ngen,\infty}(\Omega) =
\frac{1}{\ngen}\sum_{c<c'}\Zlim(w_c, w_{c'})\ell_{cc'}(\Omega)$.
Theorem~\ref{thm:second-chaos-transfer} gives the mean-square bound
$\norm{\frac{\ngen}{4}\Regge(\Omega) - U_{\ngen,\infty}(\Omega)}_{L^2(\PP)} \le
C_{\dimin,\sigw}\,\beta_\Omega\,\ngen^{-1/2}$. This $O(\ngen^{-1/2})$ rate is
sharp, since single-hyperplane perturbations of $G_B(x)$ produce a non-zero
third Hoeffding component satisfying
$\norm{\frac{\ngen}{4}P_3(\Regge(\Omega))}_{L^2(\PP)} = \Theta(\ngen^{-1/2})$
whenever $\PP(\ell_{12}(\Omega) > 0) > 0$. Hilbert--Schmidt spectral expansion
of the completely degenerate U-statistic
$U_{\ngen,\infty}(\Omega)$~\cite{RubinVitale1980,Serfling1980}
(Proposition~\ref{prop:degenerate-ustat-limit}) then yields
$\frac{\ngen}{4}\Regge(\Omega) \xrightarrow{d} \frac{1}{2}\sum_{k \ge
1}\mu_{k,\Omega}(\gamma_k^2 - 1)$ with $\gamma_k \stackrel{\mathrm{i.i.d.}}{\sim}
\Norm(0,1)$, where $(\mu_{k,\Omega})_{k \ge 1}$ are the non-zero eigenvalues of
$(T_{K_\Omega}f)(\xi) := \EE_{\xi'}[K_\Omega(\xi, \xi')f(\xi')]$ for the
slice-weighted kernel $K_\Omega(\xi_c, \xi_{c'}) :=
\Zlim(w_c, w_{c'})\ell_{cc'}(\Omega)$.

Thirdly, in the plane ($\dimin = 2$, $\Omega = \RR^2$), $T_{\Zlim}$ has
odd-harmonic eigenvalues $\lambda_m = 4/(\pi(2m-1)(2m+3))$ ($m \ge 0$,
multiplicity $2$). Thus under Gaussian normals $\frac{\ngen}{4}(2\pi -
\Theta_\infty(W)) \xrightarrow{d} \mathcal{L}_2 := \sum_{m \ge 0}\lambda_m E_m$
($E_m \stackrel{\mathrm{i.i.d.}}{\sim} \mathrm{Exp}(1)$). This limit law has variance
$\frac{1}{4}$, skewness $-3/\pi$, exact left tail $\frac{4}{3\pi}e^{-3\pi x/4}$
($x \ge 0$) and Kolmogorov rate $O(\ngen^{-1/3})$
(Theorem~\ref{thm:planar-limit-law} and Corollary~\ref{cor:planar-kolmogorov}).
Its moment generating function evaluates in closed form to
\[
  M_{\mathcal{L}_2}(s)
  = \prod_{m \ge 0} \frac{1}{1 - \lambda_m s}
  = -\sec\!\left(\sqrt{\pi^2 + \pi s}\right),
  \qquad s \in \bigl(-\tfrac{3\pi}{4}, \tfrac{5\pi}{4}\bigr).
\]
An offset-free circle-integral representation then shows that central symmetry
and $\EE[\norm{w_1}^4] < \infty$ alone imply $\frac{\ngen}{4}(2\pi -
\Theta_\infty(W)) - U_{\ngen,\infty}(\widetilde{W}) \xrightarrow{\PP} 0$,
yielding the universal limit $\frac{\tau_2(P_w)}{2}\mathcal{L}_2$ across
elliptically symmetric normal laws (Theorem~\ref{thm:planar-universality}).

Finally, in all dimensions $\dimin \ge 2$, Theorem~\ref{thm:funk-hecke}
diagonalises $T_{\Zlim}$ on odd spherical harmonics with closed-form eigenvalues
$\kappa_{2m+1}(\dimin)$, zero trace and signature $(\dimin, \infty)$. This gives an
explicit second-chaos limit for central arrangements ($b_h \equiv 0$) on balls
$\Omega = B(0, R)$. Because $\Zlim(w, w) \equiv 0$ pointwise,
$U_{\ngen,\infty}(W)$ is exactly a trace-zero quadratic form in a normalised
sum of $\ell^2$-valued vectors, with no diagonal correction. The
infinite-dimensional quadratic-form bound of G\"otze and
Zaitsev~\cite{GoetzeZaitsev2014} then gives $d_{\mathrm{Kol}}(U_{\ngen,\infty}(W),
U_\infty) \le C_{\dimin}\,\EE[\norm{w_1}^8]\,(\EE[\norm{w_1}^4])^{-2}\,\ngen^{-1}$
whenever $\EE[\norm{w_1}^8] < \infty$ (Theorem~\ref{thm:hilbert-kolmogorov}).

Taken together, these results answer the questions posed at the outset. Unsigned
tessellation functionals have non-zero means (as does the unsigned total
variation $\abs{\Regge}(\Omega)$ of Section~\ref{sec:remarks}, at scale
$\Theta(1)$) and typically obey Gaussian central limit theorems at relative
scale $\ngen^{-1/2}$, driven by $P_1$. Here orientation parity forces
$\EE[\Regge(\Omega)] = 0$ at every finite $\ngen$ (so $\EE[\Theta_\infty(W)] =
2\pi$ in the plane) and pushes fluctuations down to scale $\ngen^{-1}$, with
skewed, trace-zero second-chaos laws of signature $(\dimin, \infty)$. In
addition, although defined on the same $\ngen$ hyperplanes, $\Regge(\Omega)$ is
uncorrelated at every finite $\ngen$ and asymptotically independent as $\ngen
\to \infty$ of every converging unoriented observable.

\subsection{Relation to the literature}
\label{subsec:intro-literature}

For unsigned functionals of stationary Poisson line and hyperplane
tessellations, distribution theory, variances and central limit theorems go back
to Miles~\cite{Miles1964}, Paroux~\cite{Paroux1998} and Heinrich, Schmidt and
Schmidt~\cite{HeinrichSchmidtSchmidt2006}. These are developed systematically
in~\cite{SchneiderWeil2008,HugSchneider2024}. For the hyperbolic case, where the
Gaussian limit can fail, see Herold, Hug and
Th\"ale~\cite{HeroldHugThale2021}. Quantitative normal approximations for
U-statistics, intersection processes and binomial point processes, via
Malliavin--Stein and Chatterjee-type methods, can be found
in~\cite{ReitznerSchulte2013,LastPenroseSchulteThaele2014,LastPenrose2017,%
LachiezeReyPeccati2017}. Euclidean angle sums and face counts of
random conical tessellations and Poisson polyhedra are treated
in~\cite{HugSchneider2016,KabluchkoThaeleZaporozhets2020}. Here we
work directly with a finite, non-Poissonised arrangement of $\ngen$ hyperplanes.
The base identity $I_{\dimin}$ in~\eqref{eq:intro-metric} keeps cells with
fewer than $\dimin$ active hyperplanes uniformly non-degenerate. Each bone
deficit also couples the two hyperplanes meeting at the bone to all remaining
active hyperplanes through the inverse background Gram matrix. The key
difference, however, is that orientation oddness eliminates the first Hoeffding
component, so Gaussian central limits give way to non-Gaussian second-chaos
laws. For fixed-kernel degenerate U-statistics, spectral $\chi_1^2$-series
limits go back to Gregory~\cite{Gregory1977}, Rubin and
Vitale~\cite{RubinVitale1980}, Serfling~\cite{Serfling1980} and Dynkin and
Mandelbaum~\cite{DynkinMandelbaum1983}; Nourdin, Peccati and
Reinert~\cite{NourdinPeccatiReinert2010} establish fourth-moment universality of
Gaussian Wiener chaos for homogeneous sums, and Peccati and
Th\"ale~\cite{PeccatiThaele2013} treat non-Gaussian second-chaos limits on
Poisson spaces. Our setting is therefore quite different from the
vanishing-fourth-cumulant Gaussian regime of de Jong~\cite{deJong1987},
Lachi\`eze-Rey and Peccati~\cite{LachiezeReyPeccati2013} and D\"obler and
Peccati~\cite{DoblerPeccati2017}. Our
$O(\ngen^{-1})$ Kolmogorov rate in Theorem~\ref{thm:hilbert-kolmogorov} instead
comes from the Hilbert-space quadratic-form bound of G\"otze and
Zaitsev~\cite{GoetzeZaitsev2014} (extending~\cite{BentkusGoetze1999}), which
applies without diagonal drift because the kernel vanishes on the diagonal.

On the spectral side, diagonalising homogeneous spherical kernels via the
Funk--Hecke theorem~\cite{DaiXu2013} connects $T_{\Zlim}$ to the zonal
arc-cosine and ReLU kernels studied in machine learning
(see~\cite{ChoSaul2009,Bach2017} and the references therein). However, unlike those positive-definite
kernels, $\Zlim$ is odd in each argument and has trace zero. It has signature
$(\dimin, \infty)$ in every dimension $\dimin \ge 2$.

Finally, the ramp map~\eqref{eq:intro-map} is the graph of a single randomly
initialised ReLU layer. Combinatorial region counts and Euclidean skeleton or
intrinsic-volume functionals of ReLU arrangements are invariant under
orientation reversal (see~\cite{HaninRolnick2019,GrigsbyLindsey2022} and the
references therein). In contrast, each bone contribution $Y_{cc'}(\Omega)$
changes sign when either of its two hyperplanes is reversed, making
$\Regge(\Omega)$ orthogonal to all such unoriented observables at every finite
$\ngen$. The deterministic polyhedral rigidity and Gauss--Bonnet theory
of~\eqref{eq:intro-map} is developed in the companion
paper~\cite{NelsonPaperA2026}. The present paper is self-contained.

\medskip
The rest of the paper is organised as follows. Section~\ref{sec:model} sets up
the polyhedral model and standing assumptions.
Section~\ref{sec:local-geometry} then proves the deterministic four-arccos formula
and operator-norm perturbation bounds. Section~\ref{sec:parity-kernel} derives
sign-flip annihilation, the conditional thinned-Wishart law and the limit kernel
$\Zlim$. This is followed in Section~\ref{sec:transfer} by the cross-bone Efron--Stein
bound and transfer theorem. Sections~\ref{sec:planar} and~\ref{sec:high-dim}
establish the planar limit theory and all-dimensional Funk--Hecke and Kolmogorov
results. Finally, Section~\ref{sec:remarks} discusses implications and open
problems. Appendices~\ref{sec:app-wishart}--\ref{sec:app-spectral}
collect the thinned-Wishart remainder, third-Hoeffding sharpness and spectral
proofs.
\section{Model and standing assumptions}
\label{sec:model}

\subsection{Polyhedral graph immersion and local bone geometry}
\label{subsec:model-geometry}

Fix $\dimin \ge 2$ and $\ngen \ge 2$. Let
$W = (w_1, \dots, w_{\ngen})\transp \in \RR^{\ngen \times \dimin}$ and
$b = (b_1, \dots, b_{\ngen})\transp \in \RR^{\ngen}$ collect $\ngen$ affine
hyperplane generators $(w_h, b_h) \in \RR^{\dimin} \times \RR$. We equip
$\RR^{\dimin}$ with the standard inner product $\ip{u}{v} := \sum_{i=1}^{\dimin}
u_i v_i$, norm $\norm{u} := \sqrt{\ip{u}{u}}$, identity $I_{\dimin}$ and wedge
norm $\norm{u \wedge v} := \sqrt{\norm{u}^2\norm{v}^2 - \ip{u}{v}^2}$. With
$(z)_+ := \max(z, 0)$ applied componentwise, the ramp lift $\phi_{W,b}(x) :=
(W x + b)_+$ and graph immersion $\Phi_{W,b}(x) := (x, \phi_{W,b}(x)) \in
\RR^{\dimin} \times \RR^{\ngen}$ embed $\mathcal{M}_{W,b} :=
\Phi_{W,b}(\RR^{\dimin})$ as a Lipschitz polyhedral $\dimin$-submanifold over
the base metric $I_{\dimin}$.

Each index $h \in [\ngen] := \{1, \dots, \ngen\}$ with $w_h \ne 0$ defines an
oriented wall $\Pi_h := \{x \in \RR^{\dimin} : \ip{w_h}{x} + b_h = 0\}$ and open
half-spaces $\Pi_h^\pm := \{x \in \RR^{\dimin} : \pm(\ip{w_h}{x} + b_h) > 0\}$.
The arrangement $\arr := \{\Pi_1, \dots, \Pi_{\ngen}\}$ partitions
$\RR^{\dimin} \setminus \bigcup_{h=1}^{\ngen} \Pi_h$ into open convex polyhedral
cells $\tau \in \mathcal{C}(\arr)$. Each cell carries an active set $A(\tau) := \{h \in
[\ngen] : \ip{w_h}{x} + b_h > 0 \text{ on } \tau\}$ and a constant pullback metric
$G(\tau) = D\Phi_{W,b}(\tau)\transp D\Phi_{W,b}(\tau) = I_{\dimin} + \sum_{h \in
A(\tau)} w_h w_h\transp \in \mathrm{Sym}_{\dimin}^{++}$. For each pair $B = \{c,
c'\} \subset [\ngen]$ ($1 \le c < c' \le \ngen$) with $\norm{w_c \wedge w_{c'}}
> 0$, the codimension-$2$ bone is $B_{cc'} := \Pi_c
\cap \Pi_{c'}$. Its dimension is $\dimin - 2$, reducing to a single vertex $\{x_{cc'}\}$ when
$\dimin = 2$. The transverse interior of the bone is $B_{cc'}^{\circ} := B_{cc'} \setminus \bigcup_{h \notin \{c,
c'\}} \Pi_h$.

At a transverse point $x_0 \in B_{cc'}^{\circ}$, write $I_h(x) :=
\mathbf{1}\{\ip{w_h}{x} + b_h > 0\}$, background active set $A(x_0) := \{h \in
[\ngen] \setminus \{c, c'\} : I_h(x_0) = 1\}$ and active count $m_B(x_0) :=
\abs{A(x_0)}$. The \emph{background Gram matrix} and the four incident quadrant
metrics ($\sigma = (\sigma_c, \sigma_{c'}) \in \{\pm 1\}^2$) are
\begin{align*}
  G_B(x_0) &:= I_{\dimin} + \sum_{h \in A(x_0)} w_h w_h\transp \in \mathrm{Sym}_{\dimin}^{++}, \\
  G_\sigma(x_0) &:= G_B(x_0) + \mathbf{1}\{\sigma_c = +1\}\, w_c w_c\transp + \mathbf{1}\{\sigma_{c'} = +1\}\, w_{c'} w_{c'}\transp.
\end{align*}
Setting $L_B := [w_c \;\; w_{c'}] \in \RR^{\dimin \times 2}$, the $2 \times 2$
\emph{dual overlap matrix} and \emph{normalised Schur invariants} at $x_0$ are
\begin{equation}\label{eq:model-M-Gamma}
  \begin{gathered}
    M_B(x_0) := L_B\transp G_B(x_0)^{-1} L_B \in \mathrm{Sym}_2^{++}, \\
    a_1 := M_{B,11} > 0, \qquad a_2 := M_{B,22} > 0, \qquad
    \rho := \frac{M_{B,12}}{\sqrt{a_1 a_2}} \in (-1, 1).
  \end{gathered}
\end{equation}

\begin{definition}
\label{def:arrangement}
At $x_0 \in B_{cc'}^{\circ}$, the outward $G_\sigma(x_0)$-covariant normals
$-\sigma_c w_c$ and $-\sigma_{c'} w_{c'}$ of quadrant $\sigma = (\sigma_c,
\sigma_{c'}) \in \{\pm 1\}^2$ determine its interior dihedral angle
\begin{equation}\label{eq:model-theta-sigma}
  \theta_\sigma(x_0)
  := \arccos\!\left(
      -\frac{\sigma_c \sigma_{c'}\, w_c\transp G_\sigma(x_0)^{-1} w_{c'}}
            {\sqrt{\bigl(w_c\transp G_\sigma(x_0)^{-1} w_c\bigr)
                   \bigl(w_{c'}\transp G_\sigma(x_0)^{-1} w_{c'}\bigr)}}
    \right) \in (0, \pi),
\end{equation}
together with the pointwise Regge deficit and its quadratic proxy,
\begin{align*}
  \varepsilon_B(x_0)
  &:= 2\pi - \sum_{\sigma \in \{\pm 1\}^2} \theta_\sigma(x_0), \\
  \varepsilon_0(x_0)
  &:= -\frac{1}{4}\, a_1 a_2 \rho \sqrt{1 - \rho^2}
    = -\frac{1}{4}\, \ip{w_c}{G_B(x_0)^{-1} w_{c'}}\sqrt{\det M_B(x_0)},
\end{align*}
where $(2\pi)^{-1}\varepsilon_B(x_0)$ is Banchoff's vertex curvature at $\dimin
= 2$~\cite{Banchoff1967} and the normalised transverse
Regge--Cheeger--M\"uller--Schrader angle deficit for $\dimin \ge 3$
(Remark~\ref{rem:bone-geometry}). Setting
$\varepsilon_{cc'}(x) := 0$ off $B_{cc'}^{\circ}$ and whenever $w_c \wedge
w_{c'} = 0$, the \emph{ambient Regge curvature measure} of a Borel window $\Omega
\subseteq \RR^{\dimin}$ and its slice supremum over the affine Grassmannian
$\mathrm{Gr}_{\dimin-2}^{\mathrm{aff}}(\RR^{\dimin})$ are
\begin{equation}\label{eq:model-measures}
  \begin{gathered}
    \Regge(\Omega)
    := \sum_{1 \le c < c' \le \ngen} Y_{cc'}(\Omega),
    \qquad
    Y_{cc'}(\Omega)
    := \int_{B_{cc'} \cap \Omega} \varepsilon_{cc'}(x)\, \dd\Haus{\dimin-2}(x), \\
    \beta_\Omega
    := \sup_{F \in \mathrm{Gr}_{\dimin-2}^{\mathrm{aff}}(\RR^{\dimin})}
       \Haus{\dimin-2}(F \cap \Omega),
  \end{gathered}
\end{equation}
where $\Haus{\dimin-2}$ is $(\dimin - 2)$-dimensional Hausdorff measure
(counting measure at $\dimin = 2$) and $\ell_B(\Omega) = \ell_{cc'}(\Omega) :=
\Haus{\dimin-2}(B_{cc'} \cap \Omega)$ denotes the bone slice volume. Note that
$\beta_\Omega \le 1$ on $\RR^2$, while $\beta_\Omega < \infty$ for every bounded
Borel window $\Omega \subset \RR^{\dimin}$ when $\dimin \ge 3$.
\end{definition}

\begin{remark}
\label{rem:bone-geometry}
Across any wall $\Pi_c$, adjacent cell metrics $G_{\tau^+} = G_{\tau^-} + w_c
w_c\transp$ agree on $T\Pi_c = w_c^\perp$, so $(\RR^{\dimin}, G)$ is an
intrinsic polyhedral manifold with zero wall crease curvature. Near $x_0 \in
B_{cc'}^{\circ}$, quadrant $\sigma$ splits metrically as a Riemannian product.
One factor is $(TB_{cc'}, G_\sigma|_{TB_{cc'}})$. The other is a planar sector
in the $G_\sigma$-orthogonal complement $(TB_{cc'})^{\perp_{G_\sigma(x_0)}} =
\spn(G_\sigma(x_0)^{-1}w_c, G_\sigma(x_0)^{-1}w_{c'})$, which differs from $N_B
:= \spn(w_c, w_{c'})$ when $\dimin \ge 3$. The sector is bounded by the outward
conormals $-\sigma_c G_\sigma(x_0)^{-1}w_c$ and $-\sigma_{c'}
G_\sigma(x_0)^{-1}w_{c'}$, and its interior angle is $\theta_\sigma(x_0)$
in~\eqref{eq:model-theta-sigma}. Thus $\varepsilon_B(x_0)$
coincides with the intrinsic codimension-$2$ Regge--Lipschitz--Killing
deficit~\cite{Regge1961,Banchoff1967,CheegerMullerSchrader1984}. When $\norm{w_c
\wedge w_{c'}} > 0$, orthogonal projection of $0$ onto $B_{cc'}$ yields the
\emph{minimum-norm bone representative} $x_0^{\min}(B) = x_0^{\min}(w_c, b_c,
w_{c'}, b_{c'}) := -L_B(L_B\transp L_B)^{-1}(b_c, b_{c'})\transp \in B_{cc'}
\cap N_B$. This makes $\Regge(\Omega)$ a finite random variable whenever
$\beta_\Omega < \infty$.\footnote{Choose $Q_B \in \RR^{\dimin \times
(\dimin-2)}$ ($Q_B\transp Q_B = I_{\dimin-2}$, $\im Q_B = T B_{cc'}$) measurably
in $(w_c, b_c, w_{c'}, b_{c'})$ via Gram--Schmidt on $\{w_c \wedge w_{c'} \ne
0\}$. Set $x(y) := x_0^{\min}(B) + Q_B y$ ($y \in \RR^{\dimin-2}$). Then the
maps $(W, b, y) \mapsto I_h(x(y))$, $G_B(x(y))^{-1}$ and
$\varepsilon_B(x(y))\mathbf{1}_{B_{cc'}^\circ}(x(y))$ are jointly Borel
measurable. Since $\abs{\varepsilon_B} \le \pi/2$ and $\varepsilon_{cc'} := 0$
off $B_{cc'}^\circ$, $\Regge(\Omega)$ in~\eqref{eq:model-measures} is well
defined and bounded by $\binom{\ngen}{2}\frac{\pi}{2}\beta_\Omega$ for every
$(W, b)$.}
\end{remark}

Let $\mathcal{E}_{\mathrm{nd}}(W, b)$ denote the \emph{general-position event}
that every $k \le \dimin$ normals have rank $k$ and no $\dimin + 1$ hyperplanes
share a common point (waiving the latter condition when $b \equiv 0$ and $\dimin
\ge 3$). On $\mathcal{E}_{\mathrm{nd}}(W, b)$, every triple intersection $\Pi_c
\cap \Pi_{c'} \cap \Pi_h$ ($h \notin \{c,c'\}$) is empty or of dimension $\dimin
- 3$. Thus $\Haus{\dimin-2}(B_{cc'} \setminus B_{cc'}^{\circ}) = 0$, and
$\Regge(\Omega)$ integrates the intrinsic transverse deficit
$\varepsilon_{cc'}(x)$ against ambient Hausdorff measure $\Haus{\dimin-2}$ along
each bone. At $\dimin = 2$ this coincides with the intrinsic Gauss--Bonnet
curvature measure of $\mathcal{M}_{W,b}$. For $\dimin \ge 3$, the
intrinsic Cheeger--M\"uller--Schrader measure instead weights $B_{cc'}$ by $\sqrt{\det
G|_{TB_{cc'}}}$ (see Remark~\ref{rem:cms-intrinsic-equivalence}). Outside
$\mathcal{E}_{\mathrm{nd}}(W, b)$,~\eqref{eq:model-measures}
defines $\Regge(\Omega)$ directly by the pairwise bone sum. The sole case where
$B_{cc'}^{\circ} = \emptyset$ identically is the homogeneous planar arrangement
($\dimin = 2$, $b \equiv 0$). When $\norm{w_c \wedge w_{c'}} > 0$ for all $c <
c'$, all $\ngen$ lines meet at a conical apex at $0$. In this case, we define
\begin{equation}\label{eq:model-apex-convention}
  \Regge(\Omega) := \bigl(2\pi - \Theta_\infty(W)\bigr)\mathbf{1}\{0 \in \Omega\}
  \qquad \text{for every Borel } \Omega \subseteq \RR^2,
\end{equation}
with $\Theta_\infty(W)$ as in~\eqref{eq:planar-wedge-angle}. Whenever $b \equiv
0$ at $\dimin = 2$, we assume $\PP(w_c \wedge w_{c'} = 0) = 0$ for all $c \ne
c'$.

\subsection{Bone conditioning, Hoeffding decomposition and standing assumptions}
\label{subsec:model-assumptions}

\begin{definition}
\label{def:bone-conditioning}
Let $\xi_h := (w_h, b_h) \in \RR^{\dimin} \times \RR$ ($h \in [\ngen]$) be
mutually independent, write $\mathcal{F}_S := \sigma(\{\xi_h : h \in S\})$ for
$S \subseteq [\ngen]$, and write $\mathcal{F}_B := \sigma(w_c, b_c, w_{c'},
b_{c'})$ for the \emph{bone $\sigma$-field} of $B = \{c, c'\}$. For $X \in
L^2(\mathcal{F}_{[\ngen]})$, the \emph{Hoeffding projection} $\pi_A(X) \in
\mathcal{V}_A := \im(\pi_A) \subset L^2(\mathcal{F}_A)$ ($A \subseteq [\ngen]$)
and the \emph{$k$-th Hoeffding component} $P_k(X)$ ($0 \le k \le \ngen$) are
\begin{equation}\label{eq:model-P-k}
  \pi_A(X) := \sum_{S \subseteq A} (-1)^{\abs{A} - \abs{S}}\, \EE\bigl[X \;\big|\; \mathcal{F}_S\bigr],
  \qquad
  P_k(X) := \sum_{\substack{A \subseteq [\ngen] \\ \abs{A} = k}} \pi_A(X),
\end{equation}
satisfying $X = \sum_{k=0}^{\ngen} P_k(X)$, $\Var(X) = \sum_{k=1}^{\ngen}
\norm{P_k(X)}_{L^2(\PP)}^2$, $\mathcal{V}_A \perp \mathcal{V}_{A'}$ for $A \ne
A'$, and $\EE[\pi_A(X) \mid \mathcal{F}_B] = 0$ whenever $A \not\subseteq B$. In
particular, $P_0(X) = \EE[X]$, $P_1(X) = \sum_{h=1}^{\ngen} (\EE[X \mid \xi_h] -
\EE[X])$ and $P_2(X) = \sum_{c < c'} \pi_{\{c, c'\}}(X)$. We reserve the term
\emph{second Wiener chaos}~\cite{Janson1997} for the Gaussian quadratic-form
limits as $\ngen \to \infty$.
\end{definition}

\begin{definition}
\label{def:standing-assumptions}
Our probabilistic results refer to the following \emph{standing distributional
models} on $\xi_h = (w_h, b_h)$ ($h \in [\ngen]$):
\begin{enumerate}[label=\textup{(A\arabic*)},leftmargin=2.4em]
  \item\label{assump:A1} \emph{Independent centrally symmetric generators.}
    The pairs $\xi_1, \dots, \xi_{\ngen}$ are mutually independent in
    $\RR^{\dimin} \times \RR$, with $(-w_h, -b_h) \stackrel{d}{=} (w_h, b_h)$
    and $\PP(w_h = 0) = 0$ for each $h \in [\ngen]$ (for example when $w_h
    \stackrel{d}{=} -w_h$ and $b_h \indep w_h$ with $b_h \stackrel{d}{=} -b_h$,
    including $b_h \equiv 0$).
  \item\label{assump:A2} \emph{Gaussian normal model.}
    Assumption~\ref{assump:A1} holds and $w_1, \dots, w_{\ngen}
    \stackrel{\mathrm{i.i.d.}}{\sim} \Norm(0, \sigw^2 I_{\dimin})$ ($\sigw > 0$)
    independently of $b$. On a Borel window $\Omega \subseteq \RR^{\dimin}$ with
    $\beta_\Omega < \infty$, we say that \emph{\ref{assump:A2} holds on $\Omega$
    with symmetric i.i.d.\ offsets} when $b_1, \dots, b_{\ngen}$ are i.i.d.\
    symmetric. We require $\PP(b_1 = 0) = 0$ if $\dimin = 2$, but allow an
    arbitrary symmetric law (including $b \equiv 0$) if $\dimin \ge 3$.
  \item\label{assump:A3} \emph{Elliptically symmetric normal model.}
    The normals $w_1, \dots, w_{\ngen} \stackrel{\mathrm{i.i.d.}}{\sim} P_w$
    have covariance $\Sigma_w = \EE[w_1 w_1\transp] \in
    \mathrm{Sym}_{\dimin}^{++}$, rotationally invariant whitened vector
    $\tilde{w}_h := \Sigma_w^{-1/2} w_h$, $\PP(w_h = 0) = 0$ and finite fourth
    radial moment $\tau_2(P_w) := \EE[\norm{\tilde{w}_1}^4/\dimin^2] < \infty$.
    We call~\ref{assump:A3} \emph{isotropic} when $\Sigma_w = \sigw^2 I_{\dimin}$
    ($\sigw > 0$).
\end{enumerate}
At an $\mathcal{F}_B$-measurable point $x_0 \in B_{cc'}$, conditional boundary
non-atomicity means that for every $h \notin \{c, c'\}$,
\begin{equation}\label{eq:model-boundary-null}
  \PP\bigl(\ip{w_h}{x_0} + b_h = 0 \;\big|\; \mathcal{F}_B\bigr) = 0
  \qquad \text{almost surely.}
\end{equation}
Under~\ref{assump:A2}, this holds whenever $\PP(x_0 = 0) = 0$. Hence it holds
$\Haus{\dimin-2}$-a.e.\ on $B_{cc'}$ for any $b$ (including $b \equiv 0$) when
$\dimin \ge 3$, and at $x_{cc'}$ whenever $\PP(b_c = b_{c'} = 0) = 0$ when
$\dimin = 2$. On balls, $\beta_{B(x_*, r)} = \omega_{\dimin - 2} r^{\dimin -
2}$ with $\omega_k := \pi^{k/2}/\Gamma(k/2 + 1)$ ($\omega_0 = 1$).
\end{definition}
\section{Local bone geometry, four-arccos formula and two-sided bounds}
\label{sec:local-geometry}

Every identity and estimate in this section is deterministic and holds at a
single bone $B = B_{cc'}$ for arbitrary $G_B \in \mathrm{Sym}_{\dimin}^{++}$ and
linearly independent $w_c, w_{c'} \in \RR^{\dimin}$.

\subsection{Woodbury reduction and the four-arccos formula}
\label{subsec:woodbury-four-arccos}

\begin{lemma}
\label{lem:sector-gram}
Let $\dimin \ge 2$, $G_B \in \mathrm{Sym}_{\dimin}^{++}$, and $w_c, w_{c'} \in
\RR^{\dimin}$ with $w_c \wedge w_{c'} \ne 0$. Let $L_B = [w_c \;\; w_{c'}]$ and
let
\[
  M_B = L_B\transp G_B^{-1} L_B
  = \begin{pmatrix} a_1 & \rho\sqrt{a_1 a_2} \\ \rho\sqrt{a_1 a_2} & a_2 \end{pmatrix}
  \in \mathrm{Sym}_2^{++}
\]
be as in~\eqref{eq:model-M-Gamma}. Set $S_\sigma := \diag(\mathbf{1}\{\sigma_c =
+1\}, \mathbf{1}\{\sigma_{c'} = +1\})$ for $\sigma \in \{\pm 1\}^2$, so that
$G_\sigma = G_B + L_B S_\sigma L_B\transp$, and let $M_\sigma := L_B\transp
G_\sigma^{-1} L_B \in \mathrm{Sym}_2^{++}$. Then $I_2 + S_\sigma M_B$ is
invertible with $d_\sigma := \det(I_2 + S_\sigma M_B) \ge 1$, $\Delta :=
\det(I_2 + M_B) = 1 + a_1 + a_2 + (1 - \rho^2) a_1 a_2 > 1$, and by the Woodbury
identity,
\begin{equation}\label{eq:woodbury-sector}
  M_\sigma
  = L_B\transp G_\sigma^{-1} L_B
  = M_B (I_2 + S_\sigma M_B)^{-1}
  = (I_2 + M_B S_\sigma)^{-1} M_B.
\end{equation}
With $\alpha_1 := a_1(1 - \rho^2) > 0$ and $\alpha_2 := a_2(1 - \rho^2) > 0$, we
have $(M_\sigma)_{12} = d_\sigma^{-1}\rho\sqrt{a_1 a_2}$ and
$(M_\sigma)_{11}(M_\sigma)_{22} = d_\sigma^{-2} a_1 a_2\,(1 +
\mathbf{1}\{\sigma_{c'}=+1\}\alpha_2)\allowbreak(1 +
\mathbf{1}\{\sigma_c=+1\}\alpha_1)$.
\end{lemma}

\begin{proof}
Multiplying $G_B^{-1} - G_\sigma^{-1} = G_B^{-1} L_B S_\sigma L_B\transp
G_\sigma^{-1}$ by $L_B\transp$ and $L_B$ gives $(I_2 + M_B S_\sigma)M_\sigma =
M_B$. Since $I_2 + M_B S_\sigma$ is similar to $I_2 + M_B^{1/2} S_\sigma
M_B^{1/2} \succeq I_2$, we have $d_\sigma \ge 1$ and~\eqref{eq:woodbury-sector}.
Inverting $I_2 + M_B S_\sigma$ via the $2 \times 2$ adjugate formula with $\det
M_B = a_1\alpha_2 = a_2\alpha_1$ yields
\begin{gather*}
  M_{(-,-)} = M_B,
  \qquad
  M_{(+,+)} = \frac{1}{\Delta}
  \begin{pmatrix} a_1(1+\alpha_2) & \rho\sqrt{a_1 a_2} \\
    \rho\sqrt{a_1 a_2} & a_2(1+\alpha_1) \end{pmatrix},
  \\
  M_{(+,-)} = \frac{1}{1+a_1}
  \begin{pmatrix} a_1 & \rho\sqrt{a_1 a_2} \\
    \rho\sqrt{a_1 a_2} & a_2(1+\alpha_1) \end{pmatrix},
\end{gather*}
with $d_{(-,-)} = 1$, $d_{(+,+)} = \Delta$ and $d_{(+,-)} = 1 + a_1$
respectively, and $M_{(-,+)}$ follows by exchanging $1$ and $2$. Both entry
formulas for $M_\sigma$ follow directly.
\end{proof}

\begin{theorem}
\label{thm:four-arccos}
In the notation of Lemma~\ref{lem:sector-gram}, let $(a_1, a_2, \rho) \in (0,
\infty)^2 \times (-1, 1)$ and $\alpha_i := a_i(1 - \rho^2) > 0$ ($i = 1, 2$).
The four quadrant dihedral angles $\theta_\sigma \in (0, \pi)$
in~\eqref{eq:model-theta-sigma} are $\theta_{(-,-)} = \arccos(-\rho)$,
$\theta_{(+,-)} = \arccos(\rho/\sqrt{1+\alpha_1})$, $\theta_{(-,+)} =
\arccos(\rho/\sqrt{1+\alpha_2})$ and $\theta_{(+,+)} =
\arccos(-\rho/\sqrt{(1+\alpha_1)(1+\alpha_2)})$, and the pointwise Regge deficit
$\varepsilon_B = \Psi(a_1, a_2, \rho) := 2\pi - \sum_{\sigma \in \{\pm 1\}^2}
\theta_\sigma$ satisfies the four-arccos identity
\begin{align}
  \Psi(a_1, a_2, \rho)
  &= \arccos(\rho)
   - \arccos\bigl(\tfrac{\rho}{\sqrt{1 + \alpha_1}}\bigr)
   - \arccos\bigl(\tfrac{\rho}{\sqrt{1 + \alpha_2}}\bigr)
   + \arccos\bigl(\tfrac{\rho}{\sqrt{(1 + \alpha_1)(1 + \alpha_2)}}\bigr)
  \label{eq:four-arccos} \\
  &= -\arcsin(\rho)
   + \arcsin\bigl(\tfrac{\rho}{\sqrt{1 + \alpha_1}}\bigr)
   + \arcsin\bigl(\tfrac{\rho}{\sqrt{1 + \alpha_2}}\bigr)
   - \arcsin\bigl(\tfrac{\rho}{\sqrt{(1 + \alpha_1)(1 + \alpha_2)}}\bigr).
  \label{eq:four-arcsin}
\end{align}
\end{theorem}

\begin{proof}
Substituting Lemma~\ref{lem:sector-gram} into~\eqref{eq:model-theta-sigma}
cancels $d_\sigma^{-1}\sqrt{a_1 a_2}$ between $(M_\sigma)_{12}$ and
$\sqrt{(M_\sigma)_{11}(M_\sigma)_{22}}$, yielding $\theta_\sigma$. Applying
$\arccos(-z) = \pi - \arccos(z)$ to $\theta_{(-,-)}$ and $\theta_{(+,+)}$ gives~\eqref{eq:four-arccos},
and $\arccos(z) = \frac{\pi}{2} - \arcsin(z)$ gives~\eqref{eq:four-arcsin}.
\end{proof}

\subsection{Integral representation, two-sided bounds and smoothness}
\label{subsec:sandwich-smoothness}

\begin{proposition}
\label{prop:deficit-integral}
For $(a_1, a_2, \rho) \in (0, \infty)^2 \times (-1, 1)$, let $\Psi(a_1, a_2,
\rho)$ be given by~\eqref{eq:four-arccos}, set $\varepsilon_0(a_1, a_2, \rho) :=
-\frac{1}{4} a_1 a_2 \rho \sqrt{1 - \rho^2}$, and for $(u_1, u_2, d) \in [0,
\infty)^3$ define
\[
  \mathcal{J}(u_1, u_2, d)
  := \int_0^1\!\! \int_0^1 \frac{\dd u\, \dd v}{(1 + u_1 u + u_2 v + d u v)^{3/2}}.
\]
\begin{enumerate}[label=(\roman*)]
  \item With $\alpha_i := a_i(1-\rho^2)$ and $M = \bigl(\begin{smallmatrix}
    m_{11} & m_{12} \\ m_{12} & m_{22} \end{smallmatrix}\bigr) =
    \bigl(\begin{smallmatrix} a_1 & \rho\sqrt{a_1 a_2} \\ \rho\sqrt{a_1 a_2} &
    a_2 \end{smallmatrix}\bigr) \in \mathrm{Sym}_2^{++}$ ($\det M = a_1
    a_2(1-\rho^2)$),
    \begin{align}
      \Psi(a_1, a_2, \rho)
      &= -\frac{\rho\sqrt{1-\rho^2}}{4}\! \int_0^{a_1}\!\!\!\! \int_0^{a_2}\!
         \frac{\dd t\, \dd s}{\bigl(1 + s + t + s t (1-\rho^2)\bigr)^{3/2}}
      \label{eq:deficit-double-integral} \\
      &= \varepsilon_0(a_1, a_2, \rho)\,\mathcal{J}(a_1, a_2, \det M)
       = -\frac{1}{4}\, m_{12} \sqrt{\det M}\,\mathcal{J}(m_{11}, m_{22}, \det M)
      \notag \\
      &=: \Psi_{\mathrm{mat}}(M).
      \label{eq:deficit-Psi-M}
    \end{align}
  \item $\Psi(a_1, a_2, 0) = 0$, and for $\rho \ne 0$, $\abs{\Psi(a_1, a_2,
    \rho)}$ increases strictly in $a_1, a_2 > 0$ to $\abs{\arcsin(\rho)}$:
    \begin{equation}\label{eq:deficit-sign}
      \sgn\bigl(\Psi(a_1, a_2, \rho)\bigr)
      = \sgn\bigl(\varepsilon_0(a_1, a_2, \rho)\bigr)
      = -\sgn(\rho),
      \quad
      \abs{\Psi(a_1, a_2, \rho)}
      < \abs{\arcsin(\rho)} < \frac{\pi}{2}.
    \end{equation}
  \item For all $(a_1, a_2, \rho) \in (0, \infty)^2 \times (-1, 1)$, $0 \le 1 -
    \mathcal{J}(a_1, a_2, \det M) \le \min\{1, a_1 + a_2\}$ and
    \begin{align}
      \frac{\abs{\varepsilon_0(a_1, a_2, \rho)}}{(1+a_1)(1+a_2)}
      &\le \abs{\Psi(a_1, a_2, \rho)}
       \le \abs{\varepsilon_0(a_1, a_2, \rho)},
      \label{eq:deficit-sandwich} \\
      \bigl|\Psi(a_1, a_2, \rho) - \varepsilon_0(a_1, a_2, \rho)\bigr|
      &\le \min\{1,\, a_1 + a_2\}\,\abs{\varepsilon_0(a_1, a_2, \rho)}.
      \label{eq:deficit-linear-rem}
    \end{align}
\end{enumerate}
\end{proposition}

\begin{proof}
Setting $g(\tilde{u}, \tilde{v}) :=
\arccos\bigl(\rho\big/\sqrt{(1+\tilde{u})(1+\tilde{v})}\bigr)$ on $[0,
\infty)^2$, the identity~\eqref{eq:four-arccos} gives $\Psi(a_1, a_2, \rho) =
\int_0^{\alpha_1}\!\int_0^{\alpha_2} \partial_{\tilde{u}}\partial_{\tilde{v}}
g(\tilde{u}, \tilde{v})\,\dd\tilde{v}\,\dd\tilde{u}$ with
$\partial_{\tilde{u}}\partial_{\tilde{v}} g = -\frac{1}{4}\rho\bigl((1+\tilde{u})(1+\tilde{v})
- \rho^2\bigr)^{-3/2}$. Substituting $\tilde{u} = s(1-\rho^2), \tilde{v} =
t(1-\rho^2)$ and then $s = a_1 u, t = a_2 v$
yields~\eqref{eq:deficit-double-integral}--\eqref{eq:deficit-Psi-M}. Positivity
of the integrand and~\eqref{eq:four-arcsin} give~\eqref{eq:deficit-sign}. Using
$1 \le 1 + s + t + s t (1 - \rho^2) \le (1 + s)(1 + t)$ together with
$\int_0^{a_i} (1+s)^{-3/2}\,\dd s = 2a_i/[(1+a_i)(1 + (1+a_i)^{-1/2})] \ge
a_i/(1+a_i)$ gives~\eqref{eq:deficit-sandwich}. Finally, $0 < \mathcal{J} \le 1$
gives $1 - \mathcal{J} \le 1$, while for $0 < a_1 + a_2 < 1$ we have $\det M \le
a_1 a_2 \le \frac{1}{4}(a_1 + a_2)^2 < \frac{1}{4}(a_1 + a_2)$, so integrating
$0 \le 1 - (1 + \xi)^{-3/2} \le \frac{3}{2}\xi$ at $\xi = a_1 u + a_2 v + (\det
M)u v$ over $[0, 1]^2$ yields $1 - \mathcal{J}(a_1, a_2, \det M) \le
\frac{3}{4}(a_1 + a_2) + \frac{3}{8}\det M < a_1 + a_2$
and~\eqref{eq:deficit-linear-rem}.
\end{proof}

\begin{lemma}
\label{lem:deficit-smoothness}
By~\eqref{eq:deficit-Psi-M}, for $M = \bigl(\begin{smallmatrix} m_{11} & m_{12}
\\ m_{12} & m_{22} \end{smallmatrix}\bigr) \in \mathrm{Sym}_2^{++}$,
\begin{equation}\label{eq:simplex-Psi-Phi}
  \Psi_{\mathrm{mat}}(M)
  = -\frac{1}{4}\, m_{12}\sqrt{\det M}\,\mathcal{J}(m_{11}, m_{22}, \det M).
\end{equation}
\begin{enumerate}[label=(\roman*)]
  \item There exists a universal constant $C > 0$ such that on $[0, \infty)^3$,
    $0 < \mathcal{J}(u_1, u_2, d) \le 1$, $0 \le 1 - \mathcal{J}(u_1, u_2, d) \le
    C(u_1 + u_2 + d)$, and $\abs{\partial^\beta \mathcal{J}(u_1, u_2, d)} \le C$
    for every multi-index $\beta \in \ZZ_{\ge 0}^3$ with $1 \le \abs{\beta} \le
    3$.
  \item Let $M_0 = \bigl(\begin{smallmatrix} a_{1,0} & \rho_0\sqrt{a_{1,0} a_{2,0}} \\
    \rho_0\sqrt{a_{1,0} a_{2,0}} & a_{2,0} \end{smallmatrix}\bigr) \in
    \mathrm{Sym}_2^{++}$ satisfy $\tr(M_0) = a_{1,0} + a_{2,0} \le 1/2$, let
    $\widetilde{T}_B \in \mathrm{GL}_2(\RR)$ be any matrix with
    $\widetilde{T}_B\transp \widetilde{T}_B = M_0$ (in
    Lemma~\ref{lem:deficit-perturbation} it will be the triangular factor of a
    thin QR factorisation), and
    set $\mathcal{B}_{1/2}(I_2) := \{S \in \mathrm{Sym}_2^{++} : \norm{S -
    I_2}_{\mathrm{op}} \le 1/2\}$. There exists a universal constant $C > 0$
    such that for all $S \in \mathcal{B}_{1/2}(I_2)$ and $U_1, U_2, U_3 \in
    \mathrm{Sym}_2(\RR)$, $\abs{\Psi_{\mathrm{mat}}(\widetilde{T}_B\transp S
    \widetilde{T}_B)} \le C\,a_{1,0} a_{2,0}\sqrt{1 - \rho_0^2}$ and, for $r \in
    \{1, 2, 3\}$,
    \begin{equation}\label{eq:simplex-K-bounds}
      \bigl|D_S^r\bigl[\Psi_{\mathrm{mat}}(\widetilde{T}_B\transp S \widetilde{T}_B)\bigr][U_1, \dots, U_r]\bigr|
      \le C\, a_{1,0} a_{2,0} \sqrt{1 - \rho_0^2}\prod_{i=1}^r \norm{U_i}_{\mathrm{op}}.
    \end{equation}
\end{enumerate}
\end{lemma}

\begin{proof}
Since $1 + u_1 u + u_2 v + d u v \ge 1$ on $[0, 1]^2$, differentiating
$\mathcal{J}$ under the integral sign bounds $\abs{\partial^\beta \mathcal{J}}$
uniformly on $[0, \infty)^3$ for $1 \le \abs{\beta} \le 3$, and the mean value
theorem from $\mathcal{J}(0, 0, 0) = 1$ gives part~(i). For part~(ii), write
$\widetilde{T}_B = [\sqrt{a_{1,0}}\,\hat{r}_1 \;\; \sqrt{a_{2,0}}\,\hat{r}_2]$
with unit vectors $\hat{r}_1, \hat{r}_2 \in \RR^2$ ($\hat{r}_1\transp\hat{r}_2 =
\rho_0$, $d_0 := \det M_0 \le 1/16$), so that $\Psi_{\mathrm{mat}}(\widetilde{T}_B\transp
S \widetilde{T}_B) = -\frac{1}{4}a_{1,0}a_{2,0}\sqrt{1 - \rho_0^2}\,(\hat{r}_1\transp
S\hat{r}_2)\sqrt{\det S}\,\mathcal{J}(a_{1,0}\hat{r}_1\transp S\hat{r}_1,
a_{2,0}\hat{r}_2\transp S\hat{r}_2, d_0\det S)$ by~\eqref{eq:simplex-Psi-Phi}.
On $\mathcal{B}_{1/2}(I_2)$ the eigenvalues of $S$ lie in $[1/2, 3/2]$. Thus $\det
S \in [1/4, 9/4]$ is bounded away from zero, and $s \mapsto \sqrt{s}$ has bounded
derivatives up to order $3$. Likewise, the three polynomial arguments of
$\mathcal{J}$ map $\mathcal{B}_{1/2}(I_2)$ into $[0, \infty)^3$ with Fr\'echet
derivatives up to order $3$ bounded uniformly in $M_0$ (as $a_{1,0}, a_{2,0} \le
1/2$ and $d_0 \le 1/16$). Combining the chain rule with part~(i) and the
three-factor Leibniz rule yields~\eqref{eq:simplex-K-bounds}.
\end{proof}

\begin{lemma}
\label{lem:deficit-perturbation}
Let $\dimin \ge 2$, let $L_B = [w_c \;\; w_{c'}] \in \RR^{\dimin \times 2}$ have
$\rank(L_B) = 2$, and set $\varepsilon_B(G) := \Psi_{\mathrm{mat}}(L_B\transp
G^{-1} L_B)$ for $G \in \mathrm{Sym}_{\dimin}^{++}$. Fix $G_0 \in
\mathrm{Sym}_{\dimin}^{++}$ such that $M_0 := L_B\transp G_0^{-1} L_B =
\bigl(\begin{smallmatrix} a_{1,0} & \rho_0\sqrt{a_{1,0} a_{2,0}} \\
\rho_0\sqrt{a_{1,0} a_{2,0}} & a_{2,0} \end{smallmatrix}\bigr)$ satisfies
$\tr(M_0) = a_{1,0} + a_{2,0} \le 1/2$, and let $G_0^{-1/2} L_B = \widetilde{Q}_N
\widetilde{T}_B$ be a thin QR factorisation ($\widetilde{Q}_N\transp
\widetilde{Q}_N = I_2$, $\widetilde{T}_B\transp \widetilde{T}_B = M_0$). Set
$\Upsilon := G_0^{-1/2}(G - G_0)G_0^{-1/2}$. There exists a universal constant
$C > 0$ such that whenever $\norm{\Upsilon}_{\mathrm{op}} \le 1/3$:
\begin{enumerate}[label=(\roman*)]
  \item The first and second Fr\'echet derivatives $\mathcal{L}_{G_0,
    L_B}[\Upsilon_0]$ and $\mathcal{Q}_{G_0, L_B}[\Upsilon_0, \Upsilon_0]$ of
    $\Upsilon \mapsto \varepsilon_B(G_0^{1/2}(I_{\dimin}+\Upsilon)G_0^{1/2})$ at
    $\Upsilon = 0$ satisfy
    \begin{align}
      \abs{\mathcal{L}_{G_0, L_B}[\Upsilon_0]}
      &\le C\, a_{1,0}a_{2,0}\sqrt{1 - \rho_0^2}\,
        \norm{\widetilde{Q}_N\transp \Upsilon_0 \widetilde{Q}_N}_{\mathrm{op}},
      \label{eq:pert-taylor-L} \\
      \tfrac{1}{2}\abs{\mathcal{Q}_{G_0, L_B}[\Upsilon_0, \Upsilon_0]}
      &\le C\,a_{1,0}a_{2,0}\sqrt{1 - \rho_0^2}\,\norm{\Upsilon_0}_{\mathrm{op}}^2,
      \label{eq:pert-taylor-Q} \\
      \bigl|\varepsilon_B(G) - \varepsilon_B(G_0)
            - \mathcal{L}_{G_0, L_B}[\Upsilon]
            &\notag \\
            {}- \tfrac{1}{2}\mathcal{Q}_{G_0, L_B}[\Upsilon, \Upsilon]\bigr|
      &\le C\, a_{1,0}a_{2,0}\sqrt{1 - \rho_0^2}\,\norm{\Upsilon}_{\mathrm{op}}^3.
      \label{eq:pert-taylor-R3}
    \end{align}
  \item For $E_1, E_2 \in \mathrm{Sym}_{\dimin}(\RR)$ with $\Upsilon_i
    := G_0^{-1/2}E_i G_0^{-1/2}$ satisfying $\norm{\Upsilon_1}_{\mathrm{op}}
    + \norm{\Upsilon_2}_{\mathrm{op}} \le 1/3$, and for $w_j \in \RR^{\dimin}$
    with $\norm{G_0^{-1/2}w_j}_2^2 \le 1/3$,
    \begin{align}
      \bigl|\varepsilon_B(G_0 + E_1) - \varepsilon_B(G_0)\bigr|
      &\le C\, a_{1,0}a_{2,0}\sqrt{1 - \rho_0^2}\,\norm{\Upsilon_1}_{\mathrm{op}},
      \label{eq:pert-diff-1} \\
      \bigl|\varepsilon_B(G_0 \pm w_j w_j\transp) - \varepsilon_B(G_0)\bigr|
      &\le C\, a_{1,0}a_{2,0}\sqrt{1 - \rho_0^2}\,\norm{\widetilde{Q}_N\transp G_0^{-1/2} w_j}_2^2,
      \label{eq:pert-diff-1-rank1} \\
      \bigl|\varepsilon_B(G_0 + E_1 + E_2)
        - \varepsilon_B(G_0 + E_1)
        &\notag \\
        {}- \varepsilon_B(G_0 + E_2)
        + \varepsilon_B(G_0)\bigr|
      &\le C\, a_{1,0}a_{2,0}\sqrt{1 - \rho_0^2}\,
        \norm{\Upsilon_1}_{\mathrm{op}}\norm{\Upsilon_2}_{\mathrm{op}}.
      \label{eq:pert-diff-2}
    \end{align}
\end{enumerate}
\end{lemma}

\begin{proof}
Write $f(\Upsilon) := \varepsilon_B\bigl(G_0^{1/2}(I_{\dimin} +
\Upsilon)G_0^{1/2}\bigr) = g(S(\Upsilon))$ with $g(S) :=
\Psi_{\mathrm{mat}}(\widetilde{T}_B\transp S \widetilde{T}_B)$, $R(\Upsilon) :=
(I_{\dimin} + \Upsilon)^{-1}$ and $S(\Upsilon) := \widetilde{Q}_N\transp
R(\Upsilon) \widetilde{Q}_N$. When $\norm{\Upsilon}_{\mathrm{op}} \le 1/3$, the
Neumann series $R(\Upsilon) = \sum_{k=0}^\infty (-\Upsilon)^k$ gives
$\norm{R(\Upsilon)}_{\mathrm{op}} \le 3/2$ and $\norm{R(\Upsilon) -
I_{\dimin}}_{\mathrm{op}} \le \frac{3}{2}\norm{\Upsilon}_{\mathrm{op}} \le 1/2$,
so $S(\Upsilon) \in \mathcal{B}_{1/2}(I_2)$. Differentiating $(I_{\dimin} +
\Upsilon)R(\Upsilon) = I_{\dimin}$ inductively gives $DS(0)[\Upsilon_0] =
-\widetilde{Q}_N\transp \Upsilon_0 \widetilde{Q}_N$ and
$D^r S(\Upsilon)[\Upsilon_1, \dots, \Upsilon_r] = \allowbreak (-1)^r
\sum_{\sigma \in \mathfrak{S}_r} \allowbreak \widetilde{Q}_N\transp
R(\Upsilon)\Upsilon_{\sigma(1)} R(\Upsilon)\cdots \allowbreak
\Upsilon_{\sigma(r)} R(\Upsilon)\widetilde{Q}_N$, whence
$\norm{D^r S(\Upsilon)[\Upsilon_1, \dots, \Upsilon_r]}_{\mathrm{op}} \le
r!\,(3/2)^{r+1}\prod_{i=1}^r \norm{\Upsilon_i}_{\mathrm{op}}$ ($r \le 3$).
Combining the chain rule for $g \circ S$ with
Lemma~\ref{lem:deficit-smoothness}(ii) bounds $\mathcal{L}_{G_0,L_B}$,
$\mathcal{Q}_{G_0,L_B}$ and $D^r f$ ($r \le 3$). From this, Taylor's theorem and
integrating $Df(s\Upsilon_1)$ on $[0,1]$ and $D^2 f(s\Upsilon_1 + t\Upsilon_2)$
on $[0,1]^2$ yield~\eqref{eq:pert-taylor-R3},~\eqref{eq:pert-diff-1}
and~\eqref{eq:pert-diff-2}.

Finally, for $u_j := G_0^{-1/2}w_j$ with
$\norm{u_j}_2^2 \le 1/3$ and $t \in [0,1]$, Sherman--Morrison gives $S(\pm t u_j
u_j\transp) = I_2 \mp \frac{t}{1 \pm t\norm{u_j}_2^2}(\widetilde{Q}_N\transp
u_j)(\widetilde{Q}_N\transp u_j)\transp$. This implies $\norm{\frac{\dd}{\dd t}S(\pm t u_j
u_j\transp)}_{\mathrm{op}} = (1 \pm t\norm{u_j}_2^2)^{-2}\norm{\widetilde{Q}_N\transp
u_j}_2^2 \le \frac{9}{4}\norm{\widetilde{Q}_N\transp u_j}_2^2$. Integrating
$\frac{\dd}{\dd t}g(S(\pm t u_j u_j\transp))$ over $[0,1]$ then
yields~\eqref{eq:pert-diff-1-rank1}.
\end{proof}
\section{Antipodal parity, the thinned-Wishart law and the second Hoeffding component}
\label{sec:parity-kernel}

Throughout this section, for a bone $B = \{c, c'\}$ with $w_c \wedge w_{c'} \ne
0$, we write $R_c := \norm{w_c}/\sigw$, $R_{c'} := \norm{w_{c'}}/\sigw$,
$\vartheta_B := \arccos(\ip{w_c}{w_{c'}}/(\norm{w_c}\norm{w_{c'}})) \in (0,
\pi)$, and $\nu := m - \dimin + 2$ when $m_B(x_0) = m \ge \dimin$.

\subsection{Sign-flip oddness and the conditional thinned-Wishart law}
\label{subsec:parity-wishart}

\begin{lemma}
\label{lem:sign-flip-oddness}
Fix $B = \{c, c'\} \subset [\ngen]$ ($1 \le c
< c' \le \ngen$). For $(W, b)$ with $w_c \wedge w_{c'} \ne 0$, let
$\mathcal{I}_c\colon (w_c, b_c) \mapsto (-w_c, -b_c)$ and
$\mathcal{I}_{c'}\colon (w_{c'}, b_{c'}) \mapsto (-w_{c'}, -b_{c'})$ act on a
single generator while leaving all other $(w_h, b_h)$ unchanged. Both
involutions preserve $\Pi_c$, $\Pi_{c'}$, $B_{cc'}$, $B_{cc'}^{\circ}$,
$x_0^{\min}(B)$ and $\Haus{\dimin-2}|_{B_{cc'}}$. At every $x_0 \in
B_{cc'}^{\circ}$, they preserve $A(x_0)$, $G_B(x_0)$, $a_1(x_0)$ and $a_2(x_0)$
while negating $\rho(x_0)$, $\varepsilon_0(x_0)$ and $\varepsilon_B(x_0)$:
\[
  \varepsilon_B\bigl(\mathcal{I}_c(W, b);\, x_0\bigr)
  = \varepsilon_B\bigl(\mathcal{I}_{c'}(W, b);\, x_0\bigr)
  = -\varepsilon_B\bigl((W, b);\, x_0\bigr),
\]
so $\mathcal{I}_c \circ \mathcal{I}_{c'}$ fixes $\varepsilon_B(x_0)$ and
$Y_{cc'}(\Omega) = \int_{B_{cc'} \cap \Omega} \varepsilon_{cc'}(x)\,\dd\Haus{\dimin-2}(x)$
is odd under $\mathcal{I}_c$ and under $\mathcal{I}_{c'}$ separately.
\end{lemma}

\begin{proof}
Since $\ip{-w_c}{x} - b_c = -(\ip{w_c}{x} + b_c)$, $\mathcal{I}_c$ fixes
$\Pi_c$, $B_{cc'}$, $B_{cc'}^{\circ}$, $T B_{cc'}$ and every background
indicator $I_h(x_0)$ ($h \notin \{c, c'\}$), and so fixes $G_B(x_0)$, $a_1(x_0)$
and $a_2(x_0)$. However, $M_{B,12}(x_0) = w_c\transp G_B(x_0)^{-1} w_{c'}$
and $\rho(x_0)$ change sign, giving $\Psi(a_1, a_2, -\rho) = -\Psi(a_1, a_2,
\rho)$ by~\eqref{eq:deficit-double-integral}.
\end{proof}

\begin{lemma}
\label{lem:active-count}
Under~\ref{assump:A1}, fix a bone $B = \{c, c'\}$ and an
$\mathcal{F}_B$-measurable point $x_0 \in B_{cc'}$ on $\{w_c \wedge w_{c'} \ne
0\}$ satisfying~\eqref{eq:model-boundary-null}. Then \textup{(i)} conditionally
on $\mathcal{F}_B$, the indicators $\{I_h(x_0)\}_{h \notin \{c, c'\}}$ are
i.i.d.\ $\mathrm{Bernoulli}(1/2)$, so $m_B(x_0) \mid \mathcal{F}_B \sim
\mathrm{Bin}(\ngen - 2, \frac{1}{2})$; and \textup{(ii)} for each $h \notin \{c,
c'\}$, $w_h w_h\transp$ is conditionally independent of $I_h(x_0)$ given
$\mathcal{F}_B$, with conditional law given $\{\mathcal{F}_B, I_h(x_0) = 1\}$
equal to its unconditional law.
\end{lemma}

\begin{proof}
Conditionally on $\mathcal{F}_B$ (fixing $x_0$), the pairs $(w_h, b_h)_{h \notin
B}$ are independent with their unconditional laws, satisfy $(-w_h, -b_h)
\stackrel{d}{=} (w_h, b_h)$, and satisfy $\PP(\ip{w_h}{x_0} + b_h = 0 \mid
\mathcal{F}_B) = 0$ a.s.\ by~\eqref{eq:model-boundary-null}. For any bounded
Borel $\varphi\colon \mathrm{Sym}_{\dimin} \to \RR$, flipping $(w_h, b_h)
\mapsto (-w_h, -b_h)$ gives
\begin{align*}
  \EE\bigl[\varphi(w_h w_h\transp)\mathbf{1}\{\ip{w_h}{x_0} + b_h > 0\} \;\big|\; \mathcal{F}_B\bigr]
  &= \EE\bigl[\varphi(w_h w_h\transp)\mathbf{1}\{\ip{w_h}{x_0} + b_h < 0\} \;\big|\; \mathcal{F}_B\bigr] \\
  &= \frac{1}{2}\EE\bigl[\varphi(w_h w_h\transp)\bigr],
\end{align*}
proving \textup{(i)}--\textup{(ii)}.
\end{proof}

\begin{theorem}
\label{thm:conditional-wishart}
Under~\ref{assump:A2}, fix a bone $B = \{c, c'\}$ and an
$\mathcal{F}_B$-measurable point $x_0 \in B_{cc'}$ on $\{w_c \wedge w_{c'} \ne
0\}$ satisfying~\eqref{eq:model-boundary-null}. Let $Q = [Q_N \;\; Q_B] \in
\mathrm{O}(\dimin)$ be an $\mathcal{F}_B$-measurable orthogonal matrix whose
blocks $Q_N \in \RR^{\dimin \times 2}$ and $Q_B \in \RR^{\dimin \times (\dimin -
2)}$ span $N_B = \spn(w_c, w_{c'})$ and $T B_{cc'} = N_B^\perp$, so that $L_B =
Q_N T_B$ with $T_B := Q_N\transp L_B \in \mathrm{GL}_2(\RR)$. Then:
\begin{enumerate}[label=(\roman*)]
  \item Conditionally on $\mathcal{F}_B$, $m_B(x_0) \sim \mathrm{Bin}(\ngen - 2,
    1/2)$, and given additionally $\{m_B(x_0) = m\}$ ($0 \le m \le
    \ngen - 2$),
    \begin{equation}\label{eq:thinned-wishart}
      G_B(x_0) - I_{\dimin} \;\big|\; \{\mathcal{F}_B, m_B = m\}
      \sim \Wish_{\dimin}(m, \sigw^2 I_{\dimin}),
      \quad
      Q\transp G_B(x_0) Q
      = \left(\begin{smallmatrix}
          I_2 + \mathsf{W}_{NN} & \mathsf{W}_{NT} \\
          \mathsf{W}_{TN} & I_{\dimin - 2} + \mathsf{W}_{TT}
        \end{smallmatrix}\right),
    \end{equation}
    where $\bigl(\begin{smallmatrix} \mathsf{W}_{NN} & \mathsf{W}_{NT} \\
    \mathsf{W}_{TN} & \mathsf{W}_{TT} \end{smallmatrix}\bigr) \sim
    \Wish_{\dimin}(m, \sigw^2 I_{\dimin})$ (with $\Wish_{\dimin}(0, \sigw^2
    I_{\dimin}) := 0$).

    When $\dimin = 2$ the blocks $Q_B$, $\mathsf{W}_{NT}$
    and $\mathsf{W}_{TT}$ are empty, $Q = Q_N \in \mathrm{O}(2)$ and
    $Q\transp G_B(x_0) Q = I_2 + \mathsf{W}_{NN}$ with $\mathsf{W}_{NN} \sim
    \Wish_2(m, \sigw^2 I_2)$.
  \item For every realisation and all $m \ge 0$, the regularised transverse
    Schur complement
    \begin{equation}\label{eq:schur-normal}
      S_N(x_0)
      := (I_2 + \mathsf{W}_{NN}) - \mathsf{W}_{NT}(I_{\dimin - 2} + \mathsf{W}_{TT})^{-1} \mathsf{W}_{TN} \in \mathrm{Sym}_2^{++}
    \end{equation}
    satisfies $Q_N\transp G_B(x_0)^{-1} Q_N = S_N(x_0)^{-1}$ and $M_B(x_0) =
    T_B\transp S_N(x_0)^{-1} T_B$.
  \item For $m \ge \dimin$ ($\nu = m - \dimin + 2 \ge 2$), the unregularised
    Schur complement $S_N^{\mathrm{unreg}} := \mathsf{W}_{NN} - \mathsf{W}_{NT}
    \mathsf{W}_{TT}^{-1} \mathsf{W}_{TN}$ satisfies, conditionally on
    $\{\mathcal{F}_B, m_B(x_0) = m\}$~\cite{Muirhead1982},
    \begin{equation}\label{eq:schur-sandwich}
      S_N^{\mathrm{unreg}} \sim \Wish_2(\nu, \sigw^2 I_2)
      \,\indep\,
      \mathsf{W}_{TT} \sim \Wish_{\dimin - 2}(m, \sigw^2 I_{\dimin - 2}),
      \;
      I_2 + S_N^{\mathrm{unreg}} \preceq S_N(x_0) \preceq I_2 + \mathsf{W}_{NN}.
    \end{equation}
\end{enumerate}
\end{theorem}

\begin{proof}
Conditionally on $\mathcal{F}_B$ and $(I_h(x_0))_{h \notin B}$,
Lemma~\ref{lem:active-count}(ii) identifies $\{w_h w_h\transp : I_h(x_0) = 1\}$
as $m$ i.i.d.\ copies of $w_1 w_1\transp$, and orthogonal invariance under $Q$
gives~\eqref{eq:thinned-wishart}. Block inversion with $L_B = Q \binom{T_B}{0}$
gives~\eqref{eq:schur-normal}, while~\cite[Theorem~3.2.10]{Muirhead1982} and $0
\prec (I_{\dimin - 2} + \mathsf{W}_{TT})^{-1} \preceq \mathsf{W}_{TT}^{-1}$
give~\eqref{eq:schur-sandwich}.
\end{proof}

\subsection{Inverse-Wishart moment bounds}
\label{subsec:wishart-moments}

\begin{proposition}
\label{prop:inv-wishart-moments}
Under the hypotheses of Theorem~\ref{thm:conditional-wishart}:
\begin{enumerate}[label=(\roman*)]
  \item Although $G_B(x_0) - I_{\dimin}$ is singular whenever $m_B(x_0) <
    \dimin$ (an event of positive probability at every finite $\ngen$), every
    realisation satisfies
    \begin{equation}\label{eq:invw-barrier}
      G_B(x_0) \succeq I_{\dimin},
      \qquad
      \norm{G_B(x_0)^{-1}}_{\mathrm{op}} \le 1,
      \qquad
      \norm{S_N(x_0)^{-1}}_{\mathrm{op}} \le 1.
    \end{equation}
  \item For $m \ge \dimin$ ($\nu = m - \dimin + 2 \ge 2$), $M_B^{\mathrm{unreg}}
    := T_B\transp (S_N^{\mathrm{unreg}})^{-1} T_B \sim \Wish_2^{-1}(\nu,
    \sigw^{-2} L_B\transp L_B)$, where $\Wish_2^{-1}(\nu, \Psi)$ denotes the law
    of $S^{-1}$ for $S \sim \Wish_2(\nu, \Psi^{-1})$. We use the Bartlett
    decomposition~\cite[Theorem~3.2.14]{Muirhead1982} in its lower-triangular
    form, $\sigw^{-2}S_N^{\mathrm{unreg}} = A A\transp$ with $A =
    \bigl(\begin{smallmatrix} V_1 & 0 \\ Z & V_2 \end{smallmatrix}\bigr)$, where
    $V_1^2 \sim \chi_\nu^2$, $V_2^2 \sim \chi_{\nu - 1}^2$ and $Z \sim \Norm(0,
    1)$ are independent (the upper-triangular form $A\transp A$ would exchange
    the r\^oles of $V_1$ and $V_2$). For every $u \in \Sph{1}$,~\eqref{eq:schur-sandwich}
    and~\cite[Theorem~3.2.12]{Muirhead1982} give
    \begin{equation}\label{eq:invw-chi2}
      u\transp S_N(x_0)^{-1} u
      \le \min\Bigl(1,\, \frac{1}{\sigw^2 V_{m - \dimin + 1}}\Bigr),
      \qquad
      V_{m - \dimin + 1} := \frac{1}{\sigw^2 u\transp (S_N^{\mathrm{unreg}})^{-1} u}
      \sim \chi_{m - \dimin + 1}^2,
    \end{equation}
    where $\EE[V^{-p}] = 2^{-p} \Gamma(\frac{k}{2} - p) / \Gamma(\frac{k}{2})$
    for $V \sim \chi_k^2$ and $1 \le p < k / 2$.
  \item For every $p \ge 1$, there exists $C_{\dimin, p, \sigw} > 0$ such that
    for all $\ngen \ge 2$,
    \begin{equation}\label{eq:invw-moment-bound}
      \EE\Bigl[\norm{S_N(x_0)^{-1}}_{\mathrm{op}}^p \;\Big|\; \mathcal{F}_B\Bigr]
      \;\le\; C_{\dimin, p, \sigw}\, \ngen^{-p},
    \end{equation}
    and on $\{m_B(x_0) \ge \dimin\}$ we have $a_i(x_0) \le a_i^{\mathrm{unreg}}(x_0)
    := M_{B,ii}^{\mathrm{unreg}}(x_0)$ ($i = 1, 2$), satisfying $\EE[(a_1^{\mathrm{unreg}})^p
    \mid \mathcal{F}_B, m_B = m] \le C_{\dimin, p}\, R_c^{2p}\,\nu^{-p}$ and
    $\EE[(a_2^{\mathrm{unreg}})^p \mid \mathcal{F}_B, m_B = m] \le C_{\dimin, p}\,
    R_{c'}^{2p}\,\nu^{-p}$ for $m \ge \dimin + \lceil 2p \rceil$.
\end{enumerate}
\end{proposition}

\begin{proof}
Part~(i) is immediate from $G_B(x_0) \succeq I_{\dimin}$, and part~(ii) follows
from $S_N(x_0) \succeq S_N^{\mathrm{unreg}} \sim \Wish_2(\nu, \sigw^2 I_2)$
with~\cite[Theorems~3.2.12 and~3.2.14]{Muirhead1982}. For~\eqref{eq:invw-moment-bound},
bound $\norm{S_N(x_0)^{-1}}_{\mathrm{op}} \le e_1\transp S_N(x_0)^{-1} e_1 +
e_2\transp S_N(x_0)^{-1} e_2$ and split $m_B(x_0) \sim \mathrm{Bin}(\ngen - 2,
1/2)$ at $m_* := \max(\dimin + \lceil 2p \rceil, \lceil(\ngen - 2)/4\rceil)$.
For $m < m_*$, Hoeffding's inequality and~\eqref{eq:invw-barrier} give
$O_{\dimin,p}(\ngen^{-p})$, while for $m \ge m_*$,~\eqref{eq:invw-chi2} gives
$\EE[(u\transp S_N(x_0)^{-1} u)^p \mid \mathcal{F}_B, m_B = m] \le C_{\dimin, p,
\sigw}\ngen^{-p}$. Evaluating~\eqref{eq:invw-chi2} at $u = Q_N\transp
w_{c_i}/\norm{w_{c_i}}$ yields $a_i \le a_i^{\mathrm{unreg}} =
R_{c_i}^2/V_{m-\dimin+1}$ and its $\nu^{-p}$ moment bound.
\end{proof}

\subsection{The asymptotic curvature kernel \texorpdfstring{$\Zlim$}{Z\_infty}, second Hoeffding component and first-order annihilation}
\label{subsec:kernel-annihilation}

\begin{theorem}
\label{thm:fisher-limit-mark}
Assume the hypotheses of Proposition~\ref{prop:inv-wishart-moments}, and define
the \emph{asymptotic curvature kernel}
\begin{equation}\label{eq:limit-mark-def}
  \Zlim(w_c, w_{c'})
  := -\frac{1}{4\sigw^4}\, \ip{w_c}{w_{c'}}
     \sqrt{\norm{w_c}^2 \norm{w_{c'}}^2 - \ip{w_c}{w_{c'}}^2}
   = -\frac{R_c^2 R_{c'}^2}{4}\, \cos\vartheta_B \sin\vartheta_B.
\end{equation}
\begin{enumerate}[label=(\roman*)]
  \item For every $p \ge 1$, the dual overlap matrix $M_B(x_0) = L_B\transp
    G_B(x_0)^{-1} L_B$ satisfies
    \begin{equation}\label{eq:compressed-conv}
      \frac{\ngen \sigw^2}{2}\, M_B(x_0)
      \;\xrightarrow[\ngen \to \infty]{L^p(\PP)}\;
      L_B\transp L_B
      = \begin{pmatrix} \norm{w_c}^2 & \ip{w_c}{w_{c'}} \\ \ip{w_c}{w_{c'}} & \norm{w_{c'}}^2 \end{pmatrix}.
    \end{equation}
  \item For every $p \ge 1$, there exists $C_{\dimin, p, \sigw} > 0$ depending
    only on $(\dimin, p, \sigw)$ such that for all $\ngen \ge 2$,
    \begin{equation}\label{eq:mark-Lp-rate}
      \begin{aligned}
        \norm[\Big]{\frac{\ngen^2}{4}\,\varepsilon_B(x_0) - \Zlim(w_c, w_{c'})}_{L^p(\PP)}
        &\le C_{\dimin, p, \sigw}\, \ngen^{-1/2}, \\
        \norm[\Big]{\EE\!\left[\frac{\ngen^2}{4}\,\varepsilon_B(x_0) \;\Big|\; \mathcal{F}_B\right] - \Zlim(w_c, w_{c'})}_{L^p(\PP)}
        &\le C_{\dimin, p, \sigw}\, \ngen^{-1},
      \end{aligned}
    \end{equation}
    and for $\nu = m_B(x_0) - \dimin + 2 > 4$, $\EE[\varepsilon_0(M_B^{\mathrm{unreg}}(x_0))
    \mid \mathcal{F}_B, m_B = m] = \Zlim(w_c, w_{c'})/[(\nu - 2)(\nu - 4)]$.
\end{enumerate}
\end{theorem}

\begin{proof}
Part~(i) is proved in Appendix~\ref{sec:app-wishart} via Marcinkiewicz--Zygmund
concentration of $\frac{2}{\ngen\sigw^2}S_N(x_0) - I_2$
and~\eqref{eq:invw-moment-bound}. For part~(ii), since $\abs{\varepsilon_B} \le
\pi/2$ and $\abs{\Zlim(w_c, w_{c'})} \le \frac{1}{4}R_c^2 R_{c'}^2 \in L^p(\PP)$
(Proposition~\ref{prop:limit-mark-moments}), bounded values of $\ngen$ are
absorbed into $C_{\dimin, p, \sigw}$. Take $\ngen$ sufficiently large depending
on $(\dimin, p, \sigw)$, and let $\mathcal{G}_{\mathrm{bin}} :=
\{\abs{m_B(x_0) - \frac{\ngen - 2}{2}} \le \frac{\ngen - 2}{4}\}$. By
Hoeffding's inequality, $\PP(\mathcal{G}_{\mathrm{bin}}^c \mid \mathcal{F}_B)
\le 2\exp(-(\ngen - 2)/8)$, contributing $O_{\dimin, p}(\ngen^{-1})$ to both
lines of~\eqref{eq:mark-Lp-rate}.

On $\mathcal{G}_{\mathrm{bin}}$, $m = m_B(x_0) \ge \dimin + \lceil 8p \rceil$
and $\nu = m - \dimin + 2 \ge c_0\ngen \ge \nu_0(p)$. Writing $\PP_{\mathcal{F}_B,
m} := \PP(\cdot \mid \mathcal{F}_B, m_B = m)$, Lemmas~\ref{lem:app-bartlett}
and~\ref{lem:app-schur-shift} in Appendix~\ref{sec:app-wishart} yield
\begin{align*}
  \norm{\nu^2(\varepsilon_B(x_0) - \varepsilon_0(x_0))}_{L^p(\PP_{\mathcal{F}_B, m})}
  &\le C_{\dimin, p} R_c^2 R_{c'}^2(1 + R_c^2 + R_{c'}^2)\nu^{-1}, \\
  \norm{\nu^2\varepsilon_0(M_B(x_0)) -
  \nu^2\varepsilon_0(M_B^{\mathrm{unreg}}(x_0))}_{L^p(\PP_{\mathcal{F}_B, m})}
  &\le C_{\dimin, p, \sigw} R_c^2 R_{c'}^2 \nu^{-1}.
\end{align*}
By Lemma~\ref{lem:app-bartlett}(i),
\begin{align*}
  \nu^2\varepsilon_0(M_B^{\mathrm{unreg}}) - \Zlim(w_c, w_{c'})
  &= -\frac{1}{4}R_c^2 R_{c'}^2\Bigl[\cos\vartheta_B\sin\vartheta_B
      \bigl(\nu^2 V_1^{-3}V_2^{-1} - 1 + Z^2\nu^2 V_1^{-3}V_2^{-3}\bigr) \\
  &\qquad\qquad\qquad\quad
      - \sin^2\vartheta_B\,Z\nu^2 V_1^{-2}V_2^{-3}\Bigr] \\
  &= O_{L^p(\PP_{\mathcal{F}_B, m})}\bigl(R_c^2 R_{c'}^2 \nu^{-1/2}\bigr).
\end{align*}
Since
\[
  \norm{\ngen^2/(4\nu^2) - 1}_{L^p(\mathcal{G}_{\mathrm{bin}})}
  = O_{\dimin}(\ngen^{-1/2}),
\]
averaging over $m \sim \mathrm{Bin}(\ngen - 2, 1/2)$
and $\mathcal{F}_B$ proves the first bound in~\eqref{eq:mark-Lp-rate}.

For the conditional mean ($\nu > 4$), the Gamma identities
$\EE[V_1^{-3}]\EE[V_2^{-1}] = [(\nu-2)(\nu-3)]^{-1}$ and $\EE[V_2^{-3}] =
(\nu-4)^{-1}\EE[V_2^{-1}]$, together with $\EE[Z] = 0$ and $\EE[Z^2] = 1$, give
$\EE[\varepsilon_0(M_B^{\mathrm{unreg}}(x_0)) \mid \mathcal{F}_B, m_B = m] =
\Zlim(w_c, w_{c'})/[(\nu - 2)(\nu - 4)]$. Combining this with the two $O(\nu^{-1})$
remainder bounds (at $p = 1$) gives, on $\mathcal{G}_{\mathrm{bin}}$,
\begin{equation}\label{eq:kernel-cond-nu}
  \EE\bigl[\varepsilon_B(x_0) \;\big|\; \mathcal{F}_B, m_B = m\bigr]
  = \frac{\Zlim(w_c, w_{c'})}{(\nu - 2)(\nu - 4)}
    + O_{\dimin,\sigw}\!\bigl(R_c^2R_{c'}^2(1 + R_c^2 + R_{c'}^2)\,\nu^{-3}\bigr).
\end{equation}
It remains to average~\eqref{eq:kernel-cond-nu} over $m_B(x_0) \sim
\mathrm{Bin}(\ngen - 2, 1/2)$. Set $f(\nu) := \frac{(\ngen/2)^2}{(\nu - 2)(\nu -
4)}$ and write $\nu = \bar{\nu} + \delta$ with $\bar{\nu} := \frac{\ngen}{2} -
\dimin + 1$ and $\delta := m_B(x_0) - \frac{\ngen - 2}{2}$, so that $\EE[\delta
\mid \mathcal{F}_B] = 0$ and $\EE[\delta^2 \mid \mathcal{F}_B] = \frac{\ngen -
2}{4}$. On $\mathcal{G}_{\mathrm{bin}}$ we have $\nu \ge c_0\ngen$, and Taylor's
theorem gives
\[
  f(\nu) - 1
  = \bigl(f(\bar{\nu}) - 1\bigr) + f'(\bar{\nu})\,\delta
    + \tfrac{1}{2}f''(\tilde{\nu})\,\delta^2,
\]
for some $\tilde{\nu}$ between $\bar{\nu}$ and $\nu$, where $\abs{f(\bar{\nu})
- 1} \le C_{\dimin}/\ngen$, $\abs{f'(\bar{\nu})} \le C/\ngen$ and
$\abs{f''(\tilde{\nu})} \le C/\ngen^2$. The constant term is
$O(\ngen^{-1})$ because $(\bar{\nu} - 2)(\bar{\nu} - 4) = (\ngen/2)^2 -
O_{\dimin}(\ngen)$. The linear term has conditional mean $f'(\bar{\nu})\EE[\delta
\mathbf{1}_{\mathcal{G}_{\mathrm{bin}}} \mid \mathcal{F}_B] =
-f'(\bar{\nu})\EE[\delta \mathbf{1}_{\mathcal{G}_{\mathrm{bin}}^c} \mid
\mathcal{F}_B] = O(e^{-c'\ngen})$, since $\abs{\delta} \le \ngen$ and
$\PP(\mathcal{G}_{\mathrm{bin}}^c \mid \mathcal{F}_B) \le 2e^{-(\ngen -
2)/8}$. The quadratic term contributes at most $\frac{C}{\ngen^2}\EE[\delta^2
\mid \mathcal{F}_B] = O(\ngen^{-1})$. Hence $\EE[(f(\nu) -
1)\mathbf{1}_{\mathcal{G}_{\mathrm{bin}}} \mid \mathcal{F}_B] =
O_{\dimin}(\ngen^{-1})$. Combining this with the $O(\nu^{-3})$ remainder
in~\eqref{eq:kernel-cond-nu} (which contributes $O(\ngen^{-1})$ after
multiplication by $(\ngen/2)^2$), with the $\mathcal{G}_{\mathrm{bin}}^c$
tail, and absorbing bounded $\ngen \le H_0(\dimin, \sigw)$ via $1 \le
H_0\ngen^{-1}$, we obtain, for all linearly independent $(w_c, w_{c'})$ and
$\ngen \ge 2$,
\[
  \abs[\Big]{\EE\!\left[\frac{\ngen^2}{4}\,\varepsilon_B(x_0) \;\Big|\; \mathcal{F}_B\right] - \Zlim(w_c, w_{c'})}
  \le C_{\dimin,\sigw}\bigl(1 + R_c^2 R_{c'}^2(1 + R_c^2 + R_{c'}^2)\bigr)\,\ngen^{-1},
\]
whose $L^p(\PP)$ norm proves the second bound in~\eqref{eq:mark-Lp-rate}.
\end{proof}

\begin{proposition}
\label{prop:limit-mark-moments}
Under~\ref{assump:A2} ($\dimin \ge 2$), $R_c^2, R_{c'}^2
\stackrel{\mathrm{i.i.d.}}{\sim} \chi_{\dimin}^2$ and $\vartheta_B \in (0, \pi)$
are mutually independent, and $\vartheta_B$ has density $f_\vartheta(\vartheta) =
\wallis{\dimin - 2}^{-1} \sin^{\dimin - 2}\vartheta$ with $\wallis{k} :=
\int_0^\pi \sin^k\vartheta\,\dd\vartheta =
\sqrt{\pi}\,\Gamma(\frac{k+1}{2})/\Gamma(\frac{k+2}{2})$. \textup{(i)} For every
integer $q \ge 1$, $\EE[\Zlim(w_c, w_{c'})] = 0$ and
\begin{equation}\label{eq:limit-mark-moments}
  \EE\abs{\Zlim} = \frac{\dimin}{2\wallis{\dimin - 2}},
  \quad
  \EE[\Zlim^2] = \frac{\dimin(\dimin + 2)(\dimin - 1)}{16},
  \quad
  \EE\bigl[\abs{\Zlim}^q\bigr]
  = \frac{\Gamma(\frac{\dimin}{2} + q)\,\Gamma(\frac{q + 1}{2})\,\Gamma(\frac{\dimin - 1 + q}{2})}
          {\sqrt{\pi}\,\Gamma(\frac{\dimin}{2})\,\Gamma(\frac{\dimin - 1}{2})}.
\end{equation}
\textup{(ii)} For any ball $\Omega = B(x_*, r) \subset \RR^{\dimin}$ ($r > 0$),
every non-empty $(\dimin - 2)$-flat slice is a $(\dimin - 2)$-ball of radius at
most $r$, so $\beta_{B(x_*, r)} = \omega_{\dimin - 2} r^{\dimin - 2} < \infty$.
\end{proposition}

\begin{proof}
Antisymmetry of $\cos\vartheta$ about $\pi/2$ gives $\EE[\Zlim] = 0$. Substituting
$t = \sin^2\vartheta$ evaluates $\EE[\abs{\cos\vartheta_B}^q \sin^q\vartheta_B] =
\Gamma(\frac{\dimin - 1 + q}{2})\Gamma(\frac{q + 1}{2})/[\wallis{\dimin -
2}\Gamma(\frac{\dimin}{2} + q)]$, and multiplying by $4^{-q}(\EE[R_c^{2q}])^2$
with $\EE[R_c^{2q}] = 2^q\Gamma(\frac{\dimin}{2} + q)/\Gamma(\frac{\dimin}{2})$
gives~\eqref{eq:limit-mark-moments}. Orthogonal projection onto $F$ gives~(ii).
\end{proof}

\begin{lemma}
\label{lem:hoeffding-support}
Under~\ref{assump:A1}, write $\xi_S := (\xi_h)_{h \in S}$
for $S \subseteq [\ngen]$. For a Borel window $\Omega \subseteq \RR^{\dimin}$
and $B = \{c, c'\} \in \binom{[\ngen]}{2}$, let $Y_B(\Omega) = Y_{cc'}(\Omega) =
\int_{B_{cc'} \cap \Omega} \varepsilon_{cc'}(x)\,\dd\Haus{\dimin - 2}(x)$ be
the bone curvature functional~\eqref{eq:model-measures} on any window with
$Y_B(\Omega) \in L^2(\PP)$. Then for every $S \subseteq [\ngen]$
with $B \not\subseteq S$,
\begin{equation}\label{eq:hoeffding-cond-zero}
  \EE\bigl[Y_B(\Omega) \;\big|\; \xi_S\bigr] = 0
  \qquad \text{almost surely.}
\end{equation}
Consequently $\pi_A Y_B(\Omega) = 0$ a.s.\ whenever $B \not\subseteq A$,
$\EE[Y_B(\Omega)] = P_1 Y_B(\Omega) = 0$, and for $A \supseteq B$,
\begin{equation}\label{eq:hoeffding-P2-collapse}
  \pi_A Y_B(\Omega)
  = \sum_{B \subseteq S \subseteq A} (-1)^{\abs{A \setminus S}}\,
    \EE\bigl[Y_B(\Omega) \big| \xi_S\bigr],
  \qquad
  P_2 Y_B(\Omega)
  = \pi_B Y_B(\Omega)
  = \EE\bigl[Y_B(\Omega) \big| \mathcal{F}_B\bigr].
\end{equation}
\end{lemma}

\begin{proof}
For $j \in B \setminus S$, Lemma~\ref{lem:sign-flip-oddness} shows that
flipping $\xi_j \mapsto -\xi_j$ negates $Y_B(\Omega)$ pathwise while fixing
$\xi_S$, so~\ref{assump:A1} gives $\EE[Y_B(\Omega) \mid \xi_S] = 0$ a.s.,
yielding~\eqref{eq:hoeffding-cond-zero}--\eqref{eq:hoeffding-P2-collapse}.
\end{proof}

\begin{theorem}
\label{thm:second-chaos-ustat}
Let $\Omega \subseteq \RR^{\dimin}$ be a Borel window with $\beta_\Omega <
\infty$, and suppose~\ref{assump:A2} holds on $\Omega$ with symmetric i.i.d.\
offsets. For linearly independent $L_B = [w_c \;\; w_{c'}] \in \RR^{\dimin
\times 2}$ ($T_B\transp T_B = L_B\transp L_B$), define the fixed-$m$ conditional
mean and \emph{thinned-Wishart kernel}
\begin{equation}\label{eq:bone-kernel-def}
  \begin{aligned}
    \bar{\varepsilon}_m(L_B)
    &:= \EE\bigl[\Psi_{\mathrm{mat}}(T_B\transp S_N(x_0)^{-1} T_B) \mid \mathcal{F}_B, m_B = m\bigr], \\
    e_{\ngen}(w_c, w_{c'})
    &:= 2^{-(\ngen - 2)} \sum_{m=0}^{\ngen - 2} \binom{\ngen - 2}{m}\bar{\varepsilon}_m(L_B),
  \end{aligned}
\end{equation}
which do not depend on $x_0 \in B_{cc'}$ by
Theorem~\ref{thm:conditional-wishart}. Then \textup{(i)} $\EE[\varepsilon_B(x)
\mid \mathcal{F}_B] = e_{\ngen}(w_c, w_{c'})$ for $\Haus{\dimin - 2}$-a.e.\ $x
\in B_{cc'}$ and $\EE[Y_{cc'}(\Omega) \mid \mathcal{F}_B] = e_{\ngen}(w_c,
w_{c'})\,\ell_{cc'}(\Omega)$; and \textup{(ii)} $P_2(\Regge(\Omega))$ is a
degenerate order-$2$ U-statistic with pairwise orthogonal summands in $L^2(\PP)$,
\begin{equation}\label{eq:second-chaos-ustat}
  \begin{gathered}
    P_2(\Regge(\Omega))
    = \sum_{1 \le c < c' \le \ngen} e_{\ngen}(w_c, w_{c'})\, \ell_{cc'}(\Omega), \\
    \Var\bigl(P_2(\Regge(\Omega))\bigr)
    = \binom{\ngen}{2}\EE\bigl[e_{\ngen}(w_1, w_2)^2 \ell_{12}(\Omega)^2\bigr].
  \end{gathered}
\end{equation}
\end{theorem}

\begin{proof}
Parametrise $B_{cc'}$ $\mathcal{F}_B$-measurably by $x(y) := x_0^{\min}(B) + Q_B
y$, $y \in \RR^{\dimin - 2}$ ($y = 0$ at $\dimin = 2$, where $\PP(b_c = b_{c'} =
0) = 0$). For every $y \ne 0$ (and $y = 0$ at $\dimin = 2$), $\PP(x(y) = 0) = 0$,
so~\eqref{eq:model-boundary-null} holds at $x(y)$. By Fubini's theorem on the
jointly measurable indicators (Remark~\ref{rem:bone-geometry}),
$\Haus{\dimin - 2}(B_{cc'} \setminus B_{cc'}^{\circ}) = 0$ a.s., and for
$\mathrm{Leb}_{\dimin-2}$-a.e.\ $y$ Theorem~\ref{thm:conditional-wishart} shows
that the conditional law of $M_B(x(y))$ given $\mathcal{F}_B$ depends only on
$L_B$. Averaging over $m_B \mid \mathcal{F}_B \sim \mathrm{Bin}(\ngen - 2, 1/2)$
and integrating over $y$ gives~(i), and~\eqref{eq:hoeffding-P2-collapse} with
Hoeffding orthogonality gives~\eqref{eq:second-chaos-ustat}.
\end{proof}

\begin{proposition}
\label{thm:first-chaos-annihilation}
Under~\ref{assump:A1}, let $\Omega \subseteq
\RR^{\dimin}$ be a Borel window with $\beta_\Omega < \infty$ (more generally,
$\EE[\ell_B(\Omega)] < \infty$ in~\textup{(i)} and $\EE[\ell_B(\Omega)^2] <
\infty$ in~\textup{(ii)} for all $B$). \textup{(i)} For every $j \in [\ngen]$,
$\EE[\Regge(\Omega) \mid w_j, b_j] = 0$ a.s.\ and
\begin{equation}\label{eq:first-chaos-zero}
  \EE\bigl[\Regge(\Omega)\bigr] = 0,
  \qquad
  P_1\bigl(\Regge(\Omega)\bigr) = 0
  \qquad \text{almost surely.}
\end{equation}
\textup{(ii)} In $L^2(\PP)$, $\Regge(\Omega) = P_2(\Regge(\Omega)) + R_{\ge
3}(\Omega)$ with $R_{\ge 3}(\Omega) := \sum_{k=3}^{\ngen} P_k(\Regge(\Omega))$
and
\begin{equation}\label{eq:regge-chaos-split}
  \Var\bigl(\Regge(\Omega)\bigr)
  = \Var\Bigl(P_2\bigl(\Regge(\Omega)\bigr)\Bigr) + \Var\bigl(R_{\ge 3}(\Omega)\bigr).
\end{equation}
Both identities also hold for the homogeneous planar apex
measure~\eqref{eq:model-apex-convention} whenever the rows of $W$ are a.s.\
pairwise non-collinear (Corollary~\ref{cor:planar-moments}\textup{(i)}).
\end{proposition}

\begin{proof}
Since $\abs{Y_B(\Omega)} \le \frac{\pi}{2}\ell_B(\Omega)$,
Lemma~\ref{lem:hoeffding-support}~\eqref{eq:hoeffding-cond-zero} gives
$\EE[Y_B(\Omega) \mid \xi_S] = 0$ a.s.\ for all $B$ and $\abs{S} \le 1$. Thus
summing over bones yields~\eqref{eq:first-chaos-zero}--\eqref{eq:regge-chaos-split}
via~\eqref{eq:model-P-k} without any general-position assumption.
\end{proof}
\section{Quantitative second-order transfer theorem}
\label{sec:transfer}

Throughout this section we fix $\dimin \ge 2$, take $\ngen \ge \dimin + 2 \ge
4$, and assume that~\ref{assump:A2} holds on a Borel window $\Omega \subseteq
\RR^{\dimin}$ ($\beta_\Omega < \infty$) with symmetric i.i.d.\ offsets. By
Proposition~\ref{thm:first-chaos-annihilation}~\eqref{eq:regge-chaos-split},
$\Regge(\Omega) = P_2(\Regge(\Omega)) + R_{\ge 3}(\Omega)$ with $R_{\ge
3}(\Omega) = \sum_{B \in \binom{[\ngen]}{2}} \widetilde{Y}_B(\Omega)$ and
$\widetilde{Y}_B(\Omega) := Y_B(\Omega) - \EE[Y_B(\Omega) \mid \mathcal{F}_B]$.

\subsection{Cross-bone Efron--Stein bound}
\label{subsec:cross-bone-es}

Let $\xi_1, \dots, \xi_{\ngen}$ be independent, and let $X \in L^2(\sigma(\xi_1,
\dots, \xi_{\ngen}))$ have Hoeffding decomposition
$X = \sum_{A \subseteq [\ngen]} \pi_A X$~\eqref{eq:model-P-k}. Set $\delta_j X := X - \EE[X \mid
\mathcal{F}_{[\ngen]\setminus\{j\}}]$ and $\Delta_j X := X - X^{(j)}$, where
$X^{(j)}$ replaces $\xi_j$ by an independent copy $\xi_j'$. For $j \ne j'$, set
$\Delta_j \Delta_{j'} X := X - X^{(j)} - X^{(j')} + X^{(j,j')}$.

\begin{lemma}
\label{lem:master-es}
Let $\xi_1, \dots, \xi_{\ngen}$ be independent random variables ($\ngen \ge 4$).
\begin{enumerate}[label=(\roman*)]
  \item For $X \in L^2(\PP)$ and $j \ne j'$, $\delta_j$ and $\delta_{j'}$ are
    commuting orthogonal projections with $\delta_j X = \sum_{A \ni j} \pi_A X$,
    $\delta_j \delta_{j'} X = \sum_{A \supseteq \{j, j'\}} \pi_A X$,
    $\EE[(\delta_j X)^2] = \frac{1}{2}\EE[(\Delta_j X)^2]$ and
    $\EE[(\delta_j \delta_{j'} X)^2] = \frac{1}{4}\EE[(\Delta_j \Delta_{j'}
    X)^2]$~\cite{Hoeffding1948,EfronStein1981,BLM2013}.
  \item Let $(Y_B)_{B \in \binom{[\ngen]}{2}}$ satisfy $\pi_A Y_B = 0$ whenever
    $B \not\subseteq A$, and set $\widetilde{Y}_B := Y_B - \pi_B Y_B = \sum_{A
    \supsetneq B} \pi_A Y_B$. Then for any $B, B' \in \binom{[\ngen]}{2}$
    (writing $B' = \{d_1, d_2\}$ in the disjoint case):
    \begin{align}
      \Var(\widetilde{Y}_B)
      &= \sum_{A \supsetneq B} \EE[(\pi_A Y_B)^2]
      \le \frac{1}{2}\sum_{j \in [\ngen]\setminus B} \EE[(\Delta_j Y_B)^2]
      && (B = B'),
      \label{eq:es-diag} \\
      \abs{\Cov(\widetilde{Y}_B, \widetilde{Y}_{B'})}
      &\le \norm{\delta_{c_2} Y_B}_{L^2} \norm{\delta_{c_1} Y_{B'}}_{L^2}
      && (B = \{c, c_1\},\, B' = \{c, c_2\}),
      \label{eq:es-adj-cov} \\
      \abs{\Cov(\widetilde{Y}_B, \widetilde{Y}_{B'})}
      &\le \norm{\delta_{d_1} \delta_{d_2} Y_B}_{L^2}\norm{\delta_{c_1} \delta_{c_2} Y_{B'}}_{L^2}
      && (B \cap B' = \emptyset).
      \label{eq:es-disj-cov}
    \end{align}
  \item If $\xi_1, \dots, \xi_{\ngen}$ are i.i.d.\ and $Y_{\{c,c'\}}$ is
    exchangeable in $\{c, c'\}$, then $R_{\ge 3} := \sum_{B \in
    \binom{[\ngen]}{2}} \widetilde{Y}_B$ satisfies
    \begin{align}
      \Var(R_{\ge 3})
      &\le 9\binom{\ngen}{3} \EE\bigl[\bigl(\delta_3 Y_{\{1,2\}}\bigr)^2\bigr]
          + 6\binom{\ngen}{4} \EE\bigl[\bigl(\delta_3 \delta_4 Y_{\{1,2\}}\bigr)^2\bigr]
      \label{eq:es-master-D} \\
      &= \frac{9}{2}\binom{\ngen}{3} \EE\bigl[\bigl(\Delta_3 Y_{\{1,2\}}\bigr)^2\bigr]
        + \frac{3}{2}\binom{\ngen}{4} \EE\bigl[\bigl(\Delta_3 \Delta_4 Y_{\{1,2\}}\bigr)^2\bigr].
      \label{eq:es-master-Delta}
    \end{align}
\end{enumerate}
\end{lemma}

\begin{proof}
Part~(i) and~\eqref{eq:es-diag} follow directly from Hoeffding orthogonality
$\mathcal{V}_A \perp \mathcal{V}_{A'}$ ($A \ne A'$) and the standard
Efron--Stein jackknife identities~\cite{Hoeffding1948,EfronStein1981,BLM2013}.
For $B \ne B'$, $\Cov(\widetilde{Y}_B, \widetilde{Y}_{B'}) = \sum_{A \supseteq B
\cup B'} \EE[(\pi_A Y_B)(\pi_A Y_{B'})]$. Since $\sum_{A \supseteq
\{c,c_1,c_2\}} \pi_A Y_{\{c,c_1\}} = \delta_{c_2} Y_{\{c,c_1\}}$ and $\sum_{A
\supseteq B \cup \{d_1,d_2\}} \pi_A Y_B = \delta_{d_1}\delta_{d_2} Y_B$,
Cauchy--Schwarz gives~\eqref{eq:es-adj-cov}--\eqref{eq:es-disj-cov}. Summing
over the $\binom{\ngen}{2}(\ngen - 2) = 3\binom{\ngen}{3}$ diagonal terms,
$6\binom{\ngen}{3}$ ordered adjacent pairs and $6\binom{\ngen}{4}$ ordered
disjoint pairs gives~\eqref{eq:es-master-D}--\eqref{eq:es-master-Delta}.
\end{proof}

\begin{lemma}
\label{lem:cross-bone-es}
Under the standing hypotheses of this section, there exists a constant
$C_{\dimin} > 0$ depending only on $\dimin$ such that for all $\ngen \ge \dimin
+ 2$ and $\sigw > 0$,
\begin{equation}\label{eq:es-replacement-bounds}
  \EE\Bigl[\bigl(\Delta_3 Y_{\{1,2\}}(\Omega)\bigr)^2\Bigr]
  \le C_{\dimin}\, \beta_\Omega^2\, \ngen^{-6},
  \qquad
  \EE\Bigl[\bigl(\Delta_3 \Delta_4 Y_{\{1,2\}}(\Omega)\bigr)^2\Bigr]
  \le C_{\dimin}\, \beta_\Omega^2\, \ngen^{-8}.
\end{equation}
Thus the higher-order Hoeffding remainder satisfies
\begin{equation}\label{eq:es-variance-bound}
  \Var\bigl(R_{\ge 3}(\Omega)\bigr)
  \le C_{\dimin}\, \beta_\Omega^2\, \ngen^{-3},
  \qquad
  \norm[\Big]{\frac{\ngen}{4} R_{\ge 3}(\Omega)}_{L^2(\PP)}
  \le C_{\dimin}\, \beta_\Omega\, \ngen^{-1/2}.
\end{equation}
\end{lemma}

\begin{proof}
The basic method is the following. We parametrise the bone $B_{12}$. At
each point $x$ on this bone we compare the deficit $\varepsilon_{12}(x)$ before and
after replacing the generators $\xi_3$ and $\xi_4$. On an event where the
background metric is well conditioned, the perturbation bounds of
Lemma~\ref{lem:deficit-perturbation} control these differences by inverse
powers of $\ngen$. Off that event we use crude bounds and the exponentially
small probability of the event's complement. Integrating along the bone then
gives~\eqref{eq:es-replacement-bounds}.

Since $\abs{\varepsilon_B(x)} \le \pi/2$ bounds $\abs{\Delta_3 Y_B(\Omega)} \le
\pi\beta_\Omega$ and $\abs{\Delta_3\Delta_4 Y_B(\Omega)} \le 2\pi\beta_\Omega$,
values of $\ngen/\dimin$ below any universal constant are absorbed into
$C_{\dimin}$. Assume $\ngen/\dimin$ is sufficiently large, fix $B = \{1, 2\}$,
and parametrise $B_{12}$ $\mathcal{F}_{\{1,2\}}$-measurably by $x(y) :=
x_0^{\min}(B_{12}) + Q_{B_{12}}y$, $y \in \RR^{\dimin - 2}$
(Remark~\ref{rem:bone-geometry}, with $y = 0$ at $\dimin = 2$). We work
throughout at points $x = x(y)$ with $x \ne 0$. For $\dimin = 2$ this holds
almost surely, since $\PP(b_1 = b_2 = 0) = 0$. For $\dimin \ge 3$ it excludes at
most the single parameter value $y = 0$ (and then only when $b_1 = b_2 = 0$).
This point is $\mathrm{Leb}_{\dimin-2}$-null, so it does not affect the integrals
below. The exclusion is genuinely needed, because at $x = 0$ the indicators
$I_h(0) = \mathbf{1}\{b_h > 0\}$ are no longer fair coin flips when $P_b$ has an
atom at $0$. At every such $x$, $\PP(x = 0) = 0$ and
so~\eqref{eq:model-boundary-null} holds for every $h \notin \{1, 2\}$ and for
the independent copies $\xi_3', \xi_4'$. Set $G_0(x) := I_{\dimin} +
\sum_{h=5}^{\ngen} I_h(x) w_h w_h\transp$ and $m_0(x) := \sum_{h=5}^{\ngen}
I_h(x)$, and write $R_c := \norm{w_c}/\sigw$.

We first construct the good event. By Lemma~\ref{lem:active-count},
conditionally on $\mathcal{F}_{\{1,2\}}$ we have $m_0(x) \sim \mathrm{Bin}(\ngen
- 4, 1/2)$ and, given $m_0(x) = m$, $\widetilde{\mathsf{W}}_0(x) :=
\sum_{h=5}^{\ngen} I_h(x)(w_h/\sigw)(w_h/\sigw)\transp = Z Z\transp$ for a
$\dimin \times m$ matrix $Z$ with i.i.d.\ standard Gaussian entries, with
$\sigw^{-2} G_0(x) \succeq \widetilde{\mathsf{W}}_0(x)$. By Hoeffding's
inequality, $\PP(m_0(x) < (\ngen - 4)/4 \mid \mathcal{F}_{\{1,2\}}) \le
e^{-(\ngen - 4)/8}$. On $\{m_0(x) = m \ge (\ngen - 4)/4\}$, the
Davidson--Szarek inequality~\cite[Theorem~II.13]{DavidsonSzarek2001}
(cf.~\cite[Theorem~4.6.1]{Vershynin2018}) states that the extreme singular
values of $Z$ satisfy $\sqrt{m} - \sqrt{\dimin} - t \le s_{\min}(Z) \le
s_{\max}(Z) \le \sqrt{m} + \sqrt{\dimin} + t$ with probability at least $1 -
2e^{-t^2/2}$. Taking $t = \sqrt{m}/4$ and using $m \ge (\ngen - 4)/4 \ge
16\dimin$ for $\ngen/\dimin$ large, we obtain universal constants $c_0, c_1 > 0$
and $C_0 > 1$ such that the event
\[
  \mathcal{E}_{\mathrm{DS}}(x)
  := \bigl\{c_0\ngen \le \lambda_{\min}(\widetilde{\mathsf{W}}_0(x))
     \le \lambda_{\max}(\widetilde{\mathsf{W}}_0(x)) \le C_0\ngen\bigr\}
\]
satisfies $\PP(\mathcal{E}_{\mathrm{DS}}(x)^c \mid \mathcal{F}_{\{1,2\}}) \le
Ce^{-c_1\ngen}$, uniformly in $y$. Choosing a small universal constant
$c \in (0, c_0/6]$, define
\begin{equation}\label{eq:es-good-event}
  \mathcal{A}_B(x)
  := \mathcal{E}_{\mathrm{DS}}(x)
     \cap \bigl\{ R_1^2 + R_2^2 \le c\,\ngen \bigr\}
     \cap \bigcap_{j \in \{3, 4\}} \bigl\{ R_j^2 + R_j'^2 \le c\,\ngen \bigr\}.
\end{equation}

Next we bound the replacement differences on $\mathcal{A}_B(x)$. Around the
base matrix $G_B(x) := G_0(x) + I_3(x)w_3 w_3\transp + I_4(x)w_4 w_4\transp
\succeq G_0(x)$, the symmetric replacement perturbations $E_1 :=
I_3'(x)w_3'(w_3')\transp - I_3(x)w_3 w_3\transp$ and $E_2 :=
I_4'(x)w_4'(w_4')\transp - I_4(x)w_4 w_4\transp$ satisfy $\Delta_3\varepsilon_{12}(x)
= \varepsilon_{12}(G_B) - \varepsilon_{12}(G_B + E_1)$ and
$\Delta_3\Delta_4\varepsilon_{12}(x) = \varepsilon_{12}(G_B + E_1 + E_2)
- \varepsilon_{12}(G_B + E_1) - \varepsilon_{12}(G_B + E_2) +
\varepsilon_{12}(G_B)$. Here $\varepsilon_{12}(x)$ is evaluated at the point
$x$. This point lies on $B_{12}^\circ$ for the original and the replaced
arrangements alike, since $\PP(x \in \Pi_h) = 0$ for $h \ge 3$ and for the
copies $\xi_3', \xi_4'$ by~\eqref{eq:model-boundary-null}. For $\dimin \ge 3$
this also means that $B_{12} \cap \Pi_3$, a codimension-one affine subspace of
$B_{12}$, is $\Haus{\dimin-2}$-null. This fact justifies the use of Fubini's theorem
below. On $\mathcal{A}_B(x)$, $G_B^{-1} \preceq G_0^{-1}$ gives
$\tr(L_B\transp G_B^{-1} L_B) \le (R_1^2 + R_2^2)/(c_0\ngen) \le 1/2$,
$(w_1\transp G_B^{-1} w_1)(w_2\transp G_B^{-1} w_2) \le C\ngen^{-2} R_1^2 R_2^2$,
and $\norm{G_B^{-1/2}E_i G_B^{-1/2}}_{\mathrm{op}} \le (R_{i+2}^2 +
R_{i+2}'^2)/(c_0\ngen)$ with sum at most $1/3$.
Lemma~\ref{lem:deficit-perturbation}(ii)~\eqref{eq:pert-diff-1}
and~\eqref{eq:pert-diff-2} therefore give
\begin{equation}\label{eq:es-good-pointwise}
  \begin{gathered}
  \abs{\Delta_3 \varepsilon_{12}(x)}\,\mathbf{1}_{\mathcal{A}_B(x)}
  \le C\, \ngen^{-3}\, R_1^2 R_2^2\, (R_3^2 + R_3'^2), \\
  \abs{\Delta_3 \Delta_4 \varepsilon_{12}(x)}\,\mathbf{1}_{\mathcal{A}_B(x)}
  \le C\, \ngen^{-4}\, R_1^2 R_2^2\prod_{j=3}^4 (R_j^2 + R_j'^2).
  \end{gathered}
\end{equation}

Finally we integrate along the bone. By Cauchy--Schwarz over the slice $\{y \in
\RR^{\dimin-2} : x(y) \in \Omega\}$, which has measure $\ell_{12}(\Omega) \le
\beta_\Omega$,
\begin{align}
  \EE\Bigl[\bigl(\Delta_3 Y_{\{1,2\}}(\Omega)\bigr)^2 \;\Big|\; \mathcal{F}_{\{1,2\}}\Bigr]
  &\le \ell_{12}(\Omega)\int_{\{y \ne 0:\,x(y)\in\Omega\}}
      \EE\Bigl[\bigl(\Delta_3\varepsilon_{12}(x(y))\bigr)^2 \;\Big|\; \mathcal{F}_{\{1,2\}}\Bigr]\dd y
  \notag \\
  &\le \beta_\Omega^2 \sup_{y \ne 0:\,x(y)\in\Omega}
      \EE\Bigl[\bigl(\Delta_3\varepsilon_{12}(x(y))\bigr)^2 \;\Big|\; \mathcal{F}_{\{1,2\}}\Bigr],
  \label{eq:es-slice-cs}
\end{align}
and analogously for $\Delta_3\Delta_4 Y_{\{1,2\}}(\Omega)$ (at $\dimin = 2$ the
integral is the single value at $y = 0$). To bound the pointwise conditional
second moments at fixed $x = x(y)$, we partition the probability space into
the three disjoint events $\mathcal{A}_B(x)$, $\mathcal{E}_{\mathrm{DS}}(x)
\setminus \mathcal{A}_B(x)$ and $\mathcal{E}_{\mathrm{DS}}(x)^c$.

On $\mathcal{A}_B(x)$, squaring~\eqref{eq:es-good-pointwise} and integrating
out the independent $\chi_{\dimin}^2$ variables $R_3^2, R_3'^2, R_4^2, R_4'^2$
conditionally on $\mathcal{F}_{\{1,2\}}$ yields $C_{\dimin}\ngen^{-6} R_1^4 R_2^4$
for $(\Delta_3\varepsilon_{12}(x))^2$ and $C_{\dimin}\ngen^{-8} R_1^4 R_2^4$ for
$(\Delta_3\Delta_4\varepsilon_{12}(x))^2$. On the moderate-norm tail
$\mathcal{E}_{\mathrm{DS}}(x) \setminus \mathcal{A}_B(x)$, every perturbed
matrix $G \succeq G_0(x)$ satisfies $\lambda_{\min}(G) \ge c_0\sigw^2\ngen$. Thus
Proposition~\ref{prop:deficit-integral}\textup{(iii)} gives
$\abs{\Delta_3\varepsilon_{12}(x)} \le C\ngen^{-2} R_1^2 R_2^2$ and
$\abs{\Delta_3\Delta_4\varepsilon_{12}(x)} \le C\ngen^{-2} R_1^2 R_2^2$. At the
same time,~\eqref{eq:es-good-event} implies that $R_1^2 + R_2^2 + R_3^2 + R_3'^2
+ R_4^2 + R_4'^2 > c\ngen$ on $\mathcal{E}_{\mathrm{DS}}(x) \setminus
\mathcal{A}_B(x)$. Thus $\mathbf{1}_{\mathcal{E}_{\mathrm{DS}}(x)\setminus
\mathcal{A}_B(x)} \le (c\ngen)^{-k}(R_1^2 + R_2^2 + R_3^2 + R_3'^2 + R_4^2 +
R_4'^2)^k$ pathwise for $k \in \{1, 2\}$. Multiplying these bounds before
squaring (with $k = 1$ for $\Delta_3$ and $k = 2$ for $\Delta_3\Delta_4$) gives
\begin{align*}
  \abs{\Delta_3\varepsilon_{12}(x)}\,\mathbf{1}_{\mathcal{E}_{\mathrm{DS}}(x)\setminus\mathcal{A}_B(x)}
  &\le C\,\ngen^{-3}\,R_1^2 R_2^2\bigl(R_1^2 + R_2^2 + R_3^2 + R_3'^2 + R_4^2 + R_4'^2\bigr), \\
  \abs{\Delta_3\Delta_4\varepsilon_{12}(x)}\,\mathbf{1}_{\mathcal{E}_{\mathrm{DS}}(x)\setminus\mathcal{A}_B(x)}
  &\le C\,\ngen^{-4}\,R_1^2 R_2^2\bigl(R_1^2 + R_2^2 + R_3^2 + R_3'^2 + R_4^2 + R_4'^2\bigr)^2,
\end{align*}
whose conditional expectations after squaring are at most $C_{\dimin}\ngen^{-6}
R_1^4 R_2^4(1 + R_1^4 + R_2^4)$ and $C_{\dimin}\ngen^{-8} R_1^4 R_2^4(1 + R_1^8
+ R_2^8)$. On $\mathcal{E}_{\mathrm{DS}}(x)^c$, the deterministic bound
$\abs{\varepsilon_{12}} \le \pi/2$ gives $\abs{\Delta_3\varepsilon_{12}(x)} \le
\pi$ and $\abs{\Delta_3\Delta_4\varepsilon_{12}(x)} \le 2\pi$, while
$\PP(\mathcal{E}_{\mathrm{DS}}(x)^c \mid \mathcal{F}_{\{1,2\}}) \le
Ce^{-c_1\ngen}$. Summing the three contributions uniformly in $y \ne 0$
and taking expectations over the independent $\chi_{\dimin}^2$ variables
$(R_1^2, R_2^2)$ gives~\eqref{eq:es-replacement-bounds}, and substituting into
Lemma~\ref{lem:master-es}~\eqref{eq:es-master-Delta}
yields~\eqref{eq:es-variance-bound}.
\end{proof}

\subsection{Thinned-Wishart kernel rate and transfer theorem}
\label{subsec:kernel-rate}

\begin{lemma}
\label{lem:kernel-rate}
Under the standing hypotheses of this section, let $e_{\ngen}(w_c, w_{c'})$ be
the thinned-Wishart kernel~\eqref{eq:bone-kernel-def}, and set $R_c :=
\norm{w_c}/\sigw$, $R_{c'} := \norm{w_{c'}}/\sigw$.
\begin{enumerate}[label=(\roman*)]
  \item There exists
$C_{\dimin,\sigw} > 0$ depending only on $(\dimin, \sigw)$ such that for all
linearly independent $(w_c, w_{c'})$ and $\ngen \ge \dimin + 2$,
\begin{gather}
  \abs[\Big]{\frac{\ngen^2}{4}\, e_{\ngen}(w_c, w_{c'}) - \Zlim(w_c, w_{c'})}
  \le C_{\dimin,\sigw}\bigl(1 + R_c^2 R_{c'}^2(1 + R_c^2 + R_{c'}^2)\bigr)\, \ngen^{-1},
  \label{eq:kernel-pointwise-rate} \\
  \norm[\Big]{\Bigl(\frac{\ngen^2}{4}\, e_{\ngen}(w_c, w_{c'}) - \Zlim(w_c, w_{c'})\Bigr)\ell_{cc'}(\Omega)}_{L^2(\PP)}
  \le C_{\dimin,\sigw}\, \beta_\Omega\, \ngen^{-1}.
  \label{eq:kernel-single-rate}
\end{gather}
  \item Across all $\binom{\ngen}{2}$ bones $B = \{c, c'\}$, the centred
kernel discrepancies are pairwise orthogonal in $L^2(\PP)$, so
\begin{equation}\label{eq:kernel-sum-rate}
  \norm[\Big]{\frac{\ngen}{4} P_2\bigl(\Regge(\Omega)\bigr) - U_{\ngen,\infty}(\Omega)}_{L^2(\PP)}
  \le C_{\dimin,\sigw}\, \beta_\Omega\, \ngen^{-1}.
\end{equation}
\end{enumerate}
\end{lemma}

\begin{proof}
By~\eqref{eq:bone-kernel-def}, $e_{\ngen}(w_c, w_{c'})$ is the $\Bin(\ngen - 2,
1/2)$ average of $\bar{\varepsilon}_m(L_B)$, so~\eqref{eq:kernel-cond-nu} and
the binomial-averaging bound at the end of the proof of
Theorem~\ref{thm:fisher-limit-mark}\textup{(ii)}
give~\eqref{eq:kernel-pointwise-rate}. Taking $L^2(\PP)$ norms with
$\ell_{cc'}(\Omega) \le \beta_\Omega$ gives~\eqref{eq:kernel-single-rate}. Since
$\delta_{\ngen, cc'}(\Omega) := \bigl(\frac{\ngen^2}{4} e_{\ngen}(w_c, w_{c'}) -
\Zlim(w_c, w_{c'})\bigr)\ell_{cc'}(\Omega) \in \mathcal{V}_{\{c, c'\}}$ are
pairwise orthogonal across $\{c, c'\}$, summing $\frac{1}{\ngen}\sum_{c <
c'}\delta_{\ngen, cc'}(\Omega)$ gives~\eqref{eq:kernel-sum-rate}.
\end{proof}

\begin{theorem}[Second-chaos transfer]
\label{thm:second-chaos-transfer}
Under the standing hypotheses of this section, let $U_{\ngen,\infty}(\Omega) :=
\frac{1}{\ngen} \sum_{1 \le c < c' \le \ngen} \Zlim(w_c, w_{c'})\,
\ell_{cc'}(\Omega)$.
\begin{enumerate}[label=(\roman*)]
  \item With $\mathcal{D}_2(\Omega) := \frac{\ngen}{4} P_2(\Regge(\Omega)) -
    U_{\ngen,\infty}(\Omega) \in \bigoplus_{\abs{A}=2} \mathcal{V}_A$ and
    $\mathcal{D}_{\ge 3}(\Omega) := \frac{\ngen}{4} R_{\ge 3}(\Omega) \in
    \bigoplus_{\abs{A}\ge 3} \mathcal{V}_A$, we have $\EE[\mathcal{D}_2(\Omega)\,
    \mathcal{D}_{\ge 3}(\Omega)] = 0$ and
    \begin{equation}\label{eq:transfer-pythagorean}
      \begin{gathered}
        \frac{\ngen}{4}\Regge(\Omega) - U_{\ngen,\infty}(\Omega)
        = \mathcal{D}_2(\Omega) + \mathcal{D}_{\ge 3}(\Omega), \\
        \norm[\Big]{\frac{\ngen}{4}\Regge(\Omega) - U_{\ngen,\infty}(\Omega)}_{L^2(\PP)}^2
        = \norm{\mathcal{D}_2(\Omega)}_{L^2(\PP)}^2
          + \norm{\mathcal{D}_{\ge 3}(\Omega)}_{L^2(\PP)}^2.
      \end{gathered}
    \end{equation}
  \item There exists $C_{\dimin,\sigw} > 0$ depending only on $(\dimin, \sigw)$
    such that for all $\ngen \ge \dimin + 2$,
    \begin{align}
      \norm[\Big]{\frac{\ngen}{4}\Regge(\Omega) - U_{\ngen,\infty}(\Omega)}_{L^2(\PP)}
      &\le C_{\dimin,\sigw}\,\beta_\Omega\, \ngen^{-1/2},
      \label{eq:transfer-l2-bound} \\
      \abs[\Big]{
        \norm[\Big]{\frac{\ngen}{4}\Regge(\Omega)}_{L^2(\PP)}
        - \norm{U_{\ngen,\infty}(\Omega)}_{L^2(\PP)}
      }
      &\le C_{\dimin,\sigw}\, \beta_\Omega\, \ngen^{-1/2}.
      \label{eq:transfer-std-bound}
    \end{align}
    The $O(\ngen^{-1/2})$ rate in~\eqref{eq:transfer-l2-bound} is sharp at fixed
    $(\dimin, \sigw)$: whenever $\PP(\ell_{12}(\Omega) > 0) > 0$, one has
    $\bigl\lVert\frac{\ngen}{4}P_3(\Regge(\Omega))\bigr\rVert_{L^2(\PP)} =
    \Theta(\ngen^{-1/2})$ as $\ngen \to \infty$, with implied constants
    depending on $(\dimin, \sigw, \Omega)$.
  \item On $\mathcal{Z} := \RR^{\dimin} \times \RR$ with $\xi_h = (w_h, b_h)
    \stackrel{\mathrm{i.i.d.}}{\sim} P_\xi$, the kernel $K_\Omega(\xi_c,
    \xi_{c'}) := \Zlim(w_c, w_{c'})\ell_{cc'}(\Omega)$ is symmetric,
    square-integrable and completely degenerate ($\EE[K_\Omega(\xi_1, \xi_2) \mid
    \xi_1] = 0$ a.s.), and as $\ngen \to \infty$,
    \begin{equation}\label{eq:general-window-limit}
      \frac{\ngen}{4}\Regge(\Omega)
      \;\xrightarrow[\ngen \to \infty]{d}\;
      U_\infty(\Omega)
      := \frac{1}{2}\sum_{k=1}^\infty \mu_{k,\Omega}(\gamma_k^2 - 1),
      \qquad
      \gamma_k \stackrel{\mathrm{i.i.d.}}{\sim} \Norm(0, 1),
    \end{equation}
    where $(\mu_{k,\Omega})_{k \ge 1}$ are the non-zero eigenvalues (counted with
    multiplicity) of $(T_{K_\Omega} f)(\xi) := \EE_{\xi'}[K_\Omega(\xi, \xi')
    f(\xi')]$ on $L^2(\RR^{\dimin} \times \RR, P_\xi)$.
\end{enumerate}
\end{theorem}

\begin{proof}
Proposition~\ref{thm:first-chaos-annihilation} and Hoeffding orthogonality
give~\eqref{eq:transfer-pythagorean}, and Lemmas~\ref{lem:cross-bone-es}
and~\ref{lem:kernel-rate} give~\eqref{eq:transfer-l2-bound}--\eqref{eq:transfer-std-bound}.
The sharpness bound $\norm{\frac{\ngen}{4}P_3(\Regge(\Omega))}_{L^2(\PP)} =
\Theta(\ngen^{-1/2})$ is proved in Appendix~\ref{sec:app-wishart}. The idea is
to project $\pi_{\{1,2,3\}}\Regge(\Omega)$ onto its $\mathcal{I}_3$-even
component, which annihilates $\pi_{\{1,2,3\}}Y_{\{1,3\}}(\Omega)$ and
$\pi_{\{1,2,3\}}Y_{\{2,3\}}(\Omega)$ and leaves a leading linear fluctuation of
$\pi_{\{1,2,3\}}Y_{\{1,2\}}(\Omega)$ of exact order $\ngen^{-3}$ in $L^2(\PP)$.
Finally, Lemma~\ref{lem:sign-flip-oddness} gives complete degeneracy
of $K_\Omega$, so combining Proposition~\ref{prop:degenerate-ustat-limit}
with~\eqref{eq:transfer-l2-bound} and Slutsky's theorem
yields~\eqref{eq:general-window-limit}.
\end{proof}

\begin{remark}
\label{rem:transfer-finite-H}
The constants in Theorem~\ref{thm:second-chaos-transfer} are uniform over
$\ngen \ge \dimin + 2$, not merely for $\ngen$ large. The reason is that
$\abs{\varepsilon_{cc'}} \le \pi/2$ holds deterministically. Thus, the events on
which the active count $m_B(x)$ is too small for the Wishart estimates to apply
(for instance $\{m_B(x) < \dimin + 3\}$, which has probability at most
$Ce^{-c\ngen}$) are simply absorbed into the constant, as at
the start of the proofs of Theorem~\ref{thm:fisher-limit-mark} and
Lemma~\ref{lem:cross-bone-es}.
\end{remark}

\begin{remark}[Intrinsically weighted curvature]
\label{rem:cms-intrinsic-equivalence}
For $\dimin = 2$, $\Haus{0}$ is counting measure and $\Regge(\Omega)$ coincides
with the polyhedral curvature of Regge~\cite{Regge1961} and Cheeger, M\"uller
and Schrader~\cite{CheegerMullerSchrader1984}. For $\dimin \ge 3$, those authors weight $\varepsilon_{cc'}(x)$ by the intrinsic bone volume element $\sqrt{\det
G_B(x)|_{TB_{cc'}}}\,\dd\Haus{\dimin-2}(x) = \sqrt{\det(I_{\dimin-2} +
\mathsf{W}_{TT}(x))}\,\dd\Haus{\dimin-2}(x)$, defining a measure
$\mu_{\mathrm{CMS}}(\Omega)$. In contrast, $\Regge(\Omega)$ weights by Euclidean
$(\dimin-2)$-measure so as to isolate the transverse dihedral deficit from
tangential metric stretching along $TB_{cc'}$. The matrix $\mathsf{W}_{TT}(x) =
Q_B\transp\bigl(\sum_{h \in A(x)} w_h w_h\transp\bigr)Q_B$ depends only on
generators $h \notin \{c, c'\}$ and the subspace $TB_{cc'}$. Since both are fixed by
$\mathcal{I}_c$ and $\mathcal{I}_{c'}$, the intrinsic weight is pathwise
invariant under the orientation flips of Lemma~\ref{lem:sign-flip-oddness}.
Consequently, pathwise antipodal oddness and exact first-chaos annihilation
(Lemma~\ref{lem:hoeffding-support},
Proposition~\ref{thm:first-chaos-annihilation} and
Proposition~\ref{thm:parity-orthogonality}) hold verbatim for
$\mu_{\mathrm{CMS}}(\Omega)$. Under~\ref{assump:A2}, conditional independence of
$\mathsf{W}_{TT}(x)$ and $S_N^{\mathrm{unreg}}(x)$ in
Theorem~\ref{thm:conditional-wishart}(iii) further suggests that
$(\ngen\sigw^2/2)^{-(\dimin-2)/2}\frac{\ngen}{4}\mu_{\mathrm{CMS}}(\Omega)$
shares the second-chaos limit $U_\infty(\Omega)$. A quantitative $L^2$
transfer analysis for $\mu_{\mathrm{CMS}}(\Omega)$ is deferred to future work.
\end{remark}
\section{Planar cone angle, exact \texorpdfstring{$-\sec\sqrt{\pi^2 + \pi s}$}{-sec sqrt(pi2 + pi s)} law and fourth-moment universality}
\label{sec:planar}

When $\dimin = 2$ and $\Omega = \RR^2$, every bone $B_{cc'} = \Pi_c \cap
\Pi_{c'} = \{x_{cc'}\}$ is a point with $\ell_{cc'}(\RR^2) = 1$ and
$\beta_{\RR^2} = 1$, reducing the limiting U-statistic to the pure-mark sum
$U_{\ngen,\infty}(\RR^2) = \frac{1}{\ngen}\sum_{1 \le c < c' \le \ngen}
\Zlim(w_c, w_{c'})$. We identify $\Regge(\RR^2)$ via polyhedral Gauss--Bonnet
with the far-field cone-angle deficit $2\pi - \Theta_\infty(W)$, compute the
exact spectral law $\mathcal{L}_2$, prove fourth-moment universality via an
offset-free circle integral, and record Kolmogorov rates and finite-$\ngen$
moments.

\subsection{Planar polyhedral Gauss--Bonnet identity and offset independence}
\label{subsec:planar-gauss-bonnet}

In $\dimin = 2$, an arrangement $\arr = \{\Pi_1, \dots, \Pi_{\ngen}\}$ of
$\ngen$ lines in general position $\mathcal{E}_{\mathrm{nd}}(W, b)$ partitions
$\RR^2$ into $\binom{\ngen}{2} + \ngen + 1$ cells~\cite{Miles1964}. These comprise $\binom{\ngen - 1}{2}$
bounded polygons and $2\ngen$ unbounded wedges $C_1^\infty, \dots,
C_{2\ngen}^\infty$, ordered cyclically at infinity. Each $C_k^\infty$ is bounded
by rays along $u_k, u_{k+1} \in \Sph{1}$ ($u_{2\ngen + 1} := u_1$) parallel to
$\Pi_{c(k)}$ and $\Pi_{c(k+1)}$. It carries the constant metric $G_{(k)}(W) :=
I_2 + \sum_{h \in A_k^\infty} w_h w_h\transp$. Here $A_k^\infty := \{h \in
[\ngen] : \inner{w_h}{u} > 0 \text{ on } \mathrm{int}(C_k^\infty \cap
\Sph{1})\}$ depends \emph{only} on $W$. The intrinsic opening angle of
$C_k^\infty$ in $G_{(k)}(W)$ and the \emph{far-field cone angle} are
\begin{equation}\label{eq:planar-wedge-angle}
  \alpha_k^\infty(W)
  := \arccos\!\left(
    \frac{u_k\transp G_{(k)}(W) u_{k+1}}{\sqrt{u_k\transp G_{(k)}(W) u_k}\sqrt{u_{k+1}\transp G_{(k)}(W) u_{k+1}}}
  \right) \in (0, \pi),
  \qquad
  \Theta_\infty(W) := \sum_{k=1}^{2\ngen} \alpha_k^\infty(W).
\end{equation}

\begin{proposition}
\label{prop:planar-gauss-bonnet}
Let $\dimin = 2$ and $\ngen \ge 2$.
\begin{enumerate}[label=(\roman*)]
  \item On $\mathcal{E}_{\mathrm{nd}}(W, b)$, the total Regge curvature over
    $\RR^2$ equals the cone-angle deficit:
    \begin{equation}\label{eq:planar-gauss-bonnet}
      \Regge(\RR^2)
      = \sum_{1 \le c < c' \le \ngen} \varepsilon_{cc'}(x_{cc'})
      = 2\pi - \Theta_\infty(W),
    \end{equation}
    where $x_{cc'}$ is the unique vertex of $\Pi_c \cap \Pi_{c'}$.
  \item Across all $(W, b) \in \mathcal{E}_{\mathrm{nd}}(W, b)$, $\Regge(\RR^2)
    = 2\pi - \Theta_\infty(W)$ is independent of $b$; when $b \equiv 0$ and
    $\norm{w_c \wedge w_{c'}} > 0$ for all $c < c'$, it coincides
    with~\eqref{eq:model-apex-convention}, $\Regge(\Omega) = (2\pi -
    \Theta_\infty(W))\mathbf{1}\{0 \in \Omega\}$.
  \item If the rows of $W$ are a.s.\ pairwise non-collinear and $\tilde{b}_1,
    \dots, \tilde{b}_{\ngen} \stackrel{\mathrm{i.i.d.}}{\sim} \Norm(0, 1)$ are
    adjoined independently of $W$, then $\PP(\mathcal{E}_{\mathrm{nd}}(W,
    \tilde{b})) = 1$ and $2\pi - \Theta_\infty(W) = \Regge(\RR^2; W, \tilde{b})$
    a.s. If $w_1, \dots, w_{\ngen}$ are independent and centrally symmetric,
    then $\tilde{\xi}_h := (w_h, \tilde{b}_h)$ ($h \in [\ngen]$) are independent
    and jointly centrally symmetric.
\end{enumerate}
\end{proposition}

\begin{proof}
Choose $R > \max_{c < c'}\norm{x_{cc'}}_2$. Outside $B(0, R)$ the arrangement
consists of $2\ngen$ pairwise disjoint rays $r_k \subset \Pi_{c(k)}$ meeting
$\partial B(0, R)$ in cyclic order at $v_k(R) := r_k \cap \partial B(0, R)$. By
convexity of $C_k^\infty$, both the arc of $\partial B(0, R)$ and the chord
$e_k(R) := [v_k(R), v_{k+1}(R)]$ have relative interiors in $C_k^\infty$. The
circular segment between them meets no line $\Pi_h$ and contains no vertex
$x_{cc'}$. Hence the polygonal disc $P_R \subset \overline{B(0, R)}$ bounded by
the $2\ngen$ chords $e_k(R)$ (cf.~\cite{Banchoff1967}) contains every $x_{cc'}$
in $\mathrm{int}(P_R)$. It carries a polygonal complex with $N_{\mathrm{int}} =
\ngen^2$ interior edges, $N_\partial = 2\ngen$ boundary edges and
$N_{\mathrm{face}} = \binom{\ngen}{2} + \ngen + 1$ flat polygonal faces,
satisfying $\sum_{f=1}^{N_{\mathrm{face}}} k_f = 2N_{\mathrm{int}} +
N_\partial$.

Each interior vertex $x_{cc'}$ contributes four quadrant corners summing to
$2\pi - \varepsilon_{cc'}(x_{cc'})$. In the constant metric $G_{(k)}(W)$,
$v_k(R)$, $v_{k+1}(R)$ and $x_{c(k), c(k+1)}$ form a flat triangle with apex
angle $\alpha_k^\infty(W)$, so the two outer corners of $C_k^\infty \cap P_R$ at
$v_k(R)$ and $v_{k+1}(R)$ sum to $\pi - \alpha_k^\infty(W)$. Summing the flat
interior angle sums $(k_f - 2)\pi$ over $f = 1, \dots, N_{\mathrm{face}}$ gives
\[
  \sum_{1 \le c < c' \le \ngen}\bigl(2\pi - \varepsilon_{cc'}(x_{cc'})\bigr)
  + \sum_{k=1}^{2\ngen}\bigl(\pi - \alpha_k^\infty(W)\bigr)
  = (2N_{\mathrm{int}} + N_\partial - 2N_{\mathrm{face}})\pi
  = \Bigl(2\binom{\ngen}{2} + 2\ngen - 2\Bigr)\pi,
\]
which rearranges to~\eqref{eq:planar-gauss-bonnet}, and part~(ii) follows because
$\alpha_k^\infty(W)$ in~\eqref{eq:planar-wedge-angle} depends only on $W$.
For~(iii), triple concurrence $\tilde{b}_{h_3} = -\ip{w_{h_3}}{x_{h_1 h_2}}$ has
probability zero since $\tilde{b}_{h_3} \sim \Norm(0, 1)$ is atomless and
independent of $(W, \tilde{b}_{h_1}, \tilde{b}_{h_2})$, giving
$\PP(\mathcal{E}_{\mathrm{nd}}(W, \tilde{b})) = 1$. The independence and
symmetry of $\tilde{b}_h$ pass to $\tilde{\xi}_h$.
\end{proof}

\begin{lemma}
\label{lem:planar-offset-reduction}
Let $\dimin = 2$ and $\ngen \ge 2$, and suppose the rows of $W$ are a.s.\
pairwise non-collinear. Say that $b \in \RR^{\ngen}$ satisfies the \emph{planar
general-position condition} if either $\PP(\mathcal{E}_{\mathrm{nd}}(W, b)) = 1$
or $b \equiv 0$ under~\eqref{eq:model-apex-convention}. For any such $b$,
adjoining $\tilde{b}_1, \dots, \tilde{b}_{\ngen} \stackrel{\mathrm{i.i.d.}}{\sim}
\Norm(0, 1)$ independently of $(W, b)$ yields, almost surely,
\begin{align*}
  \bigl(\Regge(\RR^2; W, b),\, U_{\ngen,\infty}(\RR^2)\bigr)
  &= \bigl(2\pi - \Theta_\infty(W),\, U_{\ngen,\infty}(W)\bigr) \\
  &= \bigl(\Regge(\RR^2; W, \tilde{b}),\, U_{\ngen,\infty}(W)\bigr),
\end{align*}
with $\PP(\mathcal{E}_{\mathrm{nd}}(W, \tilde{b})) = 1$. Suppose $w_1, \dots,
w_{\ngen}$ are independent and the marginal law of each $w_h$ charges no affine
line in $\RR^2$ (in particular under~\ref{assump:A2}). Then the planar
general-position condition holds whenever $b \equiv 0$. It also holds whenever $b$ is
independent of $W$ (including deterministic $b \in \RR^{\ngen}$) and satisfies
$\abs{\{h : b_h = 0\}} \le 2$ a.s.
\end{lemma}

\begin{proof}
Since $\ell_{cc'}(\RR^2) = 1$, $U_{\ngen,\infty}(\RR^2) = U_{\ngen,\infty}(W)$,
and Proposition~\ref{prop:planar-gauss-bonnet}(ii)--(iii) identifies both
$\Regge(\RR^2; W, b)$ and $\Regge(\RR^2; W, \tilde{b})$ a.s.\ with $2\pi -
\Theta_\infty(W)$. When $w_1, \dots, w_{\ngen}$ are independent and each
marginal charges no affine line, $w_{h_3}$ is independent of $(w_{h_1}, b_{h_1},
w_{h_2}, b_{h_2}, b_{h_3})$. Thus $\ip{w_{h_3}}{x_{h_1 h_2}} + b_{h_3} = 0$ has
conditional probability zero whenever $(x_{h_1 h_2}, b_{h_3}) \ne (0, 0)$.
\end{proof}

\subsection{Spectral decomposition and the planar limit law \texorpdfstring{$\mathcal{L}_2$}{L\_2}}
\label{subsec:planar-spectral}

\begin{proposition}
\label{prop:planar-spectrum}
Let $\dimin = 2$ and let $P_w$ satisfy~\ref{assump:A3} on $\RR^2$ with $\Sigma_w
= \sigw^2 I_2$ ($\sigw > 0$). Write $w = \norm{w}(\cos\phi, \sin\phi)\transp$
with $\phi \sim \Unif([0, 2\pi))$ independent of $\varrho(w) :=
\frac{\norm{w}^2}{2\sigw^2}$, and recall $\tau_2(P_w) = \EE[\varrho(w)^2] <
\infty$ ($\tau_2 = 2$ under~\ref{assump:A2}).
\begin{enumerate}[label=(\roman*)]
  \item For all $w, w' \in \RR^2 \setminus \{0\}$,
    \[
      \Zlim(w, w')
      = -\frac{\norm{w}^2 \norm{w'}^2}{4\sigw^4} \cos(\phi - \phi')\,\abs{\sin(\phi - \phi')}
      = \varrho(w)\,\varrho(w')\, K_{\mathrm{ang}}(\phi - \phi'),
    \]
    where $K_{\mathrm{ang}}(\alpha) := -\cos\alpha\,\abs{\sin\alpha}$ is
    continuous, even, $2\pi$-periodic and $\pi$-antiperiodic
    ($K_{\mathrm{ang}}(\alpha + \pi) = -K_{\mathrm{ang}}(\alpha)$), with
    uniformly and absolutely convergent Fourier series
    \begin{equation}\label{eq:planar-eigenvalues}
      K_{\mathrm{ang}}(\alpha)
      = \sum_{m=0}^\infty \lambda_m \cos\bigl((2m + 1)\alpha\bigr),
      \quad
      \lambda_m
      := \frac{4}{\pi(2m - 1)(2m + 3)}
       = \frac{1}{\pi}\Bigl(\frac{1}{2m - 1} - \frac{1}{2m + 3}\Bigr).
    \end{equation}
    Here $\lambda_0 = -\frac{4}{3\pi} < 0$, $\lambda_m > 0$ for $m \ge 1$
    ($\lambda_1 = \frac{4}{5\pi}$, $\lambda_2 = \frac{4}{21\pi}$, $\dots$),
    $\sum_{m=0}^\infty \lambda_m = K_{\mathrm{ang}}(0) = 0$ and
    $\sum_{m=0}^\infty \lambda_m^2 = 2\int_0^{2\pi}
    K_{\mathrm{ang}}(\alpha)^2\frac{\dd\alpha}{2\pi} = \frac{1}{4}$.
  \item For $m \ge 0$, the functions
    \begin{equation}\label{eq:planar-mercer}
      \psi_{m,1}(w) := \sqrt{\frac{2}{\tau_2(P_w)}}\, \varrho(w)\cos((2m + 1)\phi),
      \quad
      \psi_{m,2}(w) := \sqrt{\frac{2}{\tau_2(P_w)}}\, \varrho(w)\sin((2m + 1)\phi)
    \end{equation}
    are centred and orthonormal in $L^2(\RR^2, P_w)$, with
    $\Zlim(w, w') = \allowbreak \sum_{m=0}^\infty \allowbreak
    \bigl(\frac{\tau_2(P_w)}{2}\lambda_m\bigr)\allowbreak
    [\psi_{m,1}(w)\psi_{m,1}(w') + \allowbreak \psi_{m,2}(w)\psi_{m,2}(w')]$ in
    $L^2(P_w \otimes P_w)$. Under~\ref{assump:A2} ($\tau_2 = 2$), the non-zero
    eigenvalues of $(T_{\Zlim} f)(w) = \EE_{w'}[\Zlim(w, w') f(w')]$ are
    $\{\lambda_m\}_{m \ge 0}$, each of multiplicity $2$, with $\EE[\Zlim(w,
    w')^2] = \frac{1}{2}$.
\end{enumerate}
\end{proposition}

\begin{proof}
Evenness and $\pi$-antiperiodicity eliminate all sine and even cosine harmonics.
For $m \ge 0$, $2\cos\alpha\sin\alpha\cos((2m+1)\alpha) =
\frac{1}{2}[\sin((2m+3)\alpha) - \sin((2m-1)\alpha)]$ on $[0, \pi]$ gives
\[
  \lambda_m
  = -\frac{1}{2\pi}\int_0^\pi \bigl[\sin((2m + 3)\alpha) - \sin((2m - 1)\alpha)\bigr]\,\dd\alpha
  = \frac{1}{\pi}\Bigl(\frac{1}{2m - 1} - \frac{1}{2m + 3}\Bigr),
\]
proving~\eqref{eq:planar-eigenvalues}. Independence of $\varrho(w)$ and $\phi
\sim \Unif([0, 2\pi))$ gives $\EE[\psi_{m,i}(w)\psi_{m',j}(w)] =
\delta_{m,m'}\delta_{i,j}$ and~\eqref{eq:planar-mercer}.
\end{proof}

\begin{figure}[t]
\centering
\includegraphics[width=\textwidth]{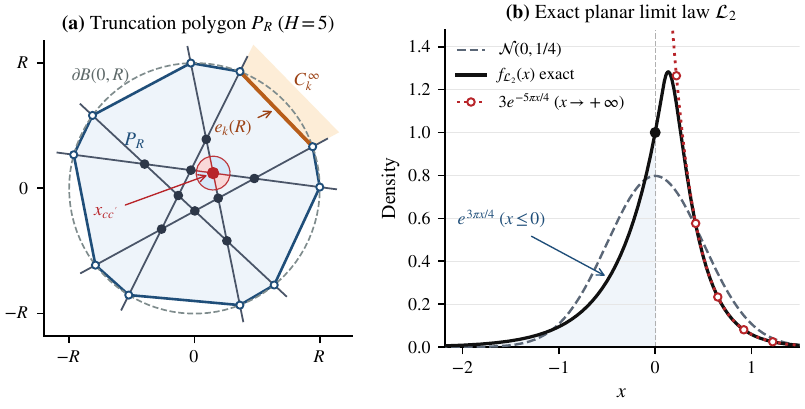}
\caption{Planar Gauss--Bonnet truncation and exact limit density $f_{\mathcal{L}_2}$ ($\dimin = 2$). (a)~Truncation polygon $P_R \subset \overline{B(0, R)}$ of Proposition~\ref{prop:planar-gauss-bonnet} for $\ngen = 5$ lines in general position, bounded by the $2\ngen$ chords $e_k(R) = [v_k(R), v_{k+1}(R)]$ across the unbounded wedges $C_k^\infty$ (with $C_k^\infty \setminus P_R$ shaded orange for one wedge) and enclosing every vertex $x_{cc'}$ together with its four quadrant sectors $\sigma \in \{\pm 1\}^2$. (b)~Exact density $f_{\mathcal{L}_2}(x)$ obtained by Fourier inversion of $\varphi_{\mathcal{L}_2}(t) = -\sec\sqrt{\pi^2 + i\pi t}$~\eqref{eq:planar-mgf}, plotted against the variance-matched Gaussian $\Norm(0, 1/4)$. Note the negative skew, the exact exponential left branch $f_{\mathcal{L}_2}(x) = e^{3\pi x/4}$ for $x \le 0$ ($f_{\mathcal{L}_2}(0) = 1$), and the lighter right tail $f_{\mathcal{L}_2}(x) = 3e^{-5\pi x/4}(1 + O(e^{-4\pi x}))$ as $x \to +\infty$~\eqref{eq:planar-tails}.}
\label{fig:planar-limit}
\end{figure}

Each multiplicity-$2$ eigenvalue $\lambda_m$ pairs two independent $\chi_1^2$
modes into $\frac{1}{2}(\gamma_{m,1}^2 + \gamma_{m,2}^2) = E_m \sim
\mathrm{Exp}(1)$, and Euler's cosine product evaluates the moment generating
function in closed form (see Figure~\ref{fig:planar-limit}).

\begin{theorem}[Planar limit law $\mathcal{L}_2$]
\label{thm:planar-limit-law}
Let $\dimin = 2$ and $\Omega = \RR^2$. Suppose $w_1, \dots, w_{\ngen}
\stackrel{\mathrm{i.i.d.}}{\sim} \Norm(0, \sigw^2 I_2)$ ($\sigw > 0$) and that
$b \in \RR^{\ngen}$ satisfies the planar general-position condition of
Lemma~\ref{lem:planar-offset-reduction} (no symmetry of $b$ is required).

Let
$\lambda_m = \frac{4}{\pi(2m - 1)(2m + 3)}$ ($m \ge 0$) and $E_m
\stackrel{\mathrm{i.i.d.}}{\sim} \mathrm{Exp}(1)$.
\begin{enumerate}[label=(\roman*)]
  \item As $\ngen \to \infty$, both $U_{\ngen,\infty}(\RR^2)$ and
    $\frac{\ngen}{4}\Regge(\RR^2) = \frac{\ngen}{4}(2\pi - \Theta_\infty(W))$
    converge in law, with $\bigl\lVert\frac{\ngen}{4}(2\pi - \Theta_\infty(W)) -
    U_{\ngen,\infty}(\RR^2)\bigr\rVert_{L^2(\PP)} \to 0$, to
    \begin{equation}\label{eq:planar-limit-series}
      \mathcal{L}_2
      := \sum_{m=0}^\infty \lambda_m (E_m - 1)
       = \sum_{m=0}^\infty \lambda_m E_m,
    \end{equation}
    where the series converges a.s.\ and in $L^2(\PP)$ ($\sum_{m=0}^\infty
    \abs{\lambda_m} = \frac{8}{3\pi}$, $\sum_{m=0}^\infty \lambda_m^2 =
    \frac{1}{4}$, $\sum_{m=0}^\infty \lambda_m = 0$).
  \item The canonical product $P(s) := \prod_{m=0}^\infty (1 - \lambda_m s)$ is
    entire of order $1/2$ and genus $0$, and $M_{\mathcal{L}_2}(s) :=
    \EE[e^{s \mathcal{L}_2}]$ is finite if and only if $s \in
    \bigl(-\frac{3\pi}{4}, \frac{5\pi}{4}\bigr)$, where
    \begin{equation}\label{eq:planar-mgf}
      \begin{gathered}
        \prod_{m=0}^\infty (1 - \lambda_m s)
        = -\cos\bigl(\sqrt{\pi^2 + \pi s}\bigr)
        \quad (\forall s \in \CC), \\
        M_{\mathcal{L}_2}(s)
        = \prod_{m=0}^\infty \frac{1}{1 - \lambda_m s}
        = -\sec\bigl(\sqrt{\pi^2 + \pi s}\bigr),
      \end{gathered}
    \end{equation}
    and $\varphi_{\mathcal{L}_2}(t) = \bigl[-\cos\bigl(\sqrt{\pi^2 + i\pi
    t}\bigr)\bigr]^{-1}$ for all $t \in \RR$.
  \item The cumulants $\kappa_r(\mathcal{L}_2) = (r - 1)!\sum_{m=0}^\infty
    \lambda_m^r$ ($r \ge 2$), skewness and excess kurtosis satisfy
    \begin{align}
      \kappa_1(\mathcal{L}_2) &= 0, \quad
      \kappa_2(\mathcal{L}_2) = \Var(\mathcal{L}_2) = \frac{1}{4}, \quad
      \kappa_3(\mathcal{L}_2) = -\frac{3}{8\pi}, \quad
      \kappa_4(\mathcal{L}_2) = \frac{1}{8} + \frac{15}{16\pi^2},
      \label{eq:planar-cumulants} \\
      \gamma_1(\mathcal{L}_2)
      &:= \frac{\kappa_3(\mathcal{L}_2)}{\kappa_2(\mathcal{L}_2)^{3/2}}
       = -\frac{3}{\pi} \approx -0.95493, \qquad
      \gamma_2(\mathcal{L}_2)
      := \frac{\kappa_4(\mathcal{L}_2)}{\kappa_2(\mathcal{L}_2)^2}
       = 2 + \frac{15}{\pi^2} \approx 3.51982.
      \label{eq:planar-skew-kurt}
    \end{align}
  \item The left and right tails satisfy, for all $x \ge 0$,
    \begin{equation}\label{eq:planar-tails}
      \PP(\mathcal{L}_2 \le -x)
      = \frac{4}{3\pi}\, e^{-\frac{3\pi}{4} x},
      \qquad
      \PP(\mathcal{L}_2 \ge x)
      = \frac{12}{5\pi}\, e^{-\frac{5\pi}{4} x}
        - \frac{20}{21\pi}\, e^{-\frac{21\pi}{4} x}
        + O\bigl(e^{-\frac{45\pi}{4} x}\bigr).
    \end{equation}
\end{enumerate}
\end{theorem}

\begin{proof}
\ppart{(i)}
Lemma~\ref{lem:planar-offset-reduction} identifies $(\Regge(\RR^2; W, b),
U_{\ngen,\infty}(\RR^2))$ a.s.\ with $(\Regge(\RR^2; W, \tilde{b}),
U_{\ngen,\infty}(W))$ under synthetic i.i.d.\ offsets $\tilde{b}_h \sim \Norm(0,
1)$ independent of $W$. Propositions~\ref{prop:planar-spectrum}(ii)
and~\ref{prop:degenerate-ustat-limit} give $U_{\ngen,\infty}(\RR^2)
\xrightarrow{d} \sum_{m=0}^\infty \lambda_m(E_m - 1)$
via~\eqref{eq:ustat-exp-pairing}. Theorem~\ref{thm:second-chaos-transfer}
transfers this limit in $L^2(\PP)$ to $\frac{\ngen}{4}(2\pi -
\Theta_\infty(W))$, and Fubini's theorem with $\sum_{m=0}^\infty
\abs{\lambda_m} = \frac{8}{3\pi}$ and $\sum_{m=0}^\infty \lambda_m = 0$
yields~\eqref{eq:planar-limit-series}.

\ppart{(ii)}
With $z(s) := 2\sqrt{1 + s/\pi}$ and $(2m - 1)(2m + 3) = (2m + 1)^2 - 4$, we
have $1 - \lambda_m s = \frac{1 - z(s)^2/(2m+1)^2}{1 - 4/(2m+1)^2}$. Since
$\prod_{m=0}^M \bigl(1 - \frac{4}{(2m+1)^2}\bigr) = -\frac{2M+3}{2M+1} \to -1 =
\cos(\pi)$, dividing Euler's product $\cos(\frac{\pi z}{2}) = \prod_{m=0}^\infty
\bigl(1 - \frac{z^2}{(2m+1)^2}\bigr)$ at $z = z(s)$ by $-1$
gives~\eqref{eq:planar-mgf} between the nearest poles $s = -\frac{3\pi}{4} =
\lambda_0^{-1}$ and $s = \frac{5\pi}{4} = \lambda_1^{-1}$.

\ppart{(iii)}
Expanding $\delta(s) := \sqrt{\pi^2 + \pi s} - \pi = \frac{s}{2} -
\frac{s^2}{8\pi} + \frac{s^3}{16\pi^2} - \frac{5 s^4}{128\pi^3} + O(s^5)$ in
$\ln M_{\mathcal{L}_2}(s) = -\ln(\cos\delta(s)) = \frac{1}{2}\delta(s)^2 +
\frac{1}{12}\delta(s)^4 + O(s^5)$
gives~\eqref{eq:planar-cumulants}--\eqref{eq:planar-skew-kurt}.

\ppart{(iv)}
Write $\mathcal{L}_2 = \lambda_0 E_0 + S_+$, with $S_+ := \sum_{m=1}^\infty
\lambda_m E_m \ge 0$ independent of $E_0 \sim \mathrm{Exp}(1)$. For $x \ge 0$,
$\PP(\mathcal{L}_2 \le -x) = \PP(E_0 \ge \frac{3\pi}{4}(x + S_+)) =
e^{-\frac{3\pi}{4}x} M_{S_+}(-\frac{3\pi}{4})$. Since $M_{S_+}(s) = (1 -
\lambda_0 s) M_{\mathcal{L}_2}(s) = -\lambda_0 (s - s_0) M_{\mathcal{L}_2}(s)$
for $s < \frac{5\pi}{4}$, and $M_{\mathcal{L}_2}$ has a simple pole of residue $1$ at
$s_0 := \lambda_0^{-1} = -\frac{3\pi}{4}$ (the residue computation below
covers $m = 0$), $M_{S_+}(-\frac{3\pi}{4}) = -\lambda_0 =
\frac{4}{3\pi}$.

For $x > 0$, $M_{\mathcal{L}_2}$ has simple poles at $s_m := \lambda_m^{-1}$
($s_0 = -\frac{3\pi}{4}$, $s_1 = \frac{5\pi}{4}$, $s_2 = \frac{21\pi}{4}$,
$s_3 = \frac{45\pi}{4}$, \dots). Along $s = \sigma + it$ with $\sigma$ bounded
and $\abs{t} \to \infty$, $\sqrt{\pi^2 + \pi s} =
\sqrt{\frac{\pi\abs{t}}{2}}(1 + i\,\sgn t) + O(\abs{t}^{-1/2})$, so
$\abs{M_{\mathcal{L}_2}(\sigma + it)/(\sigma + it)} \lesssim
\abs{t}^{-1}e^{-\sqrt{\pi\abs{t}/2}}$ uniformly in $\sigma$. We may therefore
shift the Bromwich contour $\PP(\mathcal{L}_2 \ge x) = \frac{1}{2\pi
i}\int_{c - i\infty}^{c + i\infty} \frac{M_{\mathcal{L}_2}(s)}{s}e^{-sx}\,\dd s$
from $c \in (0, s_1)$ to $c' \in (s_3, s_4)$, picking up
$-\operatorname{Res}_{s = s_m}[M_{\mathcal{L}_2}(s)e^{-sx}/s]$ at $m = 1, 2,
3$, while the shifted integral is $O(e^{-c' x})$. At $s_m$ ($m \ge 0$) we have
$\sqrt{\pi^2 + \pi s_m} = (m + \frac{1}{2})\pi$ and $\frac{\dd}{\dd
s}\sqrt{\pi^2 + \pi s}\big|_{s_m} = \frac{1}{2m+1}$, so
$\operatorname{Res}_{s = s_m} M_{\mathcal{L}_2} = (-1)^{m}(2m+1)$ and the
residue at $s_m$ contributes $(-1)^{m+1}\frac{2m+1}{s_m}e^{-s_m x}$. For $m =
1$ this is $\frac{3}{s_1}e^{-s_1 x} = \frac{12}{5\pi}e^{-5\pi x/4}$, for $m = 2$
it is $-\frac{5}{s_2}e^{-s_2 x} = -\frac{20}{21\pi}e^{-21\pi x/4}$, and the
$m = 3$ term and the shifted integral are $O(e^{-45\pi x/4})$,
proving~\eqref{eq:planar-tails}.
\end{proof}

\subsection{Offset-free circle-integral proof and universality}
\label{subsec:direct-circle-integral}

We now prove a universal planar limit theorem directly from $\Theta_\infty(W)$,
without synthetic offsets, for centrally symmetric laws with finite fourth
moment. Let $u(\phi) := (\cos\phi, \sin\phi)\transp$, $u'(\phi) := (-\sin\phi,
\cos\phi)\transp$ and $G(W, \phi) := I_2 + \sum_{h=1}^{\ngen}
\mathbf{1}\{\inner{w_h}{u(\phi)} > 0\} w_h w_h\transp$. Then
\begin{equation}\label{eq:circle-integral}
  \Theta_\infty(W)
  = \int_0^{2\pi} \omega\bigl(G(W, \phi); \phi\bigr)\, \dd\phi,
  \qquad
  \omega(G; \phi)
  := \frac{\sqrt{\det G}}{u(\phi)\transp G\, u(\phi)}.
\end{equation}
To see this, note that on each open arc where $G(W, \phi) \equiv G_{(k)}(W)$ is
constant, $\omega(G; \phi)$ is the derivative of the polar angle $\psi_G(\phi)
= \allowbreak \phi + \allowbreak \arctan\bigl(\frac{u(\phi)\transp G\, u'(\phi)}{u(\phi)\transp G\,
u(\phi) + \sqrt{\det G}}\bigr)$ of the vector $G^{1/2}u(\phi)$. The increment
of $\psi_G$ across $[\phi_k, \phi_{k+1}]$ is the intrinsic opening angle
$\alpha_k^\infty(W)$, and $\psi_G(\phi + \pi) = \psi_G(\phi) + \pi$, so summing
over the $2\ngen$ arcs gives~\eqref{eq:circle-integral}.

\begin{lemma}
\label{lem:antipodal-expansion}
Let $W = [w_1, \dots, w_{\ngen}]\transp \in \RR^{\ngen \times 2}$ have pairwise
non-collinear rows. Define
\begin{gather*}
  \bar{G}(W) := I_2 + \frac{1}{2}\sum_{h=1}^{\ngen} w_h w_h\transp \in \mathrm{Sym}_2^{++}, \\
  \Delta(W, \phi) := G(W, \phi) - \bar{G}(W)
  = \frac{1}{2}\sum_{h=1}^{\ngen} \sgn(\inner{w_h}{u(\phi)})\, w_h w_h\transp, \\
  \Upsilon(W, \phi) := \bar{G}(W)^{-1/2} \Delta(W, \phi) \bar{G}(W)^{-1/2}.
\end{gather*}
\begin{enumerate}[label=(\roman*)]
  \item For a.e.\ $\phi \in [0, 2\pi)$, $\Delta(W, \phi + \pi) = -\Delta(W,
    \phi)$, $\Upsilon(W, \phi + \pi) = -\Upsilon(W, \phi)$, $\int_0^\pi
    \omega(\bar{G}(W); \phi)\,\dd\phi = \pi$ and $\int_0^{2\pi}
    \omega(\bar{G}(W); \phi)\,\dd\phi = 2\pi$.
  \item Let $G_0 \in \mathrm{Sym}_2^{++}$ satisfy $G_0 \pm \Delta(W, \phi) \succ
    0$ for a.e.\ $\phi \in [0, 2\pi)$, and set
    \[
      \Upsilon_0(\phi) := G_0^{-1/2}\Delta(W, \phi) G_0^{-1/2},
      \qquad
      \tilde{u}_0(\phi) := \frac{G_0^{1/2}u(\phi)}{\norm{G_0^{1/2}u(\phi)}}.
    \]
    Pairing $\phi \in [0, \pi)$ with
    $\phi + \pi$ cancels all odd-order terms in $\Upsilon_0(\phi)$ and gives
    \[
      2\pi - \int_0^{2\pi} \omega\bigl(G_0 + \Delta(W, \phi); \phi\bigr)\,\dd\phi
      = \mathcal{Q}_{\mathrm{cone}}(W; G_0) + \mathcal{R}_4(W; G_0),
    \]
    recovering $2\pi - \Theta_\infty(W) = \mathcal{Q}_{\mathrm{cone}}(W;
    \bar{G}(W)) + \mathcal{R}_4(W; \bar{G}(W))$ at $G_0 = \bar{G}(W)$, where
    \begin{equation}\label{eq:cone-quadratic-def}
      \mathcal{Q}_{\mathrm{cone}}(W; G_0)
      := \int_0^{2\pi}\!\! \omega(G_0; \phi)
         \Bigl[
           \frac{1}{4}\tr(\Upsilon_0^2)
           - \frac{1}{8}(\tr \Upsilon_0)^2
           + \frac{1}{2}(\tr \Upsilon_0)\, \tilde{u}_0\transp \Upsilon_0\tilde{u}_0
           - (\tilde{u}_0\transp \Upsilon_0\tilde{u}_0)^2
         \Bigr]\dd\phi,
    \end{equation}
    and there exists a universal constant $C > 0$ such that whenever $\sup_{\phi
    \in [0, 2\pi)} \norm{\Upsilon_0(\phi)}_{\mathrm{op}} \le 1/3$,
    \begin{equation}\label{eq:cone-quartic-rem}
      \abs{\mathcal{R}_4(W; G_0)}
      \le C\int_0^{2\pi} \omega(G_0; \phi)\, \norm{\Upsilon_0(\phi)}_{\mathrm{op}}^4\, \dd\phi.
    \end{equation}
\end{enumerate}
\end{lemma}

\begin{proof}
Since $u(\phi + \pi) = -u(\phi)$, we have $\Delta(W, \phi + \pi) = -\Delta(W,
\phi)$, $\Upsilon_0(\phi + \pi) = -\Upsilon_0(\phi)$ and
\[
  \omega\bigl(G_0 + \Delta(W, \phi); \phi\bigr)
  = \omega(G_0; \phi)\,
    \frac{\sqrt{\det(I_2 + \Upsilon_0(\phi))}}{1 + \tilde{u}_0(\phi)\transp \Upsilon_0(\phi)\tilde{u}_0(\phi)}.
\]
Integrating $\omega(G_0; \phi)$ over $[0, 2\pi)$ gives $2\pi$, and pairing $\phi$
with $\phi + \pi$ removes all odd powers of $\Upsilon_0(\phi)$. Substituting
$\det(I_2 + \Upsilon_0) = 1 + \tr \Upsilon_0 + \frac{1}{2}((\tr \Upsilon_0)^2 -
\tr(\Upsilon_0^2))$ into $(1 + q)^{-1}\sqrt{1 + t_1 + t_2} = 1 + (\frac{1}{2}t_1
- q) + \allowbreak (\frac{1}{2}t_2 - \frac{1}{8}t_1^2 - \frac{1}{2}t_1 q + q^2) +
\cdots$ at $q = \tilde{u}_0\transp \Upsilon_0 \tilde{u}_0$
gives~\eqref{eq:cone-quadratic-def}. When $\norm{\Upsilon_0}_{\mathrm{op}} \le
1/3$, the even function
\[
  f(t) := \frac{1}{2}\Bigl(\frac{\sqrt{\det(I_2 + t \Upsilon_0)}}{1 + t q}
  + \frac{\sqrt{\det(I_2 - t \Upsilon_0)}}{1 - t q}\Bigr)
  = 1 + c_2 t^2 + \sum_{k=2}^\infty c_{2k} t^{2k}
\]
satisfies $\sup_{\abs{t} = (2\norm{\Upsilon_0}_{\mathrm{op}})^{-1}}\abs{f(t)} \le
3$, so Cauchy's estimate gives $\bigl|\sum_{k=2}^\infty c_{2k}\bigr| \le
C\norm{\Upsilon_0}_{\mathrm{op}}^4$ and~\eqref{eq:cone-quartic-rem}.
\end{proof}

\begin{lemma}
\label{lem:cone-ustat}
Let $W = [w_1, \dots, w_{\ngen}]\transp \in \RR^{\ngen \times 2}$ have pairwise
non-collinear rows, and let $G_0 \in \mathrm{Sym}_2^{++}$.
\textup{(i)} For every $A \in \mathrm{GL}_2(\RR)$, $\mathcal{Q}_{\mathrm{cone}}(W
A; A\transp G_0 A) = \mathcal{Q}_{\mathrm{cone}}(W; G_0)$.
\textup{(ii)} Writing $\tilde{w}_h := G_0^{-1/2} w_h =
\norm{\tilde{w}_h}(\cos\tilde{\phi}_h, \sin\tilde{\phi}_h)\transp$ and expanding
$\Upsilon_0(\phi) = \frac{1}{2}\sum_{h=1}^{\ngen}
\sgn(\inner{\tilde{w}_h}{\tilde{u}_0(\phi)}) \tilde{w}_h \tilde{w}_h\transp$
in~\eqref{eq:cone-quadratic-def}, every diagonal term ($c = c'$) vanishes and
every off-diagonal pair integrates to $\Zlim(\tilde{w}_c, \tilde{w}_{c'})$ (at
$\sigw = 1$):
\begin{equation}\label{eq:cone-exact-ustat}
  \begin{gathered}
    \mathcal{Q}_{\mathrm{cone}}(W; G_0)
    = -\frac{1}{4}\sum_{c < c'}
       \norm{\tilde{w}_c}^2 \norm{\tilde{w}_{c'}}^2
       \cos(\tilde{\phi}_c - \tilde{\phi}_{c'}) \abs{\sin(\tilde{\phi}_c - \tilde{\phi}_{c'})}, \\
    \frac{\ngen}{4} \mathcal{Q}_{\mathrm{cone}}\!\left(W;\, \frac{\ngen \sigw^2}{2} I_2\right)
    = U_{\ngen,\infty}(\RR^2).
  \end{gathered}
\end{equation}
\end{lemma}

\begin{proof}
Part~\textup{(i)} follows from the substitution $u(\phi) \mapsto
A^{-1}u(\phi)/\norm{A^{-1}u(\phi)}$, reducing~\textup{(ii)} to $G_0 = I_2$
($\omega(I_2; \phi) \equiv 1$, $\tilde{u}_0 = u$). Write $v_c := (\cos\phi_c,
\sin\phi_c)\transp$, $\alpha_c := \phi - \phi_c$ and $s_c(\phi) :=
\sgn(\inner{w_c}{u(\phi)})$, so $\Upsilon_0(\phi) =
\frac{1}{2}\sum_{c=1}^{\ngen} s_c(\phi)\norm{w_c}^2 v_c v_c\transp$. In the
four terms of~\eqref{eq:cone-quadratic-def}, the symmetrised $(c, c')$
coefficient multiplying $\norm{w_c}^2 \norm{w_{c'}}^2 s_c(\phi)s_{c'}(\phi)$ is
\begin{align*}
  &\frac{\cos^2(\phi_c - \phi_{c'})}{16}
   - \frac{1}{32}
   + \frac{\cos^2\alpha_c + \cos^2\alpha_{c'}}{16}
   - \frac{\cos^2\alpha_c \cos^2\alpha_{c'}}{4} \\
  &\qquad = \frac{1 + \cos(2(\phi_c - \phi_{c'}))}{32}
   - \frac{1}{32}
   + \frac{2 + \cos 2\alpha_c + \cos 2\alpha_{c'}}{32}
   - \frac{(1 + \cos 2\alpha_c)(1 + \cos 2\alpha_{c'})}{16} \\
  &\qquad = -\frac{1}{32}\bigl(\cos 2\alpha_c + \cos 2\alpha_{c'}
   + \cos(2(2\phi - \phi_c - \phi_{c'}))\bigr)
   = -\frac{1}{64}\, F_{cc'}'(\phi),
\end{align*}
since $\cos 2\alpha_c \cos 2\alpha_{c'} = \frac{1}{2}\cos(2(\phi_c - \phi_{c'}))
+ \frac{1}{2}\cos(2(2\phi - \phi_c - \phi_{c'}))$ cancels the
$\frac{1}{32}\cos(2(\phi_c - \phi_{c'}))$ term, where $F_{cc'}(\phi) := \sin
2\alpha_c + \sin 2\alpha_{c'} + \frac{1}{2}\sin(2(2\phi - \phi_c - \phi_{c'}))$
is $\pi$-periodic. Summing over $c, c'$ leaves
\[
  \mathcal{Q}_{\mathrm{cone}}(W; I_2)
  = -\frac{1}{64}\sum_{c, c' = 1}^{\ngen} \norm{w_c}^2 \norm{w_{c'}}^2
    \int_0^{2\pi} s_c(\phi) s_{c'}(\phi)\, F_{cc'}'(\phi)\, \dd\phi.
\] For $c = c'$,
$s_c(\phi)^2 = 1$ a.e.\ and $\int_0^{2\pi} F_{cc}'(\phi)\,\dd\phi = 0$. For $c
\ne c'$, rotate so that $\phi_c = 0$ and $\phi_{c'} = \vartheta \in (0, \pi)$.
The product $s_c(\phi)s_{c'}(\phi) = \sgn(\cos\phi)\sgn(\cos(\phi - \vartheta))$
is $\pi$-periodic, as is $F_{cc'}$, so it suffices to integrate over $[0, \pi)$
and double. Suppose first that $\vartheta < \pi/2$. On $[0, \pi)$ the product
changes sign exactly at $\pi/2$ (where $\cos\phi$ vanishes) and at $\vartheta +
\pi/2$ (where $\cos(\phi - \vartheta)$ vanishes), taking the values $+1$ on
$[0, \pi/2)$, $-1$ on $(\pi/2, \vartheta + \pi/2)$ and $+1$ on $(\vartheta +
\pi/2, \pi)$. Hence, writing $F := F_{cc'}$,
\begin{align*}
  \int_0^{2\pi} s_c s_{c'} F'\,\dd\phi
  &= 2\Bigl[\bigl(F(\tfrac{\pi}{2}) - F(0)\bigr)
    - \bigl(F(\vartheta + \tfrac{\pi}{2}) - F(\tfrac{\pi}{2})\bigr)
    + \bigl(F(\pi) - F(\vartheta + \tfrac{\pi}{2})\bigr)\Bigr] \\
  &= 4\bigl[F(\tfrac{\pi}{2}) - F(\vartheta + \tfrac{\pi}{2})\bigr],
\end{align*}
using $F(\pi) = F(0)$. Direct evaluation gives $F(\pi/2) =
\frac{1}{2}\sin 2\vartheta$ and $F(\vartheta + \pi/2) = -\frac{1}{2}\sin
2\vartheta$, so the integral equals $4\sin 2\vartheta =
8\cos\vartheta\sin\vartheta$. When $\vartheta > \pi/2$ the second sign change
occurs at $\vartheta - \pi/2 < \pi/2$ instead, and the same telescoping gives
$4[F(\pi/2) - F(\vartheta - \pi/2)] = 8\cos\vartheta\sin\vartheta$
again. Substituting into the double sum and recalling $\sin\vartheta > 0$
proves~\eqref{eq:cone-exact-ustat}.
\end{proof}

\begin{theorem}[Planar universality]
\label{thm:planar-universality}
Let $\dimin = 2$, and let $w_1, \dots, w_{\ngen}$ be i.i.d.\ in $\RR^2$ with a
centrally symmetric law $P_w$ ($-w_1 \stackrel{d}{=} w_1$) charging no line
through the origin, with $\Sigma_w := \EE[w_1 w_1\transp] \in
\mathrm{Sym}_2^{++}$ and $\EE[\norm{w_1}^4] < \infty$. Set $\tilde{w}_h :=
\Sigma_w^{-1/2} w_h$ and $\tau_2(P_w) = \frac{1}{4}\EE[\norm{\tilde{w}_1}_2^4]
\in [1, \infty)$. Then as $\ngen \to \infty$,
\begin{equation}\label{eq:univ-diff-prob}
  \frac{\ngen}{4}\bigl(2\pi - \Theta_\infty(W)\bigr)
    - U_{\ngen,\infty}(\widetilde{W})
  = \frac{\ngen}{4}\bigl(2\pi - \Theta_\infty(W)\bigr)
    - \frac{1}{\ngen}\sum_{1 \le c < c' \le \ngen} \Zlim(\tilde{w}_c, \tilde{w}_{c'})
  \xrightarrow[\ngen \to \infty]{\PP} 0,
\end{equation}
where $\Zlim(\tilde{w}_c, \tilde{w}_{c'}) =
-\frac{1}{4}\inner{\tilde{w}_c}{\tilde{w}_{c'}}\norm{\tilde{w}_c \wedge
\tilde{w}_{c'}}$. If $P_w$ also satisfies~\ref{assump:A3}, then
\begin{equation}\label{eq:univ-limit}
  \frac{\ngen}{4}\bigl(2\pi - \Theta_\infty(W)\bigr)
  \xrightarrow[\ngen \to \infty]{d}
  \frac{\tau_2(P_w)}{2}\, \mathcal{L}_2.
\end{equation}
Under the planar general-position condition of
Lemma~\ref{lem:planar-offset-reduction}, both statements hold for
$\frac{\ngen}{4}\Regge(\RR^2)$.
\end{theorem}

\begin{proof}
The proof has four parts: we standardise the normals and expand
$\Theta_\infty(W)$ via Lemma~\ref{lem:antipodal-expansion}; we show that the
expansion is valid off an event of vanishing probability; we show that the
quartic remainder is negligible; and we identify the quadratic term with the
U-statistic.

\ppart{Standardisation.}
Writing $w_h = \Sigma_w^{1/2}\tilde{w}_h$ ($\Sigma_{\tilde{w}} = I_2$)
in~\eqref{eq:circle-integral} gives $G(W, \phi) =
\Sigma_w^{1/2}\bigl(G_\Sigma(\widetilde{W}) + \Delta(\widetilde{W},
\tilde{\phi})\bigr)\Sigma_w^{1/2}$, where $u(\tilde{\phi}) :=
\Sigma_w^{1/2}u(\phi)/\norm{\Sigma_w^{1/2}u(\phi)}_2$ and $G_{\Sigma}(\widetilde{W})
:= \Sigma_w^{-1} + \frac{1}{2}\sum_{h=1}^{\ngen}\tilde{w}_h\tilde{w}_h\transp$.
For $A := \Sigma_w^{1/2}$ and $K \in \mathrm{Sym}_2^{++}$, the polar
decomposition $(AKA)^{1/2} = O K^{1/2} A$ ($O \in \mathrm{O}(2)$) implies
$\omega(AKA; \phi)\,\dd\phi = \omega(K; \tilde{\phi})\,\dd\tilde{\phi}$, and so
$\Theta_\infty(W) = \int_0^{2\pi}\omega(G_\Sigma(\widetilde{W}) +
\Delta(\widetilde{W}, \tilde{\phi}); \tilde{\phi})\,\dd\tilde{\phi}$. Applying
Lemma~\ref{lem:antipodal-expansion}(ii) at $G_0 = G_\Sigma(\widetilde{W})$ yields
$2\pi - \Theta_\infty(W) = \mathcal{Q}_{\mathrm{cone}}(\widetilde{W};
G_\Sigma(\widetilde{W})) + \mathcal{R}_4(\widetilde{W}; G_\Sigma(\widetilde{W}))$.

By central symmetry, write $\tilde{w}_h = \epsilon_h \mathring{w}_h$ with i.i.d.\
Rademacher signs $\epsilon_h \in \{\pm 1\}$ independent of $\mathring{W} =
(\mathring{w}_1, \dots, \mathring{w}_{\ngen})$, and set $\mathsf{E}_{\ngen} :=
\frac{2}{\ngen}G_\Sigma(\widetilde{W}) - I_2 =
\frac{1}{\ngen}\sum_{h=1}^{\ngen}(\mathring{w}_h\mathring{w}_h\transp - I_2) +
\frac{2}{\ngen}\Sigma_w^{-1}$ and $\Upsilon_\Sigma(\widetilde{W}, \phi) :=
G_\Sigma(\widetilde{W})^{-1/2}\Delta(\widetilde{W},
\phi)G_\Sigma(\widetilde{W})^{-1/2}$. For $\ngen \ge
8\norm{\Sigma_w^{-1}}_{\mathrm{op}}$, define
\begin{align*}
  \mathcal{E}_{\mathrm{bulk}}
  &:= \Bigl\{\frac{1}{\ngen}\sum_{h=1}^{\ngen} \norm{\mathring{w}_h}_2^4 \le 8\tau_2(P_w)\Bigr\}
     \cap \Bigl\{\norm{\mathsf{E}_{\ngen}}_{\mathrm{op}} \le \frac{1}{2}\Bigr\}, \\
  \mathcal{G}_{\ngen}(W)
  &:= \mathcal{E}_{\mathrm{bulk}} \cap \Bigl\{\sup_{\phi \in [0, 2\pi)} \norm{\Upsilon_\Sigma(\widetilde{W}, \phi)}_{\mathrm{op}} \le \frac{1}{3}\Bigr\}.
\end{align*}
On $\mathcal{E}_{\mathrm{bulk}}$, $\frac{2}{\ngen}G_\Sigma = I_2 +
\mathsf{E}_{\ngen}$ has eigenvalues in $[\frac{1}{2}, \frac{3}{2}]$, so
$\lambda_{\min}(G_\Sigma(\widetilde{W})) \ge \frac{\ngen}{4}$. Under
$\EE[\norm{w_1}^4] < \infty$, Khintchine's weak law of large numbers gives
$\mathsf{E}_{\ngen} \xrightarrow{\PP} 0$ and
$\frac{1}{\ngen}\sum_{h=1}^{\ngen}\norm{\mathring{w}_h}_2^4 \xrightarrow{\PP}
4\tau_2(P_w)$, whence $\PP(\mathcal{E}_{\mathrm{bulk}}^c) \to 0$.

\ppart{Uniform operator-norm control.}
To see that $\PP(\mathcal{G}_{\ngen}(W)^c) \to 0$ as well, we need control of
$\Upsilon_\Sigma$ uniformly in $\phi$. Any $\Upsilon \in
\mathrm{Sym}_2(\RR)$ is determined linearly by its quadratic
forms $q_i := v_i\transp \Upsilon v_i$ along $v_i := u(\phi_i)$ for $\phi_i \in
\{0, \pi/3, 2\pi/3\}$ ($\Upsilon_{11} = q_1$, $\Upsilon_{12} = (q_2 -
q_3)/\sqrt{3}$, $\Upsilon_{22} = \frac{2}{3}(q_2 + q_3) - \frac{1}{3}q_1$). Hence
$\norm{\Upsilon}_{\mathrm{op}} \le C_0\max_{i=1,2,3}\abs{q_i}$ for a universal
constant $C_0 \ge 1$.

Consequently, $\max_i\abs{q_i} \le (3C_0)^{-1}$ forces
$\norm{\Upsilon}_{\mathrm{op}} \le 1/3$. For fixed $\phi$ and $i$, $v_i\transp
\Upsilon_\Sigma(\widetilde{W}, \phi)v_i = \sum_{h=1}^{\ngen} \epsilon_h a_h$ with
$a_h := \frac{1}{2}\sgn(\ip{\mathring{w}_h}{u(\phi)})\ip{\mathring{w}_h}{G_\Sigma^{-1/2}v_i}^2$.
On $\mathcal{E}_{\mathrm{bulk}}$, $\sum_{h=1}^{\ngen} a_h^2 \le
C\lambda_{\min}(G_\Sigma)^{-2}\sum_{h=1}^{\ngen}\norm{\mathring{w}_h}_2^4 \le
C\tau_2(P_w)/\ngen$, so Hoeffding's inequality and a union bound over the
$2\ngen$ arcs of $\Sph{1}$ cut out by $\ip{\mathring{w}_h}{u} = 0$ and the $3$
directions $v_i$ give $\PP(\mathcal{E}_{\mathrm{bulk}} \setminus
\mathcal{G}_{\ngen}(W) \mid \mathring{W})\mathbf{1}_{\mathcal{E}_{\mathrm{bulk}}}
\le C\,\ngen\exp(-c\,\ngen/\tau_2(P_w))$, whence $\PP(\mathcal{G}_{\ngen}(W)^c)
\to 0$.

\ppart{Quartic remainder.}
On $\mathcal{E}_{\mathrm{bulk}}$ we have $\norm{\Upsilon_\Sigma(\widetilde{W},
\phi)}_{\mathrm{op}} \le C\norm{\ngen^{-1}\Delta(\widetilde{W},
\phi)}_{\mathrm{F}}$ and $\omega(G_\Sigma;\phi) \le \sqrt{\det
G_\Sigma}/\lambda_{\min}(G_\Sigma) \le C$ uniformly in $\phi$. Conditionally on
$\mathring{W}$, each entry of $\ngen^{-1}\Delta(\widetilde{W}, \phi)$ is a
Rademacher sum, so Khintchine's inequality gives
\[
  \EE_\epsilon\Bigl[\norm{\ngen^{-1}\Delta(\widetilde{W}, \phi)}_{\mathrm{F}}^8 \;\Big|\; \mathring{W}\Bigr]\mathbf{1}_{\mathcal{E}_{\mathrm{bulk}}}
  \le C\Bigl(\frac{1}{\ngen^2}\sum_{h=1}^{\ngen}\norm{\mathring{w}_h}_2^4\Bigr)^4\mathbf{1}_{\mathcal{E}_{\mathrm{bulk}}}
  \le C\tau_2(P_w)^4\ngen^{-4}
\]
uniformly in $\phi$. By~\eqref{eq:cone-quartic-rem} at $G_0 =
G_\Sigma(\widetilde{W})$ and Minkowski's integral inequality, we obtain
$\bigl\lVert\frac{\ngen}{4}\mathcal{R}_4(\widetilde{W}; \allowbreak
G_\Sigma(\widetilde{W}))\allowbreak\mathbf{1}_{\mathcal{G}_{\ngen}(W)}\bigr\rVert_{L^2(\PP)}
\le \allowbreak C\tau_2(P_w)^2\ngen^{-1}$. Thus
$\frac{\ngen}{4}\abs{\mathcal{R}_4(\widetilde{W}; \allowbreak
G_\Sigma(\widetilde{W}))} \xrightarrow{\PP} 0$.

\ppart{Quadratic term.}
Set $M := (I_2 + \mathsf{E}_{\ngen})^{-1}$. By Lemma~\ref{lem:cone-ustat}(ii),
\[
  \frac{\ngen}{4}\mathcal{Q}_{\mathrm{cone}}\bigl(\widetilde{W}; G_\Sigma(\widetilde{W})\bigr)
  = \frac{1}{\ngen}\sum_{1 \le c < c' \le \ngen}\epsilon_c \epsilon_{c'}\,\Zlim\bigl(M^{1/2}\mathring{w}_c, M^{1/2}\mathring{w}_{c'}\bigr),
\]
where $\Zlim(M^{1/2}y, M^{1/2}y') = -\frac{1}{4}(y\transp M y')\sqrt{\det
M}\,\norm{y \wedge y'}$. On $\mathcal{E}_{\mathrm{bulk}}$,
$\norm{\mathsf{E}_{\ngen}}_{\mathrm{op}} \le \frac{1}{2}$, and $\mathsf{E}
\mapsto (I_2 + \mathsf{E})^{-1}\det(I_2 + \mathsf{E})^{-1/2}$ is Lipschitz on
$\{\norm{\mathsf{E}}_{\mathrm{op}} \le \frac{1}{2}\}$ with value $I_2$ at $0$.
Hence $\norm{M\sqrt{\det M} - I_2}_{\mathrm{op}} \le
C\norm{\mathsf{E}_{\ngen}}_{\mathrm{op}}$ and $\abs{\Zlim(M^{1/2}y, M^{1/2}y') -
\Zlim(y, y')} \le C\norm{\mathsf{E}_{\ngen}}_{\mathrm{op}}\norm{y}_2^2\norm{y'}_2^2$.
Conditional orthogonality of $\{\epsilon_c\epsilon_{c'}\}_{c < c'}$ given
$\mathring{W}$ yields
\[
  \EE_\epsilon\Bigl[\Bigl(\frac{\ngen}{4}\mathcal{Q}_{\mathrm{cone}}\bigl(\widetilde{W}; G_\Sigma(\widetilde{W})\bigr) - U_{\ngen,\infty}(\widetilde{W})\Bigr)^2 \;\Big|\; \mathring{W}\Bigr]\mathbf{1}_{\mathcal{E}_{\mathrm{bulk}}}
  \le C\tau_2(P_w)^2\norm{\mathsf{E}_{\ngen}}_{\mathrm{op}}^2\mathbf{1}_{\mathcal{E}_{\mathrm{bulk}}}.
\]
Since $\norm{\mathsf{E}_{\ngen}}_{\mathrm{op}}\mathbf{1}_{\mathcal{E}_{\mathrm{bulk}}}
\le \frac{1}{2}$ and $\mathsf{E}_{\ngen} \xrightarrow{\PP} 0$, dominated
convergence gives
$\bigl(\frac{\ngen}{4}\mathcal{Q}_{\mathrm{cone}}(\widetilde{W}; \allowbreak
G_\Sigma(\widetilde{W})) - \allowbreak
U_{\ngen,\infty}(\widetilde{W})\bigr)\mathbf{1}_{\mathcal{E}_{\mathrm{bulk}}}
\to 0$ in $L^2(\PP)$, proving~\eqref{eq:univ-diff-prob} (with probability rate
$O_\varepsilon(\ngen^{-1})$ by Chebyshev's inequality whenever
$\EE[\norm{w_1}^8] < \infty$). Proposition~\ref{prop:planar-spectrum}(ii) and
Proposition~\ref{prop:degenerate-ustat-limit} then give~\eqref{eq:univ-limit}.
\end{proof}

\subsection{Quantitative Kolmogorov rates and finite-\texorpdfstring{$\ngen$}{H} moment asymptotics}
\label{subsec:planar-rates-moments}

\begin{corollary}
\label{cor:planar-kolmogorov}
Let $\dimin = 2$ and $\Omega = \RR^2$. Suppose $w_1, \dots, w_{\ngen}
\stackrel{\mathrm{i.i.d.}}{\sim} \Norm(0, \sigw^2 I_2)$ with fixed $\sigw > 0$,
and that $b \in \RR^{\ngen}$ satisfies the planar general-position condition of
Lemma~\ref{lem:planar-offset-reduction}. Set
$S_{\ngen} := \frac{\ngen}{4}\Regge(\RR^2) = \frac{\ngen}{4}(2\pi -
\Theta_\infty(W))$ and $U_{\ngen,\infty} := U_{\ngen,\infty}(\RR^2)$. There
exist universal $C > 0$ and $C_{\sigw} > 0$ depending only on $\sigw$ such that
\begin{align}
  d_{\mathrm{Kol}}(U_{\ngen,\infty},\, \mathcal{L}_2)
  &\le C\, \ngen^{-1}
  && (\forall \ngen \ge 2),
  \label{eq:planar-U-kolmogorov} \\
  d_{\mathrm{Kol}}\bigl(S_{\ngen},\, \mathcal{L}_2\bigr)
  &\le C_{\sigw}\, \ngen^{-1/3}
  && (\forall \ngen \ge 4).
  \label{eq:planar-S-kolmogorov}
\end{align}
\end{corollary}

The two exponents measure different things. The $O(\ngen^{-1})$ rate
in~\eqref{eq:planar-U-kolmogorov} is the quadratic-form Berry--Esseen rate for
the U-statistic alone, and is optimal in general. The $O(\ngen^{-1/3})$ rate
in~\eqref{eq:planar-S-kolmogorov} arises from smoothing the $O(\ngen^{-1/2})$
mean-square distance between $S_{\ngen}$ and $U_{\ngen,\infty}$ (sharp by
Theorem~\ref{thm:second-chaos-transfer}) against the bounded density of
$\mathcal{L}_2$. This uses Chebyshev's inequality at the scale $\eta =
\ngen^{-1/3}$ that balances the two error terms. It is not claimed to be
sharp. Under the bare fourth-moment hypothesis of
Theorem~\ref{thm:planar-universality} we only have convergence in probability,
and no rate.

\begin{proof}
For Gaussian normals, $\tau_2 = 2$ and $\EE[\varrho_1^4] = 24$, so
Theorem~\ref{thm:hilbert-kolmogorov}(ii) gives~\eqref{eq:planar-U-kolmogorov}.
Since $\abs{\varphi_{\mathcal{L}_2}(t)} \le \prod_{m=0}^2 (1 + \lambda_m^2
t^2)^{-1/2} = O(\abs{t}^{-3}) \in L^1(\RR)$, Fourier inversion bounds the
density of $\mathcal{L}_2$ by $M_{\mathrm{dens}} \le
\frac{1}{2\pi}\int_{\RR}\abs{\varphi_{\mathcal{L}_2}(t)}\,\dd t < \infty$.
Adjoining synthetic offsets via Lemma~\ref{lem:planar-offset-reduction} puts us
in the setting of Theorem~\ref{thm:second-chaos-transfer}(ii), which gives
$\delta_{\ngen} := \norm{S_{\ngen} - U_{\ngen,\infty}}_{L^2(\PP)} \le
C_{\sigw}'\ngen^{-1/2}$. Applying Chebyshev's
inequality at $x \pm \eta$ with $\eta = \delta_{\ngen}^{2/3}$ gives
$d_{\mathrm{Kol}}(S_{\ngen}, \mathcal{L}_2) \le d_{\mathrm{Kol}}(U_{\ngen,\infty},
\mathcal{L}_2) + M_{\mathrm{dens}}\eta + \eta^{-2}\delta_{\ngen}^2$,
proving~\eqref{eq:planar-S-kolmogorov}.
\end{proof}

\begin{corollary}
\label{cor:planar-moments}
Let $\dimin = 2$.
\textup{(i)} If $w_1, \dots, w_{\ngen}$ are independent and centrally symmetric
($-w_h \stackrel{d}{=} w_h$), each charging no line through the origin (with no
moment condition and possibly distinct laws across $h$), then $\EE[\Theta_\infty(W)]
= 2\pi$ and $\EE[2\pi - \Theta_\infty(W)] = 0$ for every $\ngen \ge 2$. Although
the chamber-varying pullback metrics $G_{(k)}(W)$ make $\Theta_\infty(W) \ne
2\pi$ pathwise, orientation symmetry restores in expectation the deterministic
angle sum $2\pi$ of a fixed-metric Euclidean conical tessellation~\cite{HugSchneider2016}.
\textup{(ii)} For $w_1, \dots, w_{\ngen} \stackrel{\mathrm{i.i.d.}}{\sim}
\Norm(0, \sigw^2 I_2)$ with fixed $\sigw > 0$,
\begin{equation}\label{eq:planar-variance}
  \Var\bigl(\Theta_\infty(W)\bigr)
  = \EE\bigl[(2\pi - \Theta_\infty(W))^2\bigr]
  = \frac{4}{\ngen^2} + O(\ngen^{-3}).
\end{equation}
\end{corollary}

\begin{proof}
For~(i), adjoining synthetic $\tilde{b}_h \stackrel{\mathrm{i.i.d.}}{\sim}
\Norm(0, 1)$ via Lemma~\ref{lem:planar-offset-reduction} and applying
Proposition~\ref{thm:first-chaos-annihilation}(i) give $\EE[2\pi -
\Theta_\infty(W) \mid w_j] = 0$ a.s.\ and $\EE[2\pi - \Theta_\infty(W)] = 0$
(since $0 < \Theta_\infty(W) < 2\pi\ngen$). For~(ii),
Theorem~\ref{thm:second-chaos-transfer}(i) and Lemmas~\ref{lem:cross-bone-es}
and~\ref{lem:kernel-rate} give $\EE[S_{\ngen}^2] = \norm{U_{\ngen,\infty} +
\mathcal{D}_2}_{L^2}^2 + \norm{\mathcal{D}_{\ge 3}}_{L^2}^2 = \frac{1}{4} +
O(\ngen^{-1})$, proving~\eqref{eq:planar-variance}.
\end{proof}
\section{All-dimensional Funk--Hecke spectrum, Berry--Esseen rate and parity orthogonality}
\label{sec:high-dim}

For $\dimin \ge 3$, the slice volume $\ell_{cc'}(\Omega) = \Haus{\dimin -
2}(B_{cc'} \cap \Omega)$ in Theorem~\ref{thm:second-chaos-transfer}(iii) couples
$(w_c, w_{c'})$ to $(b_c, b_{c'})$. The exception is when $b_h \equiv 0$ and $\Omega = B(0, R)$,
where every bone passes through the origin with deterministic slice volume
$\omega_{\dimin - 2} R^{\dimin - 2}$ (recall $\omega_k =
\pi^{k/2}/\Gamma(k/2 + 1)$, $\omega_0 = 1$; see Section~\ref{sec:remarks} for
general windows). We diagonalise $T_{\Zlim}$ on $\RR^{\dimin}$ via the
Funk--Hecke theorem, derive an $O(\ngen^{-1})$ Kolmogorov rate for
$U_{\ngen,\infty}(W)$, and prove finite-$\ngen$ parity orthogonality and
asymptotic independence from unoriented tessellation functionals across all
$\dimin \ge 2$.

\subsection{Funk--Hecke diagonalisation and ball-window limits}
\label{subsec:high-dim-extensions}

\begin{theorem}[Funk--Hecke spectrum]
\label{thm:funk-hecke}
Let $\dimin \ge 2$, and let $w = \norm{w} u \in \RR^{\dimin} \setminus \{0\}$
($u := w/\norm{w} \in \Sph{\dimin - 1}$) have an isotropic distribution $P_w$
satisfying~\ref{assump:A3} with $\Sigma_w = \sigw^2 I_{\dimin}$, normalised
radial variable $\varrho(w) := \frac{\norm{w}^2}{\dimin \sigw^2}$ and
$\tau_2(P_w) = \EE[\varrho(w)^2] < \infty$ ($\tau_2 = \frac{\dimin + 2}{\dimin}$
under~\ref{assump:A2}). The asymptotic curvature kernel~\eqref{eq:intro-Z-infinity}
factors as
\[
  \Zlim(w, w')
  = -\frac{\norm{w}^2 \norm{w'}^2}{4\sigw^4}\,\ip{u}{u'}\sqrt{1 - \ip{u}{u'}^2}
  = \left(\frac{\dimin^2}{4}\right) \varrho(w)\varrho(w') K_{\mathrm{sph}}(\ip{u}{u'}),
\]
where $K_{\mathrm{sph}}(t) := -t\sqrt{1 - t^2}$ for $t \in [-1, 1]$
($K_{\mathrm{ang}}(\alpha) = K_{\mathrm{sph}}(\cos\alpha)$ at $\dimin = 2$).
\begin{enumerate}[label=(\roman*)]
  \item With $\varsigma_{\dimin - 1} := \sigma_{\dimin - 1}^{-1}\Haus{\dimin -
    1}|_{\Sph{\dimin - 1}}$ denoting uniform measure on $\Sph{\dimin - 1}$, the
    spherical convolution $(T_{K_{\mathrm{sph}}} f)(u) := \int_{\Sph{\dimin - 1}}
    K_{\mathrm{sph}}(\ip{u}{v}) f(v)\,\dd\varsigma_{\dimin - 1}(v)$ has spherical
    harmonics of degree $\ell \ge 0$ as eigenfunctions, with eigenvalues
    $\kappa_{2m}(\dimin) = 0$ ($m \ge 0$) and, for odd degree $\ell = 2m + 1$
    ($m \ge 0$),
    \begin{equation}\label{eq:funk-hecke-eigenvalues}
      \kappa_{2m+1}(\dimin)
      = \frac{\Gamma(\dimin/2)^2\,\Gamma(m - \tfrac{1}{2})\,\Gamma(m + \tfrac{3}{2})}{2\pi\,\Gamma\bigl(m + \tfrac{\dimin - 1}{2}\bigr)\,\Gamma\bigl(m + \tfrac{\dimin + 3}{2}\bigr)},
      \qquad
      \kappa_1(\dimin)
      = -\frac{\Gamma(\dimin/2)^2}{2\,\Gamma\bigl(\frac{\dimin - 1}{2}\bigr)\,\Gamma\bigl(\frac{\dimin + 3}{2}\bigr)} < 0,
    \end{equation}
    while $\kappa_{2m+1}(\dimin) > 0$ for all $m \ge 1$ ($\ell \ge 3$).
  \item Let $\{Y_{2m+1, j}\}_{j=1}^{N(\dimin, 2m+1)}$ be an
    $L^2(\Sph{\dimin - 1}, \varsigma_{\dimin - 1})$-orthonormal basis of
    degree-$(2m + 1)$ spherical harmonics, with $N(\dimin, \ell) = \frac{2\ell +
    \dimin - 2}{\ell}\binom{\ell + \dimin - 3}{\ell - 1}$. Then $\psi_{2m+1,
    j}(w) := \frac{\varrho(w)}{\sqrt{\tau_2(P_w)}} Y_{2m+1, j}(u)$ are
    orthonormal in $L^2(\RR^{\dimin}, P_w)$, and the non-zero eigenvalues of
    $(T_{\Zlim} f)(w) := \EE_{w'}[\Zlim(w, w') f(w')]$ are $\Lambda_{2m+1,
    j}(\dimin) = \bigl(\frac{\dimin^2 \tau_2(P_w)}{4}\bigr)\kappa_{2m+1}(\dimin)$
    ($1 \le j \le N(\dimin, 2m+1)$). In particular, $T_{\Zlim}$ is trace class
    with $\sum_{m=0}^\infty N(\dimin, 2m+1)\abs{\kappa_{2m+1}(\dimin)} < \infty$,
    $\tr(T_{\Zlim}) = 0$ and signature
    \begin{equation}\label{eq:funk-hecke-signature}
      \operatorname{sig}(T_{\Zlim}) = (\dimin_-, \dimin_+) = (\dimin, \infty)
      \qquad (\forall \dimin \ge 2).
    \end{equation}
  \item For $(w_1, \dots, w_{\ngen})$ i.i.d.\ from isotropic $P_w$
    satisfying~\ref{assump:A3}, $b_h \equiv 0$ and $\Omega = B(0, R) \subset
    \RR^{\dimin}$ ($R > 0$, using~\eqref{eq:model-apex-convention} at $\dimin =
    2$), every bone slice has volume $\ell_{cc'}(B(0, R)) = \omega_{\dimin - 2}
    R^{\dimin - 2}$. With $U_{\ngen,\infty}(W) := \frac{1}{\ngen}\sum_{1 \le c <
    c' \le \ngen} \Zlim(w_c, w_{c'})$ and $\gamma_{m,j}
    \stackrel{\mathrm{i.i.d.}}{\sim} \Norm(0, 1)$,
    \[
      \omega_{\dimin - 2} R^{\dimin - 2}\, U_{\ngen,\infty}(W)
      \xrightarrow[\ngen \to \infty]{d}
      \omega_{\dimin - 2} R^{\dimin - 2}\,\frac{\dimin^2 \tau_2(P_w)}{8}
      \sum_{m=0}^\infty \kappa_{2m+1}(\dimin) \sum_{j=1}^{N(\dimin, 2m+1)} (\gamma_{m,j}^2 - 1),
    \]
    and $\frac{\ngen}{4}\Regge(B(0, R))$ shares this limit law when $\dimin \ge
    3$ under~\ref{assump:A2}, and when $\dimin = 2$ under isotropic~\ref{assump:A3}.
\end{enumerate}
\end{theorem}

\begin{proof}
Integrating $K_{\mathrm{sph}}(t)$ against Gegenbauer polynomials via the
Funk--Hecke theorem and the Chu--Vandermonde identity
gives~\eqref{eq:funk-hecke-eigenvalues} for $\dimin \ge 3$, while the Chebyshev
limit $\lim_{\alpha \to 0} C_\ell^\alpha(t)/C_\ell^\alpha(1) = T_\ell(t)$ yields
$\kappa_{2m+1}(2) = \lambda_m/2$ at $\dimin = 2$. Related Funk--Hecke expansions
of homogeneous spherical and ReLU kernels appear in~\cite{ChoSaul2009,Bach2017}.
Since $\kappa_{2m+1}(\dimin) = \Theta(m^{-\dimin})$ and $N(\dimin, 2m+1) =
\Theta(m^{\dimin - 2})$, part~(ii) follows by summing over odd degrees. For
part~(iii), Proposition~\ref{prop:degenerate-ustat-limit} identifies the limit
law of $U_{\ngen,\infty}(W)$ with $U_\infty$. Transferring this limit to
$\Regge(\Omega)$ then follows from Theorem~\ref{thm:second-chaos-transfer} when
$\dimin \ge 3$ and from Theorem~\ref{thm:planar-universality} when $\dimin = 2$.
Full details are given in Appendix~\ref{sec:app-spectral}.
\end{proof}

\subsection{Hilbert-space quadratic form and \texorpdfstring{$O(\ngen^{-1})$}{O(H\^{}\{-1\})} Kolmogorov rate}
\label{subsec:hilbert-kolmogorov}

Because $\Zlim(w, w) \equiv 0$, the U-statistic $U_{\ngen,\infty}(W)$ is not
merely close to a quadratic form in the normalised sum of the vectors
$\widetilde{X}(w_h)$ below; it is exactly equal to one, with no diagonal
correction, even though $\norm{\widetilde{X}(w)}_{\ell^2}$ varies radially.
This allows us to apply the Hilbert-space quadratic-form bound of
G\"otze and Zaitsev~\cite[Theorem~1.1]{GoetzeZaitsev2014} directly, rather
than the U-statistic expansions of Bentkus and
G\"otze~\cite{BentkusGoetze1999}.

\begin{theorem}[Berry--Esseen bound for $U_{\ngen,\infty}(W)$]
\label{thm:hilbert-kolmogorov}
Let $\dimin \ge 2$, and let $w_1, \dots, w_{\ngen}$ be i.i.d.\ in $\RR^{\dimin}$
from an isotropic law $P_w$ satisfying~\ref{assump:A3} with $\Sigma_w = \sigw^2
I_{\dimin}$ and $\EE[\norm{w_1}^8] < \infty$.
\begin{enumerate}[label=(\roman*)]
  \item With $\{\psi_k\}_{k \ge 1}$ and $\{\Lambda_k(\dimin)\}_{k \ge 1}$ from
    Theorem~\ref{thm:funk-hecke}\textup{(ii)}, $s_k := \sgn(\Lambda_k(\dimin))
    \in \{\pm 1\}$, $\widetilde{X}(w) := (\sqrt{\abs{\Lambda_k(\dimin)}/2}\,
    \psi_k(w))_{k \ge 1} \in \ell^2$, $J := \operatorname{diag}(s_1, s_2, \dots)$
    and $C_{\mathrm{tr}}(\dimin) := \frac{\dimin^2}{8}\sum_{m=0}^\infty N(\dimin,
    2m+1)\,\abs{\kappa_{2m+1}(\dimin)} < \infty$, one has $\Zlim(w, w) = 0$,
    \begin{equation}\label{eq:hilbert-l2-norm}
      \tr(T_{\Zlim})
      = \int_{\RR^{\dimin}} \Zlim(w, w)\,\dd P_w(w)
      = \sum_{k=1}^\infty \Lambda_k(\dimin)
      = 0,
      \quad
      \norm{\widetilde{X}(w)}_{\ell^2}^2
      = C_{\mathrm{tr}}(\dimin)\, \varrho(w)^2,
    \end{equation}
    and for $\widetilde{S}_{\ngen} := \frac{1}{\sqrt{\ngen}}\sum_{h=1}^{\ngen}
    \widetilde{X}(w_h) \in \ell^2$,
    \begin{equation}\label{eq:hilbert-quad-form}
      U_{\ngen,\infty}(W)
      = \ip{\widetilde{S}_{\ngen}}{J \widetilde{S}_{\ngen}}_{\ell^2}
      = \frac{1}{2}\sum_{k=1}^\infty \Lambda_k(\dimin)\left(\frac{1}{\sqrt{\ngen}}\sum_{h=1}^{\ngen} \psi_k(w_h)\right)^2.
    \end{equation}
  \item With $\gamma_k \stackrel{\mathrm{i.i.d.}}{\sim} \Norm(0, 1)$ and
    $U_\infty := \frac{1}{2}\sum_{k=1}^\infty \Lambda_k(\dimin)(\gamma_k^2 - 1) =
    \frac{1}{2}\sum_{k=1}^\infty \Lambda_k(\dimin)\gamma_k^2$, there exists
    $C_{\dimin} > 0$ depending only on $\dimin$ such that for all $\ngen \ge 2$,
    \begin{equation}\label{eq:hilbert-kolmogorov}
      d_{\mathrm{Kol}}\bigl(U_{\ngen,\infty}(W),\, U_\infty\bigr)
      := \sup_{t \in \RR} \abs[\big]{\PP(U_{\ngen,\infty}(W) \le t) - \PP(U_\infty \le t)}
      \le C_{\dimin}\,
        \frac{\EE[\norm{w_1}^8]}{\bigl(\EE[\norm{w_1}^4]\bigr)^2}\,\ngen^{-1}.
    \end{equation}
    The ratio on the right equals $\EE[\varrho(w_1)^4]/\tau_2(P_w)^2$, so the
    dependence on $P_w$ is only through the dimensionless kurtosis of the
    radial variable $\norm{w_1}^2$.
\end{enumerate}
\end{theorem}

\begin{proof}
\ppart{(i)}
Since $\norm{w \wedge w} = 0$, $\Zlim(w, w) = 0$ and $\tr(T_{\Zlim}) = 0$. By
the spherical harmonic addition theorem~\cite[Theorem~1.2.6]{DaiXu2013},
$\abs{\sum_j Y_{\ell, j}(u)Y_{\ell, j}(v)} \le \sum_j Y_{\ell,j}(u)^2 =
N(\dimin, \ell)$ on $\Sph{\dimin-1}$. Thus $\sum_{m \ge 0} N(\dimin,
2m+1)\abs{\kappa_{2m+1}(\dimin)} < \infty$ makes the Funk--Hecke series converge
absolutely and uniformly to the continuous kernel $K_{\mathrm{sph}}(\ip{u}{v})$
everywhere, giving $\sum_k\Lambda_k(\dimin)\psi_k(w)\psi_k(w') =
2\ip{\widetilde{X}(w)}{J\widetilde{X}(w')}_{\ell^2} = \Zlim(w, w')$ for all $w,
w' \in \RR^{\dimin}\setminus\{0\}$. Setting $u = v$ gives~\eqref{eq:hilbert-l2-norm}
and $\ip{\widetilde{X}(w_h)}{J \widetilde{X}(w_h)}_{\ell^2} =
\frac{\varrho(w_h)^2}{2\tau_2(P_w)}\sum_{k=1}^\infty\Lambda_k(\dimin) = 0$,
proving~\eqref{eq:hilbert-quad-form}. Note that the diagonal terms $h = h'$ of
$\ip{\widetilde{S}_{\ngen}}{J\widetilde{S}_{\ngen}}$ vanish individually, so
that~\eqref{eq:hilbert-quad-form} is an identity and not an approximation.

\ppart{(ii)}
We first remove the dependence on the radial law. Set $\widehat{X}(w) :=
\tau_2(P_w)^{-1/2}\widetilde{X}(w)$ and $\widehat{S}_{\ngen} :=
\frac{1}{\sqrt{\ngen}}\sum_{h=1}^{\ngen} \widehat{X}(w_h)$, so that
$U_{\ngen,\infty}(W) = \tau_2(P_w)\ip{\widehat{S}_{\ngen}}{J\widehat{S}_{\ngen}}$
and $U_\infty \stackrel{d}{=} \tau_2(P_w)\ip{G}{JG}$ for $G \sim \Norm(0,
\Sigma_0)$. By the orthonormality of $\{\psi_k\}$ and
Theorem~\ref{thm:funk-hecke}\textup{(ii)}, the operator
\[
  \Sigma_0 := \Cov\bigl(\widehat{X}(w_1)\bigr)
  = \operatorname{diag}\Bigl(\frac{\abs{\Lambda_k(\dimin)}}{2\tau_2(P_w)}\Bigr)_{k \ge 1}
  = \operatorname{diag}\Bigl(\frac{\dimin^2}{8}\abs{\kappa_{2m+1}(\dimin)}\Bigr)
\]
depends on $\dimin$ alone. Since the Kolmogorov distance is invariant under
multiplying both random variables by the same positive constant, it suffices
to bound $d_{\mathrm{Kol}}(\ip{\widehat{S}_{\ngen}}{J\widehat{S}_{\ngen}},
\ip{G}{JG})$.

We now check the hypotheses of~\cite[Theorem~1.1]{GoetzeZaitsev2014}. That
theorem concerns i.i.d.\ mean-zero vectors $X_h$ in a real separable Hilbert
space $\mathbb{H}$ of dimension $d$ with $5 \le d \le \infty$, with
$\beta_4 := \EE\norm{X_1}^4 < \infty$ and covariance operator $\mathbb{C}$, and
a bounded symmetric operator $\mathbb{Q}$ on $\mathbb{H}$. Its
non-degeneracy hypothesis is that $\mathbb{Q}\mathbb{C}$ has at least five
non-zero eigenvalues, and the constant depends on $(\mathbb{Q}, \mathbb{C})$
only through finitely many of the largest of them, so the theorem applies
verbatim when $d = \infty$ and $\mathbb{C}$ is compact. It gives $\sup_x
\abs{\PP(\ip{S_N}{\mathbb{Q}S_N} \le x) - \PP(\ip{G}{\mathbb{Q}G} \le x)} \le
c(\mathbb{Q}, \mathbb{C})\,\beta_4\,N^{-1}$ for $S_N := N^{-1/2}\sum_{h \le N}
X_h$ and $G \sim \Norm(0, \mathbb{C})$, with no Edgeworth correction term. Here
$\mathbb{H} = \ell^2$ indexed by the odd spherical-harmonic modes $(2m+1, j)$,
so $d = \infty$; $X_h = \widehat{X}(w_h)$ has mean zero by oddness,
$\psi_k(-w) = -\psi_k(w)$; $\beta_4 = \EE\norm{\widehat{X}(w_1)}^4 =
C_{\mathrm{tr}}(\dimin)^2\,\EE[\varrho(w_1)^4]/\tau_2(P_w)^2 < \infty$
by~\eqref{eq:hilbert-l2-norm}. The operator $\mathbb{Q} = J$ is a self-adjoint isometry,
and $J\Sigma_0 = \operatorname{diag}(\Lambda_k(\dimin)/2\tau_2(P_w))_{k \ge 1}$
has no zero eigenvalue at all, because $\kappa_{2m+1}(\dimin)
\ne 0$ for every $m \ge 0$ by Theorem~\ref{thm:funk-hecke}\textup{(i)} (and
$\tau_2(P_w) \ge (\EE[\varrho(w_1)])^2 = 1$ is finite). No finite-rank
truncation is needed, and the question of an $O(N^{-1})$ Edgeworth term does
not arise, both because~\eqref{eq:hilbert-quad-form} is exactly a quadratic
form in $\widehat{S}_{\ngen}$ and because the third-moment tensor
$\EE[\widehat{X}_{k_1}\widehat{X}_{k_2}\widehat{X}_{k_3}]$ vanishes identically
by oddness of $\widehat{X}$. Theorem~1.1 of~\cite{GoetzeZaitsev2014} thus
yields~\eqref{eq:hilbert-kolmogorov} with $C_{\dimin} := c(J,
\Sigma_0)\,C_{\mathrm{tr}}(\dimin)^2$, since $\EE[\varrho(w_1)^4]/\tau_2(P_w)^2
= \EE[\norm{w_1}^8]/(\EE[\norm{w_1}^4])^2$ by $\varrho(w) =
\norm{w}^2/(\dimin\sigw^2)$.
\end{proof}

\subsection{Finite-\texorpdfstring{$\ngen$}{H} parity orthogonality and asymptotic independence}
\label{subsec:parity-orthogonality}

\begin{proposition}
\label{thm:parity-orthogonality}
Let $\dimin \ge 2$, let $(W, b) = (w_h, b_h)_{h=1}^{\ngen} \in (\RR^{\dimin}
\times \RR)^{\ngen}$ satisfy~\ref{assump:A1}, and let $\Omega \subseteq
\RR^{\dimin}$ be Borel with $\beta_\Omega < \infty$, with $\Regge(\Omega)$ given
by~\eqref{eq:model-measures}. Fix a measurable hemisphere
cross-section $\mathcal{S}_+ \subset \RR^{\dimin} \setminus \{0\}$
($\mathcal{S}_+ \cap (-\mathcal{S}_+) = \emptyset$, $\mathcal{S}_+ \cup
(-\mathcal{S}_+) = \RR^{\dimin} \setminus \{0\}$), write $(w_h, b_h) =
\epsilon_h(\mathring{w}_h, \mathring{b}_h)$ with $\mathring{w}_h \in
\mathcal{S}_+$ and $\epsilon_h \in \{\pm 1\}$, and set $\mathcal{F}_{\mathrm{uns}}
:= \sigma(\{(\mathring{w}_h, \mathring{b}_h)\}_{h=1}^{\ngen}) = \sigma(\{w_h
w_h\transp, b_h w_h, b_h^2\}_{h=1}^{\ngen})$. Define $\Zlim$
by~\eqref{eq:limit-mark-def} ($\sigw > 0$) and $U_{\ngen,\infty}(W) :=
\frac{1}{\ngen}\sum_{c < c'}\Zlim(w_c, w_{c'})$.
\begin{enumerate}[label=(\roman*)]
  \item Every unoriented functional $F(W, b)$ (such as the cell count
    $N_{\mathrm{cells}}(\Omega)$, wall volume $V_{\dimin - 1}(\Omega)$ or bone
    volume $V_{\dimin - 2}(\Omega)$) is even under $\mathcal{I}_c\colon (w_c,
    b_c) \mapsto (-w_c, -b_c)$ for all $c \in [\ngen]$ and hence
    $\mathcal{F}_{\mathrm{uns}}$-measurable, whereas $Y_{cc'}(\Omega; W, b)$ and
    $\Zlim(w_c, w_{c'})$ are odd under $\mathcal{I}_c$ and $\mathcal{I}_{c'}$:
    \begin{equation}\label{eq:parity-pathwise-odd}
      \begin{gathered}
        F\bigl(\mathcal{I}_c(W, b)\bigr) = F(W, b) \; (\forall c),
        \qquad
        Y_{cc'}\bigl(\Omega; \mathcal{I}_c(W, b)\bigr) = -Y_{cc'}(\Omega; W, b), \\
        \Zlim(-w_c, w_{c'}) = -\Zlim(w_c, w_{c'}).
      \end{gathered}
    \end{equation}
  \item For every $\ngen \ge 2$ and every $\mathcal{F}_{\mathrm{uns}}$-measurable
    $F \in L^2(\PP)$, $\EE[\Regge(\Omega; W, b) \mid \mathcal{F}_{\mathrm{uns}}]
    = 0$ a.s.\ and
    \begin{equation}\label{eq:parity-finite-ortho}
      \EE\bigl[F(W, b)\, \Regge(\Omega; W, b)\bigr] = 0,
      \qquad
      \Cov\bigl(F(W, b),\, \Regge(\Omega; W, b)\bigr) = 0.
    \end{equation}
    Likewise, $\EE[U_{\ngen,\infty}(W) \mid \mathcal{F}_{\mathrm{uns}}] = 0$ a.s.\
    whenever $\EE[\norm{w_h}^2] < \infty$ ($\forall h$), and $\Cov(F(W, b),
    U_{\ngen,\infty}(W)) = 0$ whenever $\EE[\norm{w_h}^4] < \infty$ ($\forall h$).
    When $\dimin = 2$, $b \equiv 0$ and the rows of $W$ are a.s.\ pairwise
    non-collinear, $\EE[\Regge(\Omega) \mid \mathcal{F}_{\mathrm{uns}}] = 0$
    a.s.\ and~\eqref{eq:parity-finite-ortho} also hold for the homogeneous
    planar apex measure~\eqref{eq:model-apex-convention}.
\end{enumerate}
\end{proposition}

\begin{proof}
Invariance of $\Pi_c$ under $\mathcal{I}_c$ and
Lemma~\ref{lem:sign-flip-oddness} give~\eqref{eq:parity-pathwise-odd}.
Under~\ref{assump:A1}, $(\epsilon_1, \dots, \epsilon_{\ngen})$ are i.i.d.\
Rademacher variables independent of $\mathcal{F}_{\mathrm{uns}}$. For $B = \{c,
c'\}$, set $\mathcal{F}_{\mathrm{uns}}^{(B)} := \mathcal{F}_{\mathrm{uns}} \vee
\sigma(\{\epsilon_h\}_{h \notin B})$. Then $Y_{cc'}(\Omega) =
\epsilon_c\epsilon_{c'}\,\mathring{Y}_{cc'}(\Omega)$ where
$\mathring{Y}_{cc'}(\Omega)$ is $\mathcal{F}_{\mathrm{uns}}^{(B)}$-measurable and
bounded by $\frac{\pi}{2}\beta_\Omega$. This implies $\EE[Y_{cc'}(\Omega) \mid
\mathcal{F}_{\mathrm{uns}}^{(B)}] =
\mathring{Y}_{cc'}(\Omega)\,\EE[\epsilon_c\epsilon_{c'}] = 0$ a.s., giving
$\EE[\Regge(\Omega) \mid \mathcal{F}_{\mathrm{uns}}] = 0$
and~\eqref{eq:parity-finite-ortho}. Under~\eqref{eq:model-apex-convention}
(although $2\pi - \Theta_\infty(W)$ is not pathwise odd under a single flip
$\mathcal{I}_c$), adjoining $\tilde{b}_1, \dots, \tilde{b}_{\ngen}
\stackrel{\mathrm{i.i.d.}}{\sim} \Norm(0, 1)$ independently of $W$ as in
Lemma~\ref{lem:planar-offset-reduction} identifies $2\pi - \Theta_\infty(W) =
\Regge(\RR^2; W, \tilde{b})$ a.s. Here $(w_h, \tilde{b}_h)_{h=1}^{\ngen}$
satisfies~\ref{assump:A1} and $\mathcal{F}_{\mathrm{uns}} \subseteq
\sigma(\{w_h w_h\transp, \tilde{b}_h w_h, \tilde{b}_h^2\}_{h=1}^{\ngen})$, so
the tower property again yields $\EE[\Regge(\Omega) \mid
\mathcal{F}_{\mathrm{uns}}] = 0$ a.s.\ and~\eqref{eq:parity-finite-ortho}.
The claims for $U_{\ngen,\infty}(W)$ follow identically from $\Zlim(w_c, w_{c'})
= \epsilon_c\epsilon_{c'}\Zlim(\mathring{w}_c, \mathring{w}_{c'})$ and
$\abs{\Zlim(w_c, w_{c'})} \le \frac{1}{4\sigw^4}\norm{w_c}^2\norm{w_{c'}}^2$.
\end{proof}

\begin{corollary}
\label{cor:parity-asymp-indep}
In the setting of Proposition~\ref{thm:parity-orthogonality}, let $Y_{\ngen}
\xrightarrow{d} Y_\infty$ be any sequence of $\mathcal{F}_{\mathrm{uns}}$-measurable
real random variables.
\begin{enumerate}[label=(\roman*)]
  \item Assume~\ref{assump:A1} and isotropic~\ref{assump:A3} with $\Sigma_w =
    \sigw^2 I_{\dimin}$. With $U_\infty$ as in
    Theorem~\ref{thm:hilbert-kolmogorov}\textup{(ii)}, one has
    $(U_{\ngen,\infty}(W), Y_{\ngen}) \xrightarrow{d} (U_\infty, Y_\infty)$
    jointly as $\ngen \to \infty$ with $U_\infty \indep Y_\infty$:
    \begin{equation}\label{eq:parity-asymp-indep}
      \lim_{\ngen \to \infty} \EE\left[\exp\bigl(i t_1 U_{\ngen,\infty}(W) + i t_2 Y_{\ngen}\bigr)\right]
      = \EE\bigl[e^{i t_1 U_\infty}\bigr]\, \EE\bigl[e^{i t_2 Y_\infty}\bigr]
      \qquad (\forall t_1, t_2 \in \RR).
    \end{equation}
    When $\dimin = 2$ and $b$ also satisfies the planar general-position
    condition of Lemma~\ref{lem:planar-offset-reduction},
    $\bigl(\frac{\ngen}{4}\Regge(\RR^2), Y_{\ngen}\bigr) \xrightarrow{d}
    \bigl(\frac{\tau_2(P_w)}{2}\mathcal{L}_2, Y_\infty\bigr)$ with
    $\mathcal{L}_2 \indep Y_\infty$.
  \item Let $\Omega \subseteq \RR^{\dimin}$ ($\dimin \ge 2$) be any Borel window
    with $\beta_\Omega < \infty$, and suppose~\ref{assump:A2} holds on $\Omega$
    with symmetric i.i.d.\ offsets. With $U_\infty(\Omega)$ denoting the
    second-chaos limit law~\eqref{eq:general-window-limit} of
    Theorem~\ref{thm:second-chaos-transfer}\textup{(iii)}, both
    $\bigl(U_{\ngen,\infty}(\Omega), Y_{\ngen}\bigr)$ and
    $\bigl(\frac{\ngen}{4}\Regge(\Omega), Y_{\ngen}\bigr)$ converge in
    distribution as $\ngen \to \infty$ to $(U_\infty(\Omega), Y_\infty)$ with
    $U_\infty(\Omega) \indep Y_\infty$ (for $\dimin \ge 3$, $\Omega = B(0, R)$
    and $b \equiv 0$, $U_\infty(B(0, R)) \stackrel{d}{=} \omega_{\dimin - 2}
    R^{\dimin - 2} U_\infty$).
\end{enumerate}
\end{corollary}

\begin{proof}
Consider any i.i.d.\ array $z_1, \dots, z_{\ngen} \sim P_0$ on a Polish space
$\mathcal{Z}_0$ equipped with a $P_0$-preserving involution $z \mapsto -z$, such
that $z_h = \epsilon_h\mathring{z}_h$ with $\mathring{z}_h$
$\mathcal{F}_{\mathrm{uns}}$-measurable and $(\epsilon_1, \dots,
\epsilon_{\ngen})$ i.i.d.\ Rademacher independent of $\mathcal{F}_{\mathrm{uns}}$. We
take $z_h = w_h \in \RR^{\dimin}$ and $K = \Zlim$ in~\textup{(i)}, and we take $z_h
= \xi_h = (w_h, b_h) \in \RR^{\dimin} \times \RR$ and $K = K_\Omega$
in~\textup{(ii)} (where $\ell_{cc'}(\Omega)$ is invariant under $\xi_c \mapsto
-\xi_c$, so $K_\Omega(-\xi, \xi') = -K_\Omega(\xi, \xi')$). In both cases $K \in
L^2(\mathcal{Z}_0^2, P_0^{\otimes 2})$ is symmetric and satisfies $K(-z, z') =
-K(z, z')$, so $(T_K f)(-z) = -(T_K f)(z)$. Every eigenfunction of $T_K$ with
non-zero eigenvalue $\mu_k$ thus satisfies $\psi_k(-z) = -\psi_k(z)$
$P_0$-a.e. Replacing $\psi_k(z)$ by $\frac{1}{2}(\psi_k(z) - \psi_k(-z))$ makes
$\psi_k$ odd everywhere, and adjoining odd orthonormal functions if $T_K$ has
finite rank gives an orthonormal sequence $\{\psi_k\}_{k \ge 1}$ in
$L^2_{\mathrm{odd}}(\mathcal{Z}_0, P_0)$ with $K(z, z') = \sum_{k=1}^\infty \mu_k
\psi_k(z)\psi_k(z')$ in $L^2(P_0^{\otimes 2})$ and $\sum_{k=1}^\infty \mu_k^2 <
\infty$.

Set $U_{\ngen}(K) := \frac{1}{\ngen}\sum_{c < c'} K(z_c, z_{c'})$ and, for a
truncation level $r \ge 1$, $U_{\ngen}^{(r)}(K) := \allowbreak
\frac{1}{\ngen}\sum_{c < c'}\allowbreak\sum_{k=1}^r \allowbreak
\mu_k\psi_k(z_c)\psi_k(z_{c'})$ and $U_\infty^{(r)}(K) := \allowbreak
\frac{1}{2}\sum_{k=1}^r \mu_k(\gamma_k^2 - 1)$ ($\gamma_k
\stackrel{\mathrm{i.i.d.}}{\sim} \Norm(0, 1)$). Since $\EE[\psi_k(z_1)] = 0$,
the summands of $U_{\ngen}(K) - U_{\ngen}^{(r)}(K)$ are pairwise orthogonal
across $\{c, c'\}$, giving $\EE[(U_{\ngen}(K) - U_{\ngen}^{(r)}(K))^2] \le
\frac{1}{2}\sum_{k > r}\mu_k^2 \to 0$ uniformly in $\ngen \ge 2$, and likewise
$\EE[(U_\infty(K) - U_\infty^{(r)}(K))^2] = \frac{1}{2}\sum_{k > r}\mu_k^2 \to
0$.

Fix $r \ge 1$ and $t_1, t_2 \in \RR$, and set $\vec{\psi}_r(z) := (\psi_1(z),
\dots, \psi_r(z))\transp$ and $M_r := \operatorname{diag}(\mu_1, \dots, \mu_r)$.
Because $\psi_k(z_h) = \epsilon_h\psi_k(\mathring{z}_h)$ and $\epsilon_h^2 = 1$,
\[
  U_{\ngen}^{(r)}(K)
  = \frac{1}{2} S_{\ngen,r}\transp M_r S_{\ngen,r}
    - \frac{1}{2\ngen}\sum_{h=1}^{\ngen} \vec{\psi}_r(\mathring{z}_h)\transp M_r \vec{\psi}_r(\mathring{z}_h),
  \qquad
  S_{\ngen,r} := \frac{1}{\sqrt{\ngen}}\sum_{h=1}^{\ngen} \epsilon_h \vec{\psi}_r(\mathring{z}_h).
\]
By oddness, $\psi_j(\mathring{z}_h)\psi_k(\mathring{z}_h) =
\psi_j(z_h)\psi_k(z_h)$, so the $\mathcal{F}_{\mathrm{uns}}$-measurable sample
covariance $A_{\ngen} := \frac{1}{\ngen}\sum_{h=1}^{\ngen}
\vec{\psi}_r(\mathring{z}_h)\vec{\psi}_r(\mathring{z}_h)\transp$ satisfies
$A_{\ngen} \xrightarrow{\PP} I_r$ by the weak law of large numbers. For
each $\delta > 0$, the $\mathcal{F}_{\mathrm{uns}}$-measurable Lindeberg average
$L_{\ngen}(\delta) := \frac{1}{\ngen}\sum_{h=1}^{\ngen}
\norm{\vec{\psi}_r(\mathring{z}_h)}_2^2\allowbreak
\mathbf{1}\{\norm{\vec{\psi}_r(\mathring{z}_h)}_2 > \delta\sqrt{\ngen}\}$ has
mean $\EE[\norm{\vec{\psi}_r(z_1)}_2^2\mathbf{1}\{\norm{\vec{\psi}_r(z_1)}_2 >
\delta\sqrt{\ngen}\}] \to 0$ and hence converges in probability to $0$. (Only
$\psi_k \in L^2(P_0)$ is used here, so the argument does not depend on the
dimension or on any boundedness of the eigenfunctions; in case~\textup{(i)}
the $\psi_k$ are the radially weighted spherical harmonics of
Theorem~\ref{thm:funk-hecke}, square-integrable because $\tau_2(P_w) <
\infty$.) Consider any deterministic sequence of tuples $(v_{j,1}, \dots,
v_{j,\ngen_j}) \in (\RR^r)^{\ngen_j}$ with $\frac{1}{\ngen_j}\sum_h v_{j,h}
v_{j,h}\transp \to I_r$ and $\frac{1}{\ngen_j}\sum_h
\norm{v_{j,h}}_2^2\mathbf{1}\{\norm{v_{j,h}}_2 > \delta\sqrt{\ngen_j}\} \to 0$
($\forall \delta > 0$). The multivariate Lindeberg--Feller theorem for
$\frac{1}{\sqrt{\ngen_j}}\sum_h \epsilon_h v_{j,h}$ and continuous mapping send
$\EE_\epsilon[\exp(i t_1(\frac{1}{2} S\transp M_r S - \frac{1}{2\ngen_j}\sum_h
v_{j,h}\transp M_r v_{j,h}))]$ to the deterministic limit $c_1 := \EE[\exp(i
t_1 U_\infty^{(r)}(K))]$. Setting $\Phi_{\ngen} := \EE[\exp(i t_1
U_{\ngen}^{(r)}(K)) \mid \mathcal{F}_{\mathrm{uns}}]$, it follows from
$(A_{\ngen}, L_{\ngen}(1/m)) \xrightarrow{\PP} (I_r, 0)$ ($\forall m \ge 1$)
and $\abs{\Phi_{\ngen}} \le 1$ that $\Phi_{\ngen} \xrightarrow{\PP} c_1$ and
$\EE\abs{\Phi_{\ngen} - c_1} \to 0$. Since $Y_{\ngen}$ is
$\mathcal{F}_{\mathrm{uns}}$-measurable and $\abs{e^{i t_2 Y_{\ngen}}} = 1$,
\[
  \abs[\Big]{\EE\bigl[e^{i t_1 U_{\ngen}^{(r)}(K) + i t_2 Y_{\ngen}}\bigr]
    - c_1\,\EE\bigl[e^{i t_2 Y_{\ngen}}\bigr]}
  = \abs[\Big]{\EE\bigl[(\Phi_{\ngen} - c_1)\,e^{i t_2 Y_{\ngen}}\bigr]}
  \le \EE\abs{\Phi_{\ngen} - c_1}
  \xrightarrow[\ngen \to \infty]{} 0.
\]
Using $\EE[e^{i t_2 Y_{\ngen}}] \to \EE[e^{i t_2 Y_\infty}]$ and letting $r \to
\infty$ prove~\eqref{eq:parity-asymp-indep} for $U_{\ngen}(K)$. Slutsky's theorem
then transfers joint convergence to $\frac{\ngen}{4}\Regge(\RR^2)$ in~\textup{(i)}
via Lemma~\ref{lem:planar-offset-reduction} and
Theorem~\ref{thm:planar-universality} (with $U_{\ngen,\infty}(\widetilde{W}) =
U_{\ngen,\infty}(W)$ when $\Sigma_w = \sigw^2 I_2$), and to
$\frac{\ngen}{4}\Regge(\Omega)$ in~\textup{(ii)} via
Theorem~\ref{thm:second-chaos-transfer}\textup{(ii)--(iii)} (and
Theorem~\ref{thm:funk-hecke}\textup{(iii)} at $\Omega = B(0, R)$, $b \equiv 0$).
\end{proof}
\section{Concluding remarks and open problems}
\label{sec:remarks}

We have shown that joint central symmetry eliminates both the mean and the first
Hoeffding component of the signed Regge curvature $\Regge(\Omega)$ at every
finite $\ngen$. The leading fluctuations thus sit in the second Hoeffding
subspace, at scale $O(\ngen^{-1})$. In the isotropic planar case, their
limit law has the closed-form moment generating function $-\sec\sqrt{\pi^2 +
\pi s}$. The meaning of this cancellation is the following. Individual bone
curvatures are not small. Each transverse bone carries a deficit of size
$\Theta(\ngen^{-2})$, and there are $\binom{\ngen}{2}$ of them. Thus the unsigned
total variation $\abs{\Regge}(\Omega)$ has $\Theta(1)$ mean on bounded windows.
It is only the cancellation between positive and negative deficits, enforced
pathwise by orientation reversal, that removes the $\Theta(1)$ mean and the
$\Theta(\ngen^{-1/2})$ first-order fluctuations. We expect
$\abs{\Regge}(\Omega)$ itself to obey a central limit theorem at scale
$\ngen^{-1/2}$, driven by fluctuations of the active count $m_B(x)$, though
confirming this requires a separate analysis.

Several natural questions remain open. Firstly, in dimensions $\dimin \ge 3$ one
would like to characterise the spectrum of the slice-weighted operator
$T_{K_\Omega}$ in Theorem~\ref{thm:second-chaos-transfer}\textup{(iii)} beyond
central ball windows ($b_h \equiv 0$, $\Omega = B(0, R)$). One may also ask
for the limit law of spatially damped whole-space pairings $\int_{\RR^{\dimin}}
\varphi(x)\,\Regge(\dd x)$, where $\beta_{\RR^{\dimin}} = \infty$. Secondly, our
transfer bounds treat $\dimin$ as fixed while $\ngen \to \infty$, and it would
be of interest to extend them to non-Gaussian centrally symmetric laws in
$\dimin \ge 3$. The proportional regime $\ngen, \dimin \to \infty$ with
$\ngen/\dimin \to \gamma \in (1, \infty)$ is also open. There
$\mathsf{W}_{TT}(x_0)$ obeys a Marchenko--Pastur law and the typical active
count crosses $\dimin$ near $\gamma = 2$, so one should expect a genuinely
different limit.

Finally, a natural question motivated by deep neural networks is how signed
curvature behaves for multi-layer compositions of piecewise-affine maps (deep
ReLU networks), whose crease loci form bent hyperplane arrangements
(see~\cite{GrigsbyLindsey2022}). Inside each chamber of the preceding layers,
bones formed by two hyperplanes from the same layer still see a conditionally
thinned-Wishart background metric, whereas bones formed across distinct layers
undergo $O(\ngen)$ metric jumps across earlier walls. Whether those cross-layer
deficits remain at scale $O(\ngen^{-1})$ with a vanishing first Hoeffding
component would need to be established by a separate conditional symmetry
argument.

\appendix
\renewcommand{\theHsection}{\Alph{section}}
\section{Thinned-Wishart remainder estimates and third-Hoeffding sharpness}
\label{sec:app-wishart}

Throughout this appendix, fix a bone $B = \{c, c'\}$, the $\sigma$-field
$\mathcal{F}_B$, and an $\mathcal{F}_B$-measurable point $x_0 \in B_{cc'}$
satisfying~\eqref{eq:model-boundary-null}.

In the notation of
Theorem~\ref{thm:conditional-wishart}, write $\mathsf{W}_{NN} = Z_N Z_N\transp$,
$\mathsf{W}_{NT} = Z_N Z_T\transp$ and $\mathsf{W}_{TT} = Z_T Z_T\transp$, where
conditionally on $(\mathcal{F}_B, m_B(x_0) = m)$ the matrices $Z_N \in \RR^{2
\times m}$ and $Z_T \in \RR^{(\dimin - 2) \times m}$ are independent with
i.i.d.\ $\Norm(0, \sigw^2)$ entries. Let $T_B \in \RR^{2 \times 2}$ satisfy
$T_B\transp T_B = L_B\transp L_B$ as in
Proposition~\ref{prop:inv-wishart-moments}, so $(\sigw^{-1}T_B)\transp(\sigw^{-1}T_B)
= \sigw^{-2}L_B\transp L_B$ has diagonal entries $R_c^2, R_{c'}^2$ and
off-diagonal entry $R_c R_{c'}\cos\vartheta_B$.

\begin{proof}[Proof of Theorem~\ref{thm:fisher-limit-mark}\textup{(i)}]
By~\eqref{eq:schur-normal}, $\frac{2}{\ngen \sigw^2} S_N(x_0) = I_2 +
\Delta_{\ngen}$ with
\begin{align*}
  \Delta_{\ngen}
  = \frac{2}{\ngen}\Bigl[
      \sigw^{-2}\bigl(\mathsf{W}_{NN} - m_B\sigw^2 I_2\bigr)
      &+ \Bigl(m_B - \frac{\ngen - 2}{2}\Bigr) I_2 + (\sigw^{-2} - 1) I_2 \\
      &- \sigw^{-2} \mathsf{W}_{NT}(I_{\dimin - 2} + \mathsf{W}_{TT})^{-1} \mathsf{W}_{TN}
    \Bigr].
\end{align*}
By the Marcinkiewicz--Zygmund inequality, the conditional $L^{2p}$ norms of the
first two terms are $O_p(\ngen^{-1/2})$, while the third is deterministic of
order $O_{\sigw}(\ngen^{-1})$. In the fourth term, $P_T' :=
Z_T\transp(I_{\dimin - 2} + Z_T Z_T\transp)^{-1} Z_T$ is dominated in the
Loewner order by the orthogonal projection $\Pi_T'$ onto $\im(P_T')$ (of rank
$r_T \le \dimin - 2$, independent of $Z_N$), so $0 \preceq
\sigw^{-2}\mathsf{W}_{NT}(I_{\dimin - 2} + \mathsf{W}_{TT})^{-1} \mathsf{W}_{TN}
\preceq \sigw^{-2}Z_N \Pi_T' Z_N\transp \sim \Wish_2(r_T, I_2) = O_{L^{2p}}(1)$.
Hence $\EE[\norm{\Delta_{\ngen}}_{\mathrm{op}}^{2p} \mid \mathcal{F}_B]^{1/(2p)}
\le C_{\dimin, p, \sigw}\ngen^{-1/2}$. Writing $\frac{\ngen\sigw^2}{2}
S_N(x_0)^{-1} - I_2 = (I_2 + \Delta_{\ngen})^{-1} - I_2 =
-\frac{\ngen\sigw^2}{2} S_N(x_0)^{-1}\Delta_{\ngen}$ and applying Cauchy--Schwarz
in $L^p(\PP(\cdot \mid \mathcal{F}_B))$ with~\eqref{eq:invw-moment-bound} gives
$\EE[\norm{\frac{\ngen\sigw^2}{2} S_N(x_0)^{-1} - I_2}_{\mathrm{op}}^p \mid
\mathcal{F}_B] \le C'_{\dimin, p, \sigw}\ngen^{-p/2}$; multiplying by the
$\mathcal{F}_B$-measurable factor $\norm{T_B}_{\mathrm{op}}^{2p} =
\norm{L_B}_{\mathrm{op}}^{2p} \in L^1(\PP)$ and taking expectations
yields~\eqref{eq:compressed-conv}.
\end{proof}

Fix now $m_B(x_0) = m \ge \dimin$ with $\nu := m - \dimin + 2 \ge 2$, and write
$\PP_m$, $\EE_m$ for probabilities and expectations conditional on
$(\mathcal{F}_B, m_B(x_0) = m)$. Set $\widetilde{Z}_N := \sigw^{-1}Z_N$,
$\widetilde{Z}_T := \sigw^{-1}Z_T$, $\widetilde{\mathsf{W}}_{TT} :=
\widetilde{Z}_T\widetilde{Z}_T\transp \sim \Wish_{\dimin - 2}(m, I_{\dimin -
2})$, $\widetilde{S}_N := \sigw^{-2}S_N(x_0)$ and $\widetilde{S}_N^{\mathrm{unreg}}
:= \sigw^{-2}S_N^{\mathrm{unreg}}(x_0) \sim \Wish_2(\nu, I_2)$, so $M_B(x_0) =
\sigw^{-2}T_B\transp \widetilde{S}_N^{-1} T_B$ and $M_B^{\mathrm{unreg}}(x_0) =
\sigw^{-2}T_B\transp (\widetilde{S}_N^{\mathrm{unreg}})^{-1} T_B$. We use the
lower-triangular Bartlett factorisation $\widetilde{S}_N^{\mathrm{unreg}} = A
A\transp$, $A = \bigl(\begin{smallmatrix} V_1 & 0 \\ Z & V_2
\end{smallmatrix}\bigr)$, of Proposition~\ref{prop:inv-wishart-moments}\textup{(ii)},
with $(V_1^2, V_2^2, Z) \sim \chi_\nu^2 \otimes \chi_{\nu - 1}^2 \otimes
\Norm(0, 1)$. In these variables, $\EE[(\chi_k^2)^{-r/2}] = O(k^{-r/2})$ ($k > r$) and
$\EE[(\chi_k^2)^{r/2}] = O(k^{r/2})$ ($r \ge 0$) imply that any finite
combination of monomials $V_1^{-a}V_2^{-b}\abs{Z}^c$ ($c \ge 0$, $a, b \le m_0$)
homogeneous of degree $-2k$ in $(V_1, V_2, Z)$ ($a + b - c = 2k$) is
$O_{L^p}(\nu^{-k})$ for $\nu \ge \nu_0(p) := \lceil p m_0\rceil + 2$.

\begin{lemma}
\label{lem:app-bartlett}
Under Proposition~\ref{prop:inv-wishart-moments}, let $m_B(x_0) = m \ge \dimin$
($\nu = m - \dimin + 2$).
\begin{enumerate}[label=(\roman*)]
  \item In the Bartlett frame of Proposition~\ref{prop:inv-wishart-moments}\textup{(ii)}
    aligned with $Q_N\transp w_c/\norm{w_c}$,
    \[
      \varepsilon_0(M_B^{\mathrm{unreg}})
      = -\frac{R_c^2 R_{c'}^2}{4}\Bigl(
          \cos\vartheta_B\sin\vartheta_B\bigl(V_1^{-3}V_2^{-1} + Z^2 V_1^{-3}V_2^{-3}\bigr)
          - \sin^2\vartheta_B\, Z V_1^{-2}V_2^{-3}
        \Bigr),
    \]
    where the first term equals $\Zlim(w_c, w_{c'})(V_1^{-3}V_2^{-1} + Z^2
    V_1^{-3}V_2^{-3})$.
  \item For every $p \ge 1$ and $m \ge \dimin + \lceil 8p \rceil$, with
    $C_{\dimin, p} < \infty$ independent of $\sigw > 0$,
    \begin{equation}\label{eq:app-bartlett-cond-lp}
      \begin{aligned}
        \EE_m\bigl[\abs{\nu^2 \varepsilon_B(x_0)}^p\bigr]
        &\le \EE_m\bigl[\abs{\nu^2 \varepsilon_0(x_0)}^p\bigr]
        \le C_{\dimin, p}\, R_c^{2p} R_{c'}^{2p}, \\
        \EE_m\bigl[\abs{\nu^2 (\varepsilon_B(x_0) - \varepsilon_0(x_0))}^p\bigr]
        &\le C_{\dimin, p}\, R_c^{2p} R_{c'}^{2p}(1 + R_c^2 + R_{c'}^2)^p\, \nu^{-p}.
      \end{aligned}
    \end{equation}
\end{enumerate}
\end{lemma}

\begin{proof}
For~\textup{(i)}, in the frame where $\sigw^{-1}T_B = \bigl(\begin{smallmatrix}
R_c & R_{c'}\cos\vartheta_B \\ 0 & R_{c'}\sin\vartheta_B
\end{smallmatrix}\bigr)$, inverting the lower-triangular factor gives
$(\widetilde{S}_N^{\mathrm{unreg}})^{-1} = A^{-\top}A^{-1}$ with entries
$V_1^{-2} + Z^2V_1^{-2}V_2^{-2}$, $-ZV_1^{-1}V_2^{-2}$ and $V_2^{-2}$, and
substituting into $\varepsilon_0(M) = -\frac{1}{4}m_{12}\sqrt{\det M}$ with
$\det M_B^{\mathrm{unreg}} = R_c^2R_{c'}^2\sin^2\vartheta_B\,V_1^{-2}V_2^{-2}$
yields the stated formula. For~\textup{(ii)},
Proposition~\ref{prop:deficit-integral}~\eqref{eq:deficit-sandwich}--\eqref{eq:deficit-linear-rem}
gives $\abs{\varepsilon_B(x_0)} \le \abs{\varepsilon_0(x_0)} \le \frac{1}{4}
a_1(x_0) a_2(x_0)$ and $\abs{\varepsilon_B - \varepsilon_0} \le (a_1 +
a_2)\abs{\varepsilon_0} \le \frac{1}{4} a_1 a_2(a_1 + a_2)$. By~\eqref{eq:schur-sandwich},
$M_B(x_0) \preceq M_B^{\mathrm{unreg}}(x_0)$ and $a_i \le
a_i^{\mathrm{unreg}} = R_{c_i}^2 / V_{m - \dimin + 1}$, so
Proposition~\ref{prop:inv-wishart-moments}\textup{(iii)} and H\"older's
inequality yield~\eqref{eq:app-bartlett-cond-lp}.
\end{proof}

\begin{lemma}
\label{lem:app-schur-shift}
Let $\Pi_T := \widetilde{Z}_T\transp\widetilde{\mathsf{W}}_{TT}^{-1}\widetilde{Z}_T
\in \RR^{m \times m}$ be the orthogonal projection onto the row space of
$\widetilde{Z}_T$ (of rank $\dimin - 2$ a.s.) and set $P_T :=
\widetilde{Z}_T\transp\widetilde{\mathsf{W}}_{TT}^{-1}(\sigw^{-2}I_{\dimin - 2}
+ \widetilde{\mathsf{W}}_{TT})^{-1}\widetilde{Z}_T$.
\begin{enumerate}[label=(\roman*)]
  \item One has $\widetilde{S}_N^{\mathrm{unreg}} = \widetilde{Z}_N(I_m -
    \Pi_T)\widetilde{Z}_N\transp$ and $\widetilde{S}_N =
    \widetilde{S}_N^{\mathrm{unreg}} + \sigw^{-2}(I_2 + N_T)$, where $N_T :=
    \widetilde{Z}_N P_T \widetilde{Z}_N\transp \succeq 0$ ($N_T \equiv 0$ at
    $\dimin = 2$) is independent of $\widetilde{S}_N^{\mathrm{unreg}} \sim
    \Wish_2(\nu, I_2)$ under $\PP_m$, with $\EE_m[\norm{N_T}_{\mathrm{op}}^q]^{1/q}
    \le C_{\dimin, q}\nu^{-1}$ for $\nu \ge 2q$.
  \item For every $p \ge 1$ and $\nu \ge \nu_0(p)$,
    \begin{equation}\label{eq:app-reg-Lp}
      \norm[\big]{\nu^2\varepsilon_0(M_B(x_0)) - \nu^2\varepsilon_0(M_B^{\mathrm{unreg}}(x_0))}_{L^p(\PP_m)}
      \le C_{\dimin, p, \sigw}\,R_c^2R_{c'}^2\,\nu^{-1}.
    \end{equation}
\end{enumerate}
\end{lemma}

\begin{proof}
\ppart{(i)}
Dividing $S_N - S_N^{\mathrm{unreg}} = I_2 +
\mathsf{W}_{NT}\mathsf{W}_{TT}^{-1}(I_{\dimin - 2} +
\mathsf{W}_{TT})^{-1}\mathsf{W}_{TN}$ by $\sigw^2$ gives $\widetilde{S}_N =
\widetilde{S}_N^{\mathrm{unreg}} + \sigw^{-2}(I_2 + N_T)$. Since $P_T = \Pi_T P_T
\Pi_T$, conditioning on $\widetilde{Z}_T$ splits $\widetilde{Z}_N$ into
independent Gaussian projections $\widetilde{Z}_N(I_m - \Pi_T)$ and
$\widetilde{Z}_N\Pi_T$, so $\widetilde{S}_N^{\mathrm{unreg}} \sim \Wish_2(\nu,
I_2)$ is independent of $(\widetilde{Z}_T, N_T)$ and $\norm{\widetilde{Z}_N\Pi_T}_{\mathrm{F}}^2
\sim \chi_{2(\dimin - 2)}^2$ is independent of $\widetilde{Z}_T$. Since $\widetilde{\mathsf{W}}_{TT}^{-1}$ and $(\sigw^{-2}I_{\dimin - 2} +
\widetilde{\mathsf{W}}_{TT})^{-1} \preceq \widetilde{\mathsf{W}}_{TT}^{-1}$
commute, $0 \preceq P_T \preceq
\widetilde{Z}_T\transp\widetilde{\mathsf{W}}_{TT}^{-2}\widetilde{Z}_T$, so we
have $\norm{N_T}_{\mathrm{op}} \le \norm{\widetilde{Z}_N\Pi_T}_{\mathrm{F}}^2
\tr(\widetilde{\mathsf{W}}_{TT}^{-1})$, and each diagonal entry of
$\widetilde{\mathsf{W}}_{TT}^{-1}$ has law $1/\chi_{\nu + 1}^2$, yielding the
$L^q$ bound.

\ppart{(ii)}
With $J_\perp := \bigl(\begin{smallmatrix} 0 & -1 \\ 1 & 0 \end{smallmatrix}\bigr)$,
the $2 \times 2$ adjugate is $\operatorname{adj}(S) = J_\perp\transp S J_\perp$,
so for $x := \sigw^{-1}J_\perp\transp T_B e_1$ and $y := \sigw^{-1}J_\perp\transp
T_B e_2$ (which satisfy $\norm{x} = R_c$ and $\norm{y} = R_{c'}$) the
off-diagonal entry of $\sigw^{-2}T_B\transp S^{-1}T_B$ is $x\transp S y/\det
S$, while its determinant is $R_c^2R_{c'}^2\sin^2\vartheta_B/\det S$. Hence
\[
  \varepsilon_0(\sigw^{-2}T_B\transp S^{-1} T_B)
  = -\tfrac{1}{4}R_c R_{c'}\sin\vartheta_B\,f(S),
  \qquad
  f(S) := x\transp S y\,(\det S)^{-3/2},
  \quad S \in \mathrm{Sym}_2^{++}.
\]
Now let $\Delta \succeq 0$ and $t := \tr(S^{-1}\Delta)$. Since the eigenvalues
of $S^{-1/2}\Delta S^{-1/2}$ are non-negative with sum $t$, we have $\det(S +
\Delta) = \det S\,(1 + u)$ with $t \le u \le t + t^2/4$, and therefore $0 \le 1
- (1 + u)^{-3/2} \le \min\{1, \tfrac{3}{2}u\} \le 3t$. Writing
\[
  f(S + \Delta) - f(S)
  = \frac{x\transp \Delta y}{\det(S + \Delta)^{3/2}}
    - \frac{x\transp S y}{(\det S)^{3/2}}\bigl(1 - (1 + u)^{-3/2}\bigr)
\]
and using $\det(S + \Delta) \ge \det S$ gives $\abs{f(S + \Delta) - f(S)} \le
(\det S)^{-3/2}\bigl[R_c R_{c'}\norm{\Delta}_{\mathrm{op}} + 3\abs{x\transp S
y}\,t\bigr]$. Evaluating at $S = \widetilde{S}_N^{\mathrm{unreg}}$ and $\Delta =
\sigw^{-2}(I_2 + N_T)$, so that $t \le \sigw^{-2}(1 +
\norm{N_T}_{\mathrm{op}})\tr((\widetilde{S}_N^{\mathrm{unreg}})^{-1})$, bounds
$\abs{\varepsilon_0(M_B(x_0)) - \varepsilon_0(M_B^{\mathrm{unreg}}(x_0))}$ by
$\sigw^{-2}(1 + \norm{N_T}_{\mathrm{op}})R_c^2 R_{c'}^2$ times a combination of
Bartlett monomials of degree $-6$, whose $L^p(\PP_m)$ norm is $O(\nu^{-3})$ by
the remark preceding Lemma~\ref{lem:app-bartlett}. Multiplying by $\nu^2$ and
using the $L^q$ bound on $\norm{N_T}_{\mathrm{op}}$ from part~\textup{(i)} with
H\"older's inequality proves~\eqref{eq:app-reg-Lp}.
\end{proof}

\begin{proof}[Proof of the sharpness assertion in Theorem~\ref{thm:second-chaos-transfer}\textup{(ii)}]
The upper bound $O(\ngen^{-1/2})$ holds by~\eqref{eq:transfer-pythagorean} and
Lemma~\ref{lem:cross-bone-es}~\eqref{eq:es-variance-bound}. Assume now that
$\PP(\ell_{12}(\Omega) > 0) > 0$. We reduce
the third Hoeffding component to the single bone $B_{12}$ by symmetrising
under the orientation flip of the third generator. This flip kills the other two
bones through $\{1, 2, 3\}$. The surviving contribution is the even part, in $w_3$, of the
response of $\varepsilon_{12}$ to the rank-one perturbation $w_3 w_3\transp$ of
the background metric. We use a second-order expansion to isolate its leading term.
This term is linear in $w_3 w_3\transp - \sigw^2 I_{\dimin}$ and is of order
$\ngen^{-3}$. A direct calculation then shows that this term has non-zero
conditional variance.

\ppart{Reduction to a single bone.}
By exchangeability across the $\binom{\ngen}{3}$ triples and
Lemma~\ref{lem:hoeffding-support},
\[
  \Var\bigl(P_3(\Regge(\Omega))\bigr)
  = \binom{\ngen}{3}\norm{\pi_{\{1,2,3\}}\Regge(\Omega)}_{L^2(\PP)}^2,
\]
with $\pi_{\{1,2,3\}}\Regge(\Omega) = \pi_{\{1,2,3\}}Y_{\{1,2\}}(\Omega) +
\pi_{\{1,2,3\}}Y_{\{1,3\}}(\Omega) + \pi_{\{1,2,3\}}Y_{\{2,3\}}(\Omega)$. Under
the third-coordinate sign flip $\mathcal{I}_3\colon (w_3, b_3) \mapsto (-w_3,
-b_3)$, Lemma~\ref{lem:sign-flip-oddness} shows that $Y_{\{1,3\}}(\Omega)$ and
$Y_{\{2,3\}}(\Omega)$ are odd, and because $\pi_{\{1,2,3\}}$ commutes with
$\mathcal{I}_3$ their projections $\pi_{\{1,2,3\}}Y_{\{1,3\}}(\Omega)$ and
$\pi_{\{1,2,3\}}Y_{\{2,3\}}(\Omega)$ are also odd under $\mathcal{I}_3$.
Applying the $L^2(\PP)$-orthogonal projection $\frac{1}{2}(I + \mathcal{I}_3)$
annihilates both terms and gives
\[
  \norm{\pi_{\{1,2,3\}}\Regge(\Omega)}_{L^2(\PP)} \ge \norm{T}_{L^2(\PP)},
  \qquad
  T := \frac{1}{2}(I + \mathcal{I}_3)\pi_{\{1,2,3\}}Y_{\{1,2\}}(\Omega).
\]

\ppart{Isolating the rank-one perturbation.}
By~\eqref{eq:hoeffding-P2-collapse}, $\pi_{\{1,2,3\}}Y_{\{1,2\}}(\Omega) =
\EE[Y_{\{1,2\}}(\Omega) \mid \xi_1, \xi_2, \xi_3] - \EE[Y_{\{1,2\}}(\Omega) \mid
\xi_1, \xi_2]$. Consider $x \in B_{12}^\circ \setminus (\Pi_3 \cup \{0\})$. The excluded set is $\Haus{\dimin-2}$-null a.s.\ because, by~\eqref{eq:model-boundary-null},
$B_{12} \cap \Pi_3$ is a.s.\ a proper affine subspace of $B_{12}$. Meanwhile, the
point $0$ lies on $B_{12}$ only when $b_1 = b_2 = 0$, which has probability
zero at $\dimin = 2$ and for $\dimin \ge 3$ removes a single point from the
$(\dimin-2)$-flat $B_{12}$. At such $x$, the map
$\mathcal{I}_3$ exchanges $I_3(x)$ with $1 - I_3(x)$ while leaving the
remaining background metric $G_{-3}(x) := I_{\dimin} + \sum_{h=4}^{\ngen} I_h(x)
w_h w_h\transp$ unchanged. Thus the symmetrised pointwise deficit
\[
  \varepsilon^{\mathrm{sym}}_{12}(x)
  := \frac{1}{2}(I + \mathcal{I}_3)\varepsilon_{12}(x)
   = \frac{1}{2}\Bigl[\varepsilon_{12}\bigl(G_{-3}(x)\bigr)
     + \varepsilon_{12}\bigl(G_{-3}(x) + w_3 w_3\transp\bigr)\Bigr]
\]
depends on $\xi_3$ only through $w_3 w_3\transp$. Defining the linear operator
\[
  \mathcal{P}_3(Z)
  := \EE\bigl[Z \;\big|\; \xi_1, \xi_2, w_3 w_3\transp\bigr]
     - \EE\bigl[Z \;\big|\; \xi_1, \xi_2\bigr],
\]
we have $T = \mathcal{P}_3\bigl(\int_{B_{12} \cap \Omega}
\varepsilon^{\mathrm{sym}}_{12}(x)\,\dd\Haus{\dimin-2}(x)\bigr)$. Because the
two conditioning $\sigma$-fields are nested, $\norm{\mathcal{P}_3(Z)}_{L^2(\PP)}
\le \norm{Z}_{L^2(\PP)}$ for every $Z \in L^2(\PP)$, and $\mathcal{P}_3(Z) = 0$
whenever $Z$ is conditionally independent of $w_3$ given $(\xi_1, \xi_2)$.

\ppart{Second-order expansion and remainders.}
Set $\bar{g} := 1 + \frac{\ngen - 2}{2}\sigw^2$ and choose a small universal
constant $c \in (0, \frac{1}{12}]$. Define $\mathcal{G}_{12} := \{\norm{w_1}^2 +
\norm{w_2}^2 \le c\bar{g}\} \in \mathcal{F}_{\{1,2\}}$ and write
\[
  G_{-3}(x) = \bar{g}\bigl(I_{\dimin} + \Upsilon(x) + \Upsilon_-\bigr),
  \qquad
  G_{-3}(x) + w_3 w_3\transp = \bar{g}\bigl(I_{\dimin} + \Upsilon(x) + \Upsilon_+\bigr),
\]
where $\Upsilon(x) := \bar{g}^{-1}\bigl(\sum_{h=4}^{\ngen} I_h(x) w_h w_h\transp
- \frac{\ngen - 3}{2}\sigw^2 I_{\dimin}\bigr)$, $\Upsilon_- :=
-\frac{\sigw^2}{2\bar{g}}I_{\dimin}$ and $\Upsilon_+ := \frac{1}{2\bar{g}}(2w_3
w_3\transp - \sigw^2 I_{\dimin})$, with average $\Upsilon_3 :=
\frac{1}{2}(\Upsilon_- + \Upsilon_+) = \frac{1}{2\bar{g}}(w_3 w_3\transp -
\sigw^2 I_{\dimin})$. Conditionally on $(\xi_1, \xi_2)$, $\Upsilon(x)$ is
centred and independent of $w_3$: indeed $\EE[\Upsilon(x) \mid \xi_1, \xi_2] =
0$ because~\eqref{eq:model-boundary-null} holds at $x \ne 0$, so that
Lemma~\ref{lem:active-count} makes each $I_h(x)$ ($h \ge 4$) a fair coin flip
independent of $w_h w_h\transp$. For $\ngen$ large enough that
$\norm{\Upsilon_-}_{\mathrm{op}} \le c$, set
\[
  \mathcal{E}(x)
  := \mathcal{G}_{12}
     \cap \bigl\{\norm{\Upsilon(x)}_{\mathrm{op}} \le c\bigr\}
     \cap \bigl\{\norm{\Upsilon_+}_{\mathrm{op}} \le c\bigr\}.
\]
Writing $\mathcal{L} := \mathcal{L}_{\bar{g}I_{\dimin}, L_{12}}$ and
$\mathcal{Q} := \mathcal{Q}_{\bar{g}I_{\dimin}, L_{12}}$,
Lemma~\ref{lem:deficit-perturbation}\textup{(i)} at $G_0 = \bar{g}I_{\dimin}$
expands $\varepsilon^{\mathrm{sym}}_{12}(x)$ on $\mathcal{E}(x)$ as
$\varepsilon^{\mathrm{sym}}_{12}(x) = P_{\mathrm{poly}}(x) +
R_{\mathrm{cub}}(x)$, where
\begin{align*}
  P_{\mathrm{poly}}(x)
  &:= \varepsilon_{12}(\bar{g}I_{\dimin}) + \mathcal{L}[\Upsilon(x)]
     + \mathcal{L}[\Upsilon_3] + \tfrac{1}{2}\mathcal{Q}[\Upsilon(x), \Upsilon(x)]
     + \mathcal{Q}[\Upsilon(x), \Upsilon_3] \\
  &\qquad + \tfrac{1}{4}\bigl(\mathcal{Q}[\Upsilon_+, \Upsilon_+]
       + \mathcal{Q}[\Upsilon_-, \Upsilon_-]\bigr)
\end{align*}
and
\[
  \abs{R_{\mathrm{cub}}(x)}\mathbf{1}_{\mathcal{E}(x)}
  \le C\bar{g}^{-2}\norm{w_1}^2\norm{w_2}^2
    \bigl(\norm{\Upsilon(x)}_{\mathrm{op}} + \norm{\Upsilon_+}_{\mathrm{op}}
    + \norm{\Upsilon_-}_{\mathrm{op}}\bigr)^3.
\]
Split
\[
  \varepsilon^{\mathrm{sym}}_{12}(x)
  = P_{\mathrm{poly}}(x)\mathbf{1}_{\mathcal{G}_{12}}
    + R_{\mathrm{cub}}(x)\mathbf{1}_{\mathcal{E}(x)}
    + \bigl(\varepsilon^{\mathrm{sym}}_{12}(x) - P_{\mathrm{poly}}(x)\bigr)\mathbf{1}_{\mathcal{G}_{12} \setminus \mathcal{E}(x)}
    + \varepsilon^{\mathrm{sym}}_{12}(x)\mathbf{1}_{\mathcal{G}_{12}^c}.
\]
We integrate each term over $B_{12} \cap \Omega$ and apply $\mathcal{P}_3$.
Inside $P_{\mathrm{poly}}(x)\mathbf{1}_{\mathcal{G}_{12}}$, the operator
$\mathcal{P}_3$ annihilates all $w_3$-free summands and the bilinear term
$\mathcal{Q}[\Upsilon(x), \Upsilon_3]\mathbf{1}_{\mathcal{G}_{12}}$ (since
$\EE[\Upsilon(x) \mid \xi_1, \xi_2, w_3 w_3\transp] = \EE[\Upsilon(x) \mid
\xi_1, \xi_2] = 0$), while preserving $\mathcal{L}[\Upsilon_3]\mathbf{1}_{\mathcal{G}_{12}}$
unchanged. It also sends $\frac{1}{4}\mathcal{Q}[\Upsilon_+,
\Upsilon_+]\mathbf{1}_{\mathcal{G}_{12}}$ to an error bounded pointwise by
$C\bar{g}^{-2}\norm{w_1}^2\norm{w_2}^2\norm{\Upsilon_+}_{\mathrm{op}}^2$, which
is $O_{L^2(\PP)}(\ngen^{-4})$ because $\bar{g}^{-2} = O(\ngen^{-2})$ and
$\norm{\Upsilon_+}_{\mathrm{op}} = O_{L^p}(\ngen^{-1})$. For the remaining three
terms, under~\ref{assump:A2} the conditional thinned-Wishart concentration of
$\Upsilon(x)$ (Theorem~\ref{thm:conditional-wishart}
and~\cite{DavidsonSzarek2001}) gives
\[
  \EE\bigl[\norm{\Upsilon(x)}_{\mathrm{op}}^6 \;\big|\; \xi_1, \xi_2\bigr]
  \le C_{\dimin}\ngen^{-3},
  \qquad
  \PP\bigl(\norm{\Upsilon(x)}_{\mathrm{op}} > c \;\big|\; \xi_1, \xi_2\bigr)
  \le C_{\dimin} e^{-c_1\ngen},
\]
while Gaussian chi-squared tails give $\PP(\mathcal{G}_{12}^c)
+ \PP(\norm{\Upsilon_+}_{\mathrm{op}} > c) \le C_{\dimin} e^{-c_1\ngen}$. Thus,
on $\mathcal{G}_{12}$, the cubic term has conditional $L^2$ norm at most
\[
  C\bar{g}^{-2}\norm{w_1}^2\norm{w_2}^2 \cdot O(\ngen^{-3/2}) = O(\ngen^{-7/2}).
\]
Since $\abs{\varepsilon^{\mathrm{sym}}_{12}(x)} \le
\pi/2$ deterministically, $\varepsilon^{\mathrm{sym}}_{12}(x)$ contributes
$O(e^{-c_1\ngen/2})$ in $L^2(\PP)$ on both $\mathcal{G}_{12}^c$ and
$\mathcal{G}_{12} \setminus \mathcal{E}(x)$; meanwhile
\[
  \abs{P_{\mathrm{poly}}(x)}\mathbf{1}_{\mathcal{G}_{12}}
  \le C\bar{g}^{-2}\norm{w_1}^2\norm{w_2}^2\bigl(1 + \norm{\Upsilon(x)}_{\mathrm{op}}^2
    + \norm{\Upsilon_+}_{\mathrm{op}}^2
    + \norm{\Upsilon_-}_{\mathrm{op}}^2\bigr)\mathbf{1}_{\mathcal{G}_{12}}
\]
has conditional $L^4$ norm $O(\ngen^{-2})$, so Cauchy--Schwarz gives
\[
  \norm{P_{\mathrm{poly}}(x)\mathbf{1}_{\mathcal{G}_{12} \setminus
  \mathcal{E}(x)}}_{L^2(\PP)} \le C_{\dimin}\ngen^{-2}e^{-c_1\ngen/4}.
\]
By
Minkowski's integral inequality on $\ell_{12}(\Omega) \le \beta_\Omega$ and the
contractivity of $\mathcal{P}_3$, we are left with
\[
  T = \mathcal{L}_{\bar{g}I_{\dimin}, L_{12}}[\Upsilon_3]\,\ell_{12}(\Omega)\,\mathbf{1}_{\mathcal{G}_{12}}
      + O_{L^2(\PP)}\bigl(\beta_\Omega\,\ngen^{-7/2}\bigr).
\]

\ppart{Non-degeneracy of the leading term.}
Set $\Gamma := L_{12}\transp L_{12}$, $M_0 := \bar{g}^{-1}\Gamma$ and $\gamma :=
\sigw^{-1}L_{12}\transp w_3$, so that $\gamma \mid \mathcal{F}_{\{1,2\}} \sim
\Norm(0, \Gamma)$. Write $\nabla f(M) \in \mathrm{Sym}_2(\RR)$ for the
symmetric Frobenius gradient ($Df(M)[\Xi] = \ip{\nabla f(M)}{\Xi}_{\mathrm{F}}$
for $\Xi \in \mathrm{Sym}_2(\RR)$, meaning $(\nabla f)_{11} = \partial_{m_{11}} f$,
$(\nabla f)_{22} = \partial_{m_{22}} f$ and $(\nabla f)_{12} =
\frac{1}{2}\partial_{m_{12}} f$). The chain rule then gives
$\mathcal{L}_{\bar{g}I_{\dimin}, L_{12}}[\Upsilon_3] = \allowbreak
-\frac{\sigw^2}{2\bar{g}^2}\allowbreak\ip{\nabla\Psi_{\mathrm{mat}}(M_0)}{\gamma\gamma\transp
- \Gamma}_{\mathrm{F}}$. We use $\Psi_{\mathrm{mat}}(M) =
\varepsilon_0(M)\mathcal{J}(M)$, the degree-$2$ homogeneity
$\nabla\varepsilon_0(M_0) = \bar{g}^{-1}\nabla\varepsilon_0(\Gamma)$, the estimate $0 \le 1 -
\mathcal{J}(M_0) \le \tr M_0$
(Proposition~\ref{prop:deficit-integral}\textup{(iii)}), and the bound
$\norm{\nabla\mathcal{J}(M_0)}_{\mathrm{F}} \le C$
(Lemma~\ref{lem:deficit-smoothness}\textup{(i)}). Together, these yield
\[
  \norm[\big]{\Gamma^{1/2}\nabla\Psi_{\mathrm{mat}}(M_0)\Gamma^{1/2}
    - \bar{g}^{-1}\Gamma^{1/2}\nabla\varepsilon_0(\Gamma)\Gamma^{1/2}}_{\mathrm{F}}
  \le C\bar{g}^{-2}\norm{w_1}^2\norm{w_2}^2\bigl(\norm{w_1}^2 + \norm{w_2}^2\bigr)
\]
on $\mathcal{G}_{12}$, and dominated convergence yields
\[
  \bar{g}^3 T
  \;\xrightarrow[\ngen \to \infty]{L^2(\PP)}\;
  -\frac{\sigw^2}{2}\,\ell_{12}(\Omega)\,\bigl\langle\nabla\varepsilon_0(\Gamma),\, \gamma\gamma\transp - \Gamma\bigr\rangle_{\mathrm{F}}.
\]
Conditionally on $\mathcal{F}_{\{1,2\}}$,
$\ip{\nabla\varepsilon_0(\Gamma)}{\gamma\gamma\transp - \Gamma}_{\mathrm{F}}$
has variance $2\norm{\Gamma^{1/2}\nabla\varepsilon_0(\Gamma)\Gamma^{1/2}}_{\mathrm{F}}^2$,
which is strictly positive whenever $w_1 \wedge w_2 \ne 0$: indeed,
$(\nabla\varepsilon_0(\Gamma))_{11} =
-\frac{\Gamma_{12}\Gamma_{22}}{8\sqrt{\det\Gamma}} \ne 0$ when $\Gamma_{12} \ne
0$, whereas $(\nabla\varepsilon_0(\Gamma))_{12} =
\frac{1}{2}\partial_{\Gamma_{12}}\varepsilon_0(\Gamma)\big|_{\Gamma_{12}=0} =
-\frac{1}{8}\sqrt{\Gamma_{11}\Gamma_{22}} \ne 0$ when $\Gamma_{12} = 0$. Since
$\PP(\ell_{12}(\Omega) > 0) > 0$, we conclude that $\lim_{\ngen \to \infty}
\bar{g}^6 \EE[T^2] > 0$, whence $\Var(P_3(\Regge(\Omega))) \ge
\binom{\ngen}{3}\EE[T^2] = \Theta(\ngen^{-3})$ and
$\norm{\frac{\ngen}{4}P_3(\Regge(\Omega))}_{L^2(\PP)} = \Theta(\ngen^{-1/2})$.
\end{proof}

\section{Funk--Hecke spectral proof and degenerate U-statistic limit}
\label{sec:app-spectral}

\begin{proof}[Proof of Theorem~\ref{thm:funk-hecke}]
\ppart{(i)}
By the Funk--Hecke theorem~\cite{DaiXu2013} with $\alpha = (\dimin - 2)/2$ and
$\sigma_{k-1} = 2\pi^{k/2}/\Gamma(k/2)$,
\[
  \kappa_\ell(\dimin)
  = \frac{\sigma_{\dimin - 2}}{\sigma_{\dimin - 1}}
    \int_{-1}^1 K_{\mathrm{sph}}(t)\frac{C_\ell^\alpha(t)}{C_\ell^\alpha(1)}(1 - t^2)^{\alpha - 1/2}\,\dd t,
\]
which vanishes for $\ell = 2m$ by parity. For $\ell = 2m + 1$ and $\dimin \ge 3$
($\alpha > 0$), reversing the descending power series~\cite[eq.~(18.5.10)]{DLMF}
writes the odd Gegenbauer polynomial in hypergeometric form as
\[
  C_{2m+1}^\alpha(t)
  = (-1)^m \frac{2\Gamma(m+\alpha+1)}{m!\Gamma(\alpha)}
    \sum_{k=0}^m \frac{(-m)_k(m+\alpha+1)_k}{(3/2)_k k!}t^{2k+1}.
\]
We integrate this against $-t(1 - t^2)^\alpha$ using the identity $\int_0^1 s^{k+1/2}(1-s)^\alpha\,\dd s
= \frac{\sqrt{\pi}\Gamma(\alpha+1)}{2\Gamma(\alpha+5/2)}\frac{(3/2)_k}{(\alpha+5/2)_k}$, which
cancels $(3/2)_k$ and leaves the terminating series ${}_2F_1(a, b; c; 1)$ with
parameters
\[
  a = -m, \qquad b = m + \alpha + 1, \qquad c = \alpha + \tfrac{5}{2},
  \qquad\text{so that}\qquad c - b = \tfrac{3}{2} - m.
\]
By the Chu--Vandermonde identity~\cite[Corollary~2.2.3]{AndrewsAskeyRoy1999},
${}_2F_1(-m, b; c; 1) = (c - b)_m/(c)_m$, and here $(\frac{3}{2} - m)_m =
\Gamma(\frac{3}{2})/\Gamma(\frac{3}{2} - m) = (-1)^{m+1}\Gamma(m -
\frac{1}{2})/(2\sqrt{\pi})$ by the reflection formula, while $(c)_m =
\Gamma(m + \alpha + \frac{5}{2})/\Gamma(\alpha + \frac{5}{2})$. Hence
\[
  {}_2F_1\bigl(-m, m+\alpha+1; \alpha + \tfrac{5}{2}; 1\bigr)
  = (-1)^{m+1}\frac{\Gamma(m - \frac{1}{2})\,\Gamma(\alpha + \frac{5}{2})}
                  {2\sqrt{\pi}\,\Gamma(m + \alpha + \frac{5}{2})}.
\]
Applying Legendre's
duplication formula to $C_{2m+1}^\alpha(1) =
\Gamma(2m+2\alpha+1)/[\Gamma(2m+2)\Gamma(2\alpha)]$
yields~\eqref{eq:funk-hecke-eigenvalues} for $\dimin \ge 3$, and letting $\alpha
\to 0^+$ ($\Gamma(2\alpha)/\Gamma(\alpha) \to \frac{1}{2}$,
$C_\ell^\alpha(t)/C_\ell^\alpha(1) \to T_\ell(t)$) gives $\kappa_{2m+1}(2) =
\lambda_m/2$.

\ppart{(ii) and (iii)}
Stirling's formula gives $\kappa_{2m+1}(\dimin) = \Theta(m^{-\dimin})$ and
$N(\dimin, 2m+1) = \Theta(m^{\dimin - 2})$ as $m \to \infty$, so $\sum_{m=0}^\infty
N(\dimin, 2m+1)\abs{\kappa_{2m+1}(\dimin)} \asymp \sum_{m=1}^\infty m^{-2} <
\infty$ and $\tr(T_{\Zlim}) = \frac{\dimin^2\tau_2(P_w)}{4}K_{\mathrm{sph}}(1) =
0$. For part~\textup{(iii)}, Proposition~\ref{prop:degenerate-ustat-limit} gives
the limit of $U_{\ngen,\infty}(W)$, and
Theorems~\ref{thm:second-chaos-transfer} ($\dimin \ge 3$)
and~\ref{thm:planar-universality} ($\dimin = 2$) transfer it to $\Regge(B(0, R))$.
\end{proof}

\begin{proposition}[Rubin and Vitale~{\cite{RubinVitale1980}}, Serfling~{\cite{Serfling1980}}]
\label{prop:degenerate-ustat-limit}
Let $Z_1, \dots, Z_{\ngen}$ be i.i.d.\ in a Polish space $(\mathcal{Z},
\mathcal{B}(\mathcal{Z}), P_Z)$, and let $K\colon \mathcal{Z} \times \mathcal{Z}
\to \RR$ be a symmetric, square-integrable, completely degenerate kernel
($\EE[K(Z_1, Z_2)^2] < \infty$, $\EE[K(z, Z_2)] = 0$ $P_Z$-a.e.) whose
Hilbert--Schmidt operator $(T_K f)(z) := \int_{\mathcal{Z}} K(z, z') f(z')\, \dd
P_Z(z')$ on $L^2(\mathcal{Z}, P_Z)$ has non-zero eigenvalues $(\mu_m)_{m \ge 1}$
(counted with multiplicity). Then $\EE[U_{\ngen}(K)] = 0$, $\Var(U_{\ngen}(K)) =
\frac{\ngen - 1}{2\ngen}\sum_{m=1}^\infty \mu_m^2$ and
\[
  U_{\ngen}(K) := \frac{1}{\ngen}\sum_{1 \le c < c' \le \ngen} K(Z_c, Z_{c'})
  \xrightarrow[\ngen \to \infty]{d}
  U_\infty(K) := \frac{1}{2}\sum_{m=1}^\infty \mu_m (\gamma_m^2 - 1),
  \qquad
  \gamma_m \stackrel{\mathrm{i.i.d.}}{\sim} \Norm(0, 1).
\]
If every non-zero eigenvalue of $T_K$ has multiplicity $2$ with pair values
$(\lambda_m)_{m \ge 0}$, then
\begin{equation}\label{eq:ustat-exp-pairing}
  U_\infty(K) \stackrel{d}{=} \sum_{m=0}^\infty \lambda_m (E_m - 1),
  \qquad
  E_m \stackrel{\mathrm{i.i.d.}}{\sim} \mathrm{Exp}(1),
  \qquad
  \Var(U_\infty(K)) = \sum_{m=0}^\infty \lambda_m^2.
\end{equation}
\end{proposition}

\section*{Declaration of generative AI use}
During the preparation of this manuscript, large language models (Google Gemini
and Anthropic Claude) were used to assist with drafting and editing prose,
drafting and condensing proof write-ups from the author's derivations,
structural organisation and internal error-checking passes. Every definition,
theorem statement and proof, together with all calculations, was independently
verified line by line by the author, who takes full responsibility for the
content of this article.


\end{document}